

\magnification =\magstep1
\baselineskip =14pt

\centerline {\bf CHOW QUOTIENTS OF GRASSMANNIANS I.}

\vskip 1cm

\centerline {\bf M.M.Kapranov}

\vskip 1.5cm

\hfill\hfill {\sl To Izrail Moiseevich Gelfand}

\vskip 1.5cm

\centerline {\bf Contents:}
\vskip .5cm

\noindent {\bf Introduction.}
\vskip .3cm

\noindent {\bf Chapter 0. Chow quotients.}

\noindent 0.1. Chow varieties and Chow quotients.

\noindent 0.2. Torus action on a projective space: secondary polytopes.

\noindent 0.3.  Structure of cycles from the Chow quotient.

\noindent 0.4. Relation to Mumford quotients.

\noindent 0.5. Hilbert quotients (the Bialynicki-Birula-Sommese
construction).

\vskip .3cm

\noindent {\bf Chapter 1. Generalized Lie complexes.}

\noindent 1.1. Lie complexes and the Chow quotient of Grassmannian.

\noindent 1.2. Chow strata and matroid decompositions of hypersimplex.

\noindent 1.3. An example: matroid decompositions of the
 hypersimplex $\Delta (2,n)$.

\noindent 1.4. Relation to the secondary variety for the
 product of two simplices.

\noindent 1.5. Relation to Hilbert quotients.

\noindent 1.6. The (hyper-) simplicial structure on the collection of
$G(k,n)//H$.

\vskip .3cm

\noindent {\bf Chapter 2. Projective configurations and the
 Gelfand- MacPherson isomorphism.}

\noindent 2.1. Projective configurations and their Chow quotient.

\noindent 2.2. The Gelfand-MacPherson correspondence.

\noindent 2.3. Duality (or association).

\noindent 2.4. Gelfand-MacPherson correspondence and Mumford quotients.

\vskip .3cm

\noindent {\bf Chapter 3. Visible contours of (generalized)
Lie complexes and Veronese varieties. }

\noindent 3.1. Visible contours and the logarithmic Gauss map.

\noindent 3.2. Bundles of logarithmic forms on $P^{k-1}$ and
visible contours.

\noindent 3.3. Visible contours as the Veronese variety in the Grassmannian.

\noindent 3.4. Properties of special Veronese varieties.

\noindent 3.5. Steiner construction of Veronese varieties
in the Grassmannian.

\noindent 3.6. The sweep of a special Veronese variety.

\noindent 3.7. An example: Visible contours and sweeps of
 Lie complexes in $G(3,6)$.

\noindent 3.8. Chordal varieties of Veronese curves.

\noindent 3.9. The homology class of a Veronese variety in
Grassmannian and dimensions of representations of $GL(n-k)$.

\vskip .3cm

\noindent {\bf Chapter 4. The Chow quotient of $G(2,n)$ and
 the Grothendieck - Knudsen moduli space $\overline{M_{0,n}}$.}

\noindent 4.1. The space $G(2,n)//H$ and stable curves.

\noindent 4.2. The birational maps $\sigma_i :
 \overline{M_{0,n}}\rightarrow P^{n-3}$.

\noindent 4.3. Representation of the space
 $\overline{M_{0,n}} = G(2,n)//H$ as a blow-up.

\vskip .3cm

\noindent {\bf References.}
\vskip 1cm

\centerline {\bf Introduction.}
\vskip 1cm

The study of the action of the maximal torus $H\i GL(n)$
on the Grassmann variety $G(k,n)$
is connected with numerous questions of geometry and analysis.
Among these let us mention the general theory of hypergeometric
 functions [18],
K-theory [7], combinatorial constructions of characteristic classes
[17,20,40]. It was noted by I.M.Gelfand and R.W.MacPherson [20,40]
 that this problem is equivalent to the classical problem of study
 the projective equivalence classes of configurations of $n$ points
 in the projective space $P^{k-1}$.

In the present paper we propose a geometric approach to this problem
 which is based on the study of the behavior of orbit closures.
 Namely, the closures of  generic orbits are compact varieties
"of the same type". Now there is a beautiful construction in
 algebraic geometry --- that of Chow varieties [50].
It produces compact varieties whose points parametrize algebraic
cycles (= positive integral combinations of irreducible subvarieties)
in a given variety with given dimension and degree. In particular,
any one-parameter family of subvarieties "of the same type"  has a
limit in the Chow variety. We define the {\it Chow quotient}
$G(k,n)//H$ to be the space of such limits of closures of generic orbits.
 Any point of $G(k,n)//H$ represents an $(n-1)$ -dimensional family
 of $(k-1)$ -dimensional projective subspaces in $P^{n-1}$.

Families of subvarieties in a variety $X$ whose parameter space
 has the same dimension as $X$ are classically known as {\it complexes}.
 We call closures of generic orbits in $G(k,n)$ Lie complexes and
 their limit positions- generalized Lie complexes. In the simplest
 case of $G(2,4)$ Lie complexes are the so-called tetrahedral
complexes of lines in $P^3$ which have a long history
 (see the bibliography in [20]).
\vskip .2cm

 The variety  $G(k,n)//H$ can be defined in two more ways:

a) As the space of limits of closures of generic orbits of
the group $GL(k)$ in the Cartesian power $(P^{k-1})^n$ (Theorem 2.2.4).

b) As the space of limits of special Veronese varieties in
the Grassmannian
$G(k-1,n-1) = G(k-1, {\bf h})$, where {\bf h} is the Lie
algebra of the torus $H$(Theorem 3.3.14).
\vskip .2cm

According to interpretation a), the space $G(k,n)//H$ is
obtained by adding some "ideal" elements to the space of
projective equivalence classes of configurations of $n$ points
(or hyperplanes) in $P^{k-1}$ in general position. These ideal
 elements are more subtle than just non-general configurations:
a limit position of closures of generic orbits can be the
union of several orbit closures, each of which represents
 some configuration.

In fact, it turns
out that these  elements behave in many respects as if they
actually were
configurations in general position. In particular, there are
restriction and projection maps (Theorem 1.6.6)
$$G(k,n-1)//({\bf C}^*)^{n-1}\buildrel \tilde b_i\over
\longleftarrow G(k,n)//({\bf C}^*)^n \buildrel
\tilde a_i\over\longrightarrow G(k-1, n-1)//({\bf C}^*)^{n-1}$$
The map $\tilde a_i$ corresponds to the restriction of the generic
 hyperplane
configuration \hfill\break
$(M_1, ..., M_n)$ to the hyperplane $M_i$. The map $\tilde b_i$
 corresponds to deleting $M_i$. The maps $\tilde a_i, \tilde b_i$
 (regarded on the
 generic part of Grassmannian)
were at the origin of the use of Grassmannians in the problem of
combinatorial calculation of Pontryagin classes [20,40]. These
maps were also used in [7] to define
the so-called Grassmannian complexes serving as an approximation
 to K-theory.
Note that $\tilde a_i$ cannot be, in general, extended to
arbitrary configurations. It is the Chow quotient approach
 that permits this extension.
\vskip .2cm

Veronese varieties (which make their appearance in the
interpretation b) above) are defined classically as images
 of a projective space in the projective embedding given by
all homogeneous polynomials of given degree. In particular,
 Veronese curves are just rational normal curves and  posess
 a lot of remarkable geometric properties,see [46].
In our situation Veronese varieties in fact lie on a
 Grassmannian (in Pl\"ucker embedding): a $(k-1)$ -
dimensional variety lies in $G(k-1, {\bf h})$ so that for
$k=2$ we obtain the curve in a projective space.

The Veronese variety associated to a Lie complex $Z$ is
obtained as its {\it visible contour} i.e. the locus of
subspaces from $Z$ which contain a given generic point $p$,
say $p=(1:...:1)\in P^{n-1}$. The consideration of visible
contours is a classical method of analyzing  complexes of
subspaces [24,27] which, to our knowledge, has yet never been
 applied to Lie complexes.

Along with the visible contour of $Z$ lying in $G(k-1,{\bf h})$
 we consider the so-called {\it visible sweep} $Sw(Z)$. This is
a subvariety in the projective space $P({\bf h})$ which is the
union of all subspaces from the visible contour. This variety
 can be found very explicitly. Thus, if our  Lie complex has
the form $Z=\overline {H.L}$ where $L\in G(k,n)$ is the graph
 of a linear operator $A=||a_{ij}||:{\bf C}^k\rightarrow
{\bf C}^{n-k}, \quad k\leq n-k$ (this assumption does not
restrict the generality) then the sweep is the projectivization
 of the following determinantal cone:
$$\{(t_1,...,t_n)\in {\bf h} :\quad {\rm rank} \,\,
||a_{ij}(t_i-t_j)||_{i=1,...,k, j=1,...,n-k} \leq k\}.$$
The linear forms $t_i-t_j$ entering the matrix above
 are  roots of {\bf h}, the Cartan subalgebra of ${\bf pgl}(n)$.
 More precisely, we encounter those roots  which enter into the
 weight decomposition of the parabolic subalgebra defining the Grassmannian.

Veronese varieties in Grassmannian which arise naturally in our
 constructions, seem th be "right" generalizations of Veronese
curves in projective spaces. We show in \S 3.5 that these varieties
 admit a Steiner-type construction, which is well known for curves [24].

The homology class in the Grassmannian of such a variety is given
by an extremely beautiful formula (Theorem 3.9.8). To state it,
 recall that the homology of Grassmannian is freely generated
 by  Schubert cycle $\sigma_\alpha$ which correspond to Young
diagrams. It turns out that the multiplicity of the cycle
$\sigma_\alpha$ in the class of Veronese variety equals the
 dimension of the space $\Sigma^{\alpha^*}({\bf C}^{n-k})$,
 the irreducible representation of the group $GL(n-k)$
corresponding to the Young diagram dual to $\alpha$.

In the particular case $k=2$   the construction b) realizes
points of our variety  as limit positions of rational normal
 curves (Veronese curves, for short) in $P^{n-2}$ through a
fixed set of $n$ points in general position. We deduce from
that  that $G(2,n)//H$ is isomorphic to the Grothendieck-Knudsen
 moduli space $\overline{M_{0,n}}$ of stable $n$ -pointed curves
 of genus 0 (Theorem 4.1.8). This is certainly the most natural
compactification of the space $M_{0,n}$ of projective equivalence
 classes of $n$ -tuples of distinct points on $P^1$. It is smooth
 and the complement to $M_{0,n}$ is a divisor with normal crossings.

In general, Veronese varieties in Grassmannians arising in construction b)
are Grassmannian embeddings of $P^{k-1}$ corresponding to
 vector bundles on  of the form $\Omega^1(\log M)$, where
 $M=(M_1,...,M_n)$ is a configuration of hyperplanes in
 general position. They become Veronese varieties after
the Pl\"ucker embedding of the Grassmannian. The configuration
of hyperplanes on $P^{k-1}$ can be read off the corresponding
Veronese variety by intersecting it with natural sub-Grassmannians.
 This is explained in \S 3.
\vskip .2cm

 The Chow quotients of toric varieties by the action of a subtorus of
 the defining torus
were studied in [30,31]. It was found that this provides a natural
setting for the theory of secondary polytopes introduced in [22,23] and
their generalizations - fiber polytopes [8]. By considering the
Pl\"ucker embedding of Grassmannian we can apply the results of [30,31].
 This gives a description of possible degenerations of orbit closures
in $G(k,n)$ (i.e. of generalized Lie complexes) in terms of polyhedral
decomposition of a certain polytope $\Delta(k,n)$ called the hypersimplex
 [17,19,20]. For the case
of $G(2,n)$ these decompositions are in bijection with  trees
which describe combinatorics of stable curves.

It is now well-known [19,21,40] that various types of closures
 of torus orbits in $G(k,n)$ correspond to {\it matroids} i.e.
 types of combinatorial behavior of a configurations of $n$
hyperplanes in projective space. To each matroid the corresponding
{\it stratum} in $G(k,n)$ is associated. It is  formed by orbits of
 the given type [18,19,21]. Matroids are in one-to-one
correspondence with certain polytopes in the hypersimpex
[19,21]. From our point of view, however, a natural object
is not an individual matroid but a collection of matroids
such that the corresponding polytopes form a polyhedral d
ecomposition of hypersimplex. We call such collection
{\it matroid decompositions}.

\vskip .2cm

It should be said that our approach differs considerably
 from that of geometric invariant theory developed by
 D.Mumford [42]. In particular, Mumford's quotients of
$G(2,n)$ by $H$ and of $(P^1)^n$ by $GL(2)$ though isomorphic
 to each other, do not coincide with the space $\overline{M_{0,n}}$
 which provides a finer compactification. Note that Mumford's
quotient depends upon a choice of a projective embedding
(and this is felt for varieties like $(P^{k-1})^n$ with
large Picard group) and upon the choice of linearization
i.e. the extension of the action to the graded coordinate
ring of the embedding
(and this is felt for groups like the torus).
We prove in \S 0 that the Chow quotient always maps to any
 Mumford quotient by a regular birational map.

Instead of Chow variety one can use Hilbert schemes and obtain
a different compactification. Such a construction was considered
in 1985 by A.Bialynicki-Birula and A.J.Sommese [10] and later by
 Y.Hu [26]. The advantage of Hilbert schemes is that they represent
 an easily described functor so it is easy to construct morphisms
into them. In our particular example of the torus action on $G(k,n)$
 both constructions lead to the same answer.

\vskip .2cm

In the forthcoming second part of this paper we shall study the
 degenerations of Veronese varieties which provide a
higher-dimensional analog of stable curves of \hfill\break
Grothendieck and Knudsen. The main idea of stable pointed
curves is that points are never allowed to coincide. When
they try to do so, the topology of the curve changes in such
 a way that the points remain distinct. In higher dimensions
instead of a collection of points we have a divisor on a variety.
 The analog of the condition that points are distinct is that the
 divisor has normal crossings.
In particular, a generic configuration of hyperplanes defines
a divisor with normal crossings. When the hyperplanes "try" to
 intersect non-normally, the  corresponding Veronese variety
 degenerates in such a way as to preserve the normal crossing.

I am grateful to W.Fulton who suggested that the space
$\overline{M_{0,n}}$ might be related to Chow quotients
 and informed me on his joint work  with R.MacPherson [15]
 on a related subject. I am also grateful to Y.Hu for informing
 me about his work [26] and  about earlier work of
A.Bialynicki-Birula and A.J.Sommese [10].

\vskip .3cm

I am happy to be able to dedicate this paper to Izrail Moiseevich Gelfand.

\hfill\vfill\eject

\beginsection Chapter 0. CHOW QUOTIENTS.

\vskip 1cm

\centerline {\bf (0.1) Chow varieties and Chow quotients.}
\vskip .5cm

Let $H$ be an algebraic group acting on a complex projective variety $X$.
We shall describe in this section an approach to constructing of
 the algebraic "coset space" of $X$ by $H$ which was introduced
in [31]. (A similar approach was  introduced earlier in [10],
 see \S 0.5 below).
\vskip .3cm

 \noindent (0.1.1)\underbar {Setup of the approach.} For any
point $x\in X$ we consider the orbit closure $\overline {H.x}$
 which is a compact subvariety in $X$. For some sufficiently small
 Zariski open subset $U\i X$ of "generic"  points all these
varieties have the same dimension, say, $r$ and  represent
the same homology class $\delta\in H_{2r}(X,{\bf Z})$.
 The set $U$ may be supposed $H$- invariant. Moreover,
since we are free to delete bad orbits from $U$, the
construction of the quotient $U/H$ presents no difficulty
and the problem is to construct a "right" compactification
 of $U/H$. A natural approach to this is to study the limit
 positions of the varieties $\overline {H.x}$
when $x$ tends to the infinity of $U$ (i.e. ceases to be generic).
One of the precise ways to speak about such "limits" is provided by
  Chow varieties of algebraic cycles. Before proceeding further we
 recall the main definitions (cf. [30,50]).
\vskip .3cm

\noindent (0.1.2) By a (positive) $r$- dimensional algebraic
cycle on $X$ we shall understand  a finite formal non-negative
 integral combination $Z=\sum c_i Z_i$ where $c_i\in {\bf Z}_+$
 and $Z_i$ are irreducible $r$-dimensional  closed algebraic
subvarieties in $X$. Denote by ${\cal C}_r(X,\delta)$ the set
of all $r$- dimensional algebraic cycles in $X$ which have
 the homology class $\delta$. It is known   that
${\cal C}_r(X,\delta)$ is canonically equipped with a
 structure of a \underbar {projective} (in particular,
 compact) algebraic variety (called the Chow variety).
In this form this result is due to D.Barlet [3].
\vskip .3cm

\noindent (0.1.3) A more classical approach to Chow varieties is that
of Chow forms [50].
This approach first gives the projective embedding of
${\cal C}_r(P(V),d)$, the Chow variety of $r$-dimensional
 cycles of degree $d$ in the projective space
$P(V)$. The Chow form of any cycle $Z, dim(Z)=r, deg(Z)=d$,
is a polynomial
$R_Z(l_0,...,l_r)$ in the coefficients of $r+1$ indeterminate
linear
forms $l_i\in V^*$ which is defined, up to a constant factor,
 by the following
 properties,see [30,50]:

\item{(0.1.3.1)} $R_{Z+W}=R_Z.R_W$.
\item{(0.1.3.2)} If $Z$ is an irreducible subvariety then $R_Z$ is an
irreducible polynomial which vanishes for given $l_0,...,l_r$ if
and only if the projective subspace $\Pi(l_0,...,l_r)=\{l_0=...=l_r=0\}$
of codimension $r+1$ intersects $Z$.
\vskip .3cm

\noindent (0.1.4) It is now classical [50] that the correspondence
$Z\mapsto R_Z$ identifies
${\cal C}_r(P(V),d)$ with a Zariski closed subset of the projective
 space of
polynomials $F(l_0,...,l_r)$ homogeneous of degree $d$ in each $l_i$.
\vskip .3cm

\noindent (0.1.5) If $X\i P(V)$ is any projective subvariety and
$\delta\in H_{2r}(X,{\bf Z})$, then
${\cal C}_r(X,\delta)$ becomes a  subset of ${\cal C}_r(P(V),d)$
 where $d\in H_{2r}(P(V), {\bf Z})={\bf Z}$ is the image of $\delta$.
The   result of Barlet mentioned in (0.1.2) shows,
in particular, that (over a field of complex numbers) this subset
 is Zariski
closed and the resulting structure of algebraic variety on
${\cal C}_r(X,\delta)$ does not depend on the projective embedding.
(The fact that for a given projective subvariety $X$ the set of
 $Z\in {\cal C}_r(P(V),d)$ lying on $X$ is Zariski
closed, is classical).

In the case of base field of characteristic $p$, which we do not
consider here,
the situation is more subtle, see [43,44].
\vskip .3cm

\noindent (0.1.6) Let us return to the situation (0.1.1) of the
group $H$ acting on $X$.
 We see that for  $x\in U$ as in (0.1.1) the subvariety $\overline {H.x}$
is a point of the variety ${\cal C}_r(X,\delta)$. The correspondence
$x\mapsto \overline {H.x}$ defines therefore an embedding of the
quotient variety $U/H$ into ${\cal C}_r(X,\delta)$.

\proclaim (0.1.7) Definition. {\rm [31]} The Chow quotient $X//H$ is the
 closure of $U/H$ in ${\cal C}_r(X,\delta)$.

Thus $X//H$ is a projective algebraic variety compactifying $U/H$.
"Infinite" points of $X//H$ are some algebraic cycles in $X$ which
are limits (or "degenerations") of generic orbit closures.

\noindent {\bf (0.1.8) Remarks.} a) Definition (0.1.7) does not depend
 on the freedom in the choice of $U$ since deletion from $U$ of orbits
which are already "generic" results in their reappearance as points in
the closure.

b) The notion of "genericity" used in  Definition (0.1.7) is usually much
 more restrictive than Mumford's notion of  stability [42]. In fact (0.1.7)
 makes no appeal to stability and is defined entirely in terms of $X$ and
 the action of $H$.

\vskip 1cm

\centerline {\bf (0.2)  Torus action on a projective space:
 secondary polytopes.}
\vskip .5cm

If $H$ is an algebraic torus acting on a projective variety
 $X$, then $X$ may be equivariantly embedded into a projective
 space with $H$ -action.
 The case of torus action on a projective space recalled in
 this subsection
 will be basic for our study of more general torus actions in this paper.
\vskip .3cm

\noindent (0.2.1) Let  $H$ be an algebraic torus $({\bf C}^*)^k$.
 A character of $H$ is the same as a Laurent monomial
$t^\omega= t_1^{\omega_1}...t_k^{\omega_k}$ where $\omega
= (\omega^{(1)},...,\omega^{(k)})\in {\bf Z}^k$ is an
integer vector. A collection $A= \{\omega_1,..., \omega_N\}$
of vectors from ${\bf Z}^k$ defines therefore a diagonal
homomorphism from $H$ to $GL(N)$.
It is well-known that any representation of the torus can
 be brought
into a diagonal form.
\vskip .3cm

\noindent (0.2.2) The homomorphism $H \rightarrow GL(N)$ constructed
from the set $A$ above
defines an $H$- action on the projective space
 $P^{N-1}$. The homogeneous coordinates in $P^{N-1}$ are naturally
labelled by elements of $A$. So we shall denote this space by $P(A)$
 and the coordinates by $(x_\omega)_{\omega\in A}$ thus dropping
the (unnatural) numeration of $\omega$ 's.

 The Chow quotient $P(A)//H$ was described in [30,31].
 We shall use this description so we recall it here.

\vskip .3cm

\noindent (0.2.3)
First of all, $P(A)//H$ is a projective toric variety.

To see this, we note that the
"big" torus $({\bf C}^*)^A$  acts on $P(A)$
(by dilation of homogeneous coordinates)  commuting with
 $H$ (which is just a subtorus of $({\bf C}^*)^A$).
Therefore $P(A)//H$ is the closure, in the Chow variety,
of the $({\bf C}^*)^A$ - orbit of the variety $X_A =\overline{H.x}$,
 where $x\in P(A)$ is the point with all coordinates equal to 1.
\vskip .3cm

\noindent (0.2.5) It is known that projective toric
varieties are classified by lattice polytopes, see [49].
 In what follows we shall describe the polytope corresponding
 to the toric variety $P(A)//H$.

\vskip .3cm
\noindent (0.2.6) Let $Q\i {\bf R}^k$ be the convex hull of the set $A$.
A {\it triangulation} of the pair $(Q,A)$ is a collection of simplices
 in $Q$ whose vertices lie in $A$  intersecting only along common
faces and covering $Q$. To any such triangulation $T$ we associate
its {\it characteristic function} $\phi_T :A\rightarrow {\bf Z}$ as follows.
By definition, the value of $\phi_T$ on $\omega\in A$ is the sum
of volumes of all simplices of
 $Q$ for which $\omega$ is a vertex.
The volume form is normalized by the condition that the smallest possible
volume of a lattice simplex equals 1.
\vskip .3cm

 \noindent (0.2.7) The {\it secondary polytope} $\Sigma (A)$
 is, by definition, the convex hull of all characteristic functions
$\phi_T$ in the space ${\bf R}^A$.

Secondary polytopes were introduced in [22,23] in connection with
Newton polytopes of multi-dimensional discriminants.
It was shown in [22] that the vertices
 of $\Sigma (A)$ are precisely functions $\phi_T$ where the triangulation
$T$ is {\it regular} i.e. posesses a strictly convex piecewise- linear
function.

\vskip .3cm
\proclaim (0.2.8) Theorem. {\rm [30]} The toric variety
 $P(A)//H$ corresponds to the convex lattice polytope $\Sigma (A)$.

\vskip .3cm

\noindent (0.2.9) \underbar {Complements.}
All the faces of the secondary polytope $\Sigma (A)$ posess a
complete description.
We shall use in this paper only the case when  elements of $A$ are
exactly vertices of $Q$ so we shall restrict ourselves to this case,
see [22] for general case.
Let us call a {\it polyhedral decomposition} of $Q$ a collection of
 convex polytopes in $Q$ whose vertices lie in $A$, which intersect
 only along common faces and cover $Q$. A polyhedral decomposition
 ${\cal D}$ is called {\it regular}, if it posesses a strictly convex
 pieciwise- linear function. It was shown in [22] that vertices of
 $\Sigma (A)$ are in bijection with regular polyhedral decompositions
 of $Q$. Vertices of the face of $\Sigma(A)$ corresponding to such a
 decomposition ${\cal D}$ are precisely functions $\phi_T$ for all
 regular triangulations $T$ refining ${\cal D}$.

\vskip .3cm

\noindent (0.2.10) It was shown  in [31] that any cycle from
$P(A)//H$ is a sum of toric
subvarieties (closures of $H$ -orbits), see [31], Proposition 1.1.
 In particular, a regular triangulation $T$ represents a 0-dimensional
torus orbit in $P(A)//H$ i.e. some algebraic cycle. This cycle has
 the form $\sum_{\sigma\in T} Vol(\sigma).{\bf L}(\sigma)$ where
${\bf L}(\sigma)$ is the coordinate $k$ -dimensional projective
subspace in $P(A)$ spanned by basis vectors corresponding to
vertices of $\sigma$.

\vskip 1cm

\centerline {\bf (0.3) Structure of cycles from Chow quotient.}

\vskip .5cm

\proclaim (0.3.1) Theorem. Let $H$ be a reductive group acting on a smooth
 projective variety $X$. Suppose that the stationary subgroups
$H_x, x\in X$ are trivial for  generic $x$ and are never unipotent.
 Then any component $Z_i$ of any cycle $Z=\sum c_iZ_i\in X//H$ is a
closure of a single $H$- orbit.

\noindent {\sl Proof:} For the case when $H$ is a torus, this statement
 follows from results of [30,31]. Indeed, we can take an equivariant
embedding of $X$ into  a projective space $P^N$ with $H$ -action in such
a way that the dimension of a generic $H$ -orbit on $X$ is the same as
the dimension of a generic $H$ -orbit on $P^N$. The degeneration of torus
 orbits on a projective space (and, more generally, on toric varieties)
 was studied in [30,31] where it was found that any orbit degenerates
in a union of finitely many orbits ([31], Proposition 1.1).

 Consider now the general case. Let $C(t), t\neq 0$, be a 1-parameter
family of closures of generic orbits, $C(0)= \lim_{t\rightarrow 0}C(t)$
 is their limit in the Chow variety. Let $C$ be any component of $C(0)$.
Suppose, contrary to our statement, that $C$ is not a closure of a single
 orbit. Then, for all points of $C$ the stabilizer $H_x$ has positive
dimension. By our assumption, these stabilizers  all non-unipotent.
Let $x$ be some fixed generic point of $C$. Then $H_x$ contains some
 torus $T$. Include $x$ in a 1-parameter family of points $x(t)\in C(t)$
 such that for $t\neq 0$ the point $x(t)$ lies in the  orbit open in $C(t)$.
 Consider the closures of orbits $\overline {T.x(t)}\i C(t)$.
When $t\rightarrow 0$ these closure should degenerate into some cycle
$Z$ whose support contains $x$. But we know that each component of $Z$
is a closure of one $T$- orbit and so $x$ should lie on the intersection
of components of $Z$. This means that each generic point $x\in C$ lies
in the closure of some orbit $H.y$ not coinciding with $H.x$. This is
impossible.

\noindent {\bf (0.3.2) Example.} The assumption that the stabilizers of
points are
never unipotent in the formulation of
Theorem 0.3.1 can not be dropped. To construct an example, consider
the group $H=SL(2,{\bf C})$. It has a standard action on ${\bf C}^2$.
Let $X_0$ be the product ${\bf C}^2\times {\bf C}^2$ on which $H$ acts
diagonally. Let $X= P^2\times P^2$ be the natural compactification
of $X_0$ with the obviously extended action of $H=SL(2)$.
Generic $H$ -
orbits on $X$ are 3-dimensional : two pairs of independent vectors
$(e_1, e_2)$ and $(f_1,f_2)$ can be brought to each other by a unique
transformations from $SL(2)$ if and only if $det(e_1,e_2) = det(f_1,f_2)$.
Thus a generic orbit depends on one parameter namely $det (e_1,e_2)$.
However, when this parameter approaches 0,
the orbit degenerates into the 3-dimensional variety  of proportional
 pairs $(e_1,e_2)$. This variety is a union of a 1-parametric family
of 2-dimensional orbits $O_\lambda = \{(e_1,e_2):e_1=\lambda e_2\}$.
The stabilizer of each of this orbit
is  a unipotent subgroup in $H=SL(2)$.

\vskip 1cm

\centerline {\bf (0.4) Relation to Mumford quotients.}

\vskip .5cm

It is useful to have a comparison of the Chow quotient with the more standard
constructions, namely Mumford's geometric invariant theory quotients [42].

\vskip .3cm

\noindent (0.4.1) To define
 the Mumford's quotient, we should choose an $H$ -equivariant
projective embedding of $X$ and extend the $H$ -action to the
 homogeneous coordinate ring ${\bf C}[X]$ of $X$ with respect
 to this embedding.
This is equivalent to extending the action to the ample line
 bundle ${\cal L}$ defining the embedding. Such an extension
is called {\it linearization}. Denote the chosen linearization by $\alpha$.
 The Mumford's quotient $(X/H)_\alpha$ or $(X/H)_{{\cal L}, \alpha}$
corresponding to ${\cal L}, \alpha$ is defined as $Proj {\bf C}[X]^H$,
the projective spectrum of the invariant subring [42]. Thus there are
 two choices in the definition
 of Mumford quotient: that of an ample line bundle and that of
 extension of the action to the chosen line bundle.
\vskip .3cm

 \noindent (0.4.2) By general theory of  [42], points of
 $(X/H)_{{\cal L}, \alpha}$ are equivalence classes of
 $\alpha$ -semistable orbits in $X$. More precisely, two
 semistable orbits $O,O'$ are equivalent if any invariant
 homogeneous function vanishing on $O$ does so
 on $O'$ and conversely.
A Zariski open set in $(X/H)_{{\cal L}, \alpha}$ is formed
by $\alpha$ -stable orbits [42].
The have the property that no two $\alpha$ -stable orbits are
 equivalent. We shall say that the linearization $\alpha$ is
 non-degenerate if there are $\alpha$ -stable orbits.

The following result was proven in [31] for the case of a torus acting
on a toric
variety and, (independently and simultaneously) in [26] for
 torus action on an arbitrary variety.

\proclaim (0.4.3) Theorem. Let $H$ be a reductive group
 acting on a projective
variety $X$, ${\cal L}$ - an ample line bundle on $X$
and $\alpha$ be a
 linearization i.e. an extension of the $H$ -action on $X$
to ${\cal L}$.
 Suppose that $\alpha$ is non-degenerate.  Then there is
 a regular birational morphism $p_\alpha:X//H\rightarrow
(X/H)_{{\cal L}, \alpha}$.

 For any algebraic cycle $Z=\sum n_iZ_i$ in $X$ we shall call its support
 and denote $supp (Z)$ the union $\cup Z_i$.
The proof of Theorem 0.4.3 consists of three steps:
\vskip .3cm

\noindent{\bf (0.4.4)} For any cycle $Z\in X//H$ as
above there is at least one orbit
 in $supp (Z)$ which is $\alpha$ - semistable.
\vskip .3cm

\noindent{\bf (0.4.5)} All the $\alpha$ - semistable
orbits in $supp (Z)$ are equivalent
i.e. represent the same point of the Mumford quotient
 $(X/H)_{{\cal L}, \alpha}$.
\vskip .3cm

\noindent{\bf (0.4.6)} The map $p_\alpha:X//H\rightarrow
(X/H)_{{\cal L}, \alpha}$ which takes
$Z\in X//H$ to the point of $(X/H)_{{\cal L}, \alpha}$
represented by any of
the semistable orbit in $supp (Z)$, is a morphism of algebraic
varieties.
\vskip .3cm

\noindent (0.4.7) \underbar{ Proof of {\bf (0.4.4)}}:
 We use the interpretation of semistability via the moment map
[1,35]. Let $H_c$ be the compact real form of $H$ with Lie algebra
${\cal H}$ and $\mu:X\rightarrow {\cal H}^*$ be the moment map
associated to an $H_c$ -
invariant K\"ahler form on $X$. Then an orbit $O$ is semistable
 if and only if
$\mu(\bar O)$ contains zero element of ${\cal H}^*$. Let $O(t),
 t\neq 0$, be a 1-parameter family of generic orbits and $Z(t)$
 be the closure of $O(t)$. Let $Z(0)$ be the limit of $Z(t)$ is
 the Chow variety. Since $\mu(Z(t))$ is, for $t\neq 0$, a closed
set containing 0, the set $\mu(Z(0))$ also contains 0
thus proving that at least one orbit constituting $Z(0)$, is semistable.

\vskip .3cm

\noindent (0.4.8) \underbar { Proof of {\bf (0.4.5)}}: Denote by
 ${\cal L}$ the equivariant
 ample invertible
sheaf given by the linearization $\alpha$.
By definition, two semistable orbits $O$ and $O'$ are equivalent in
 $(X/H)_{{\cal L},\alpha}$ if any section of ${\cal L}^{\otimes k}$
 vanishing
at $O$, does so at $O'$.
But the cycle $Z$ is a limit position of closures of single orbits.
So our assertion follows by continuity.
\vskip .3cm

\noindent (0.4.9) \underbar { Proof of {\bf (0.4.6)}}:
Let $X\i P(V)$ be the equivariant
projective embedding given by the linearization. We shall use
the  approach to Chow variety via Chow forms, see (0.1.3).
Let $d$ be the degree of a generic orbit closure
 $\overline {H.x}, x\in X$.
 Recall that the Chow form of any cycle $Z, dim(Z)=r, deg(Z)=d$,
 is a polynomial
$R_Z(l_0,...,l_r)$ in the coefficients of $r+1$ indeterminate linear
forms $l_i\in V^*$.
\vskip .3cm

  \noindent (0.4.9.1) Since the property of being a morphism
is local, it suffices to prove it in a suitable open covering of $X//H$.
More precisely, we are reduced to the following situation.
\vskip .3cm

\noindent (0.4.9.2) Let $f$ be an
invariant rational function on $V$ homogeneous of degree 0 (so it represents
 a regular
 function on some open set of $(X/H)_{{\cal L}, \alpha}$).
 We must express the (constant)
 value of $f$
on a generic orbit $O$ as a rational function of the coefficients of the
Chow form $R_Z$, where $Z=\overline{O}$.

\vskip .3cm

\noindent (0.4.9.3) We can write
 (for characteristic 0 only!)
$$f|_Z = {1\over d} \sum_{{\bf x}\in Z\cap L} f({\bf x}),$$
where $L$ is a generic projective subspace in $P(V)$ of codimension $r$.
On the other hand, let $l_1,...,l_r$ be equations of $L$. Then we have
the equality of polynomials in $l\in V^*$:
$$ R_Z(l,l_1,...,l_r) = c.\prod_{{\bf x}\in Z\cap L}l({\bf x})
\leqno (0.4.9.4)$$
 where $c$ is a non-zero number depending on $l_1,...,l_r$.
\vskip .3cm

\noindent (0.4.9.4) Let $V={\bf C}^{N+1}$ with coordinates $x_0,...,x_N$
 and let $\xi_i$ be dual
coordinates in $V^*$, so an indeterminate linear form on $V$
is $ (\xi,{\bf x})=\sum \xi_ix_i$ for some $\xi_0,...,\xi_N$.
The Chow form of any 0-cycle $W=  {\bf x}^{(1)}+...+{\bf x}^{(d)}$ is
the polynomial $\prod ({\bf x}^{(i)},\xi)$. Let us restrict
the considerations to the affine chart, say,  ${\bf C}^N=\{x_0\neq 0\}$
in $P(V)=P^N$.
The coordinate $x_0$ can then set to be 1 and we can set $\xi_0=1$ as well
thus obtaining  the Chow form of a 0-cycle $W\i A^N$ as before in the form
$$\Phi_W (\xi) =\prod (1+x^{(i)}_1\xi_1+...+x^{(i)}_N\xi_N).$$
 The coefficients of this polynomial at various monomials in
$\xi$'s are known an elementary symmetric functions in $d$ v
ector variables ${\bf x}^{(1)},...,{\bf x}^{(d)}$, see [28,39].
By formula (0.1)  elementary symmetric functions of the $d$ points of
intersection $Z\cap L$, for any generic $L$ of codimension $r$, can be
 polynomially expressed through the coefficients of $R_Z$.
. Therefore we are reduced to the following lemma.

\proclaim (0.4.9.5) Lemma. Let $f({\bf x}), {\bf x}=(x_1,...,x_N)$ be a
rational function in $N$ variables and $d>0$. Then there is a
rational function $U_f = U_f(\Phi)$ in the coefficients of an
indeterminate homogeneous polynomial $\Phi(\xi_1,...,\xi_N), deg(\Phi)=d$
 satisfying the following property. If ${\bf x}^{(1)},...,{\bf x}^{(d)}$
 are points not lying on the polar locus of $f$ then
$$\sum f({\bf x}^{(i)}) = U_f(\Phi_{W})$$
where $\Phi_W = \prod (1+({\bf x}^{(i)},\xi))$.

\noindent {\sl Proof:} It is known since P.A.MacMahon [28,39] that any
 symmetric polynomial in ${\bf x}^{(i)}$ (in characteristic 0) can be
 polynomially expressed via elementary symmetric polynomials (in many
different ways, if $N>1$). If $f({\bf x})=P({\bf x})/Q({\bf x})$,
where $P,Q$ are relatively prime polynomials, then
$$\sum f({\bf x}^{(i)})={1\over\prod Q({\bf x}^{(i)})}
\sum_i P({\bf x}^{(i)})\prod_{j\neq i}Q({\bf x}^{(j)})$$
is a ratio of two symmetric polynomials and the assertion follows.

The proof of Theorem 0.4.3 is finished.
\vskip .3cm

{\noindent {\bf (0.4.10) Remark.} In [31] it was shown that for the case
 of torus action on a toric variety
the Chow quotient $X//H$ is, in some sense, the "least common multiple" of
all Mumford's quotients corresponding to different linearizations. From  the
point of view of general reductive groups a more typical  case is   when the
 group
has 0-dimensional center and hence there is only one linearization.
 However, we shall see that
 in this case
Chow quotient still differs drastically from the Mumford one.
The reason for
this, as we would like to suggest, is that the Chow quotient
takes into account  not only  Mumford quotients corresponding
to various linearizations,
 but also more general
symplectic quotients [1] corresponding to coadjoint orbits of $H_c$.
 For a torus, a different choice of a coadjoint orbit amounts to a change
of a linearization, see [31] so all the symplectic quotients are reduced
to Mumford's ones. In general case the symplectic quotients corresponding
 to non-zero orbits
may not have an immediate algebro-geometric interpretation [1]. Nevertheless,
their presense is somehow felt in $X//H$.

\vskip 1cm

\centerline {\bf (0.5) Hilbert quotients (the
 Bialynicki-Birula-Sommese construction).}

\vskip .5cm

A different way of speaking about limit positions of generic orbit closures
is that of Hilbert schemes. Such a construction was considered by
A. Bialynicki - Birula and A.J. Sommese [10] and later by Y.Hu [26].
\vskip .3cm

\noindent (0.5.1) Recall [25,47] that for any projective variety
$X$ there is the {\it Hilbert scheme } ${\cal H} _X$ parametrizing
all subschemes in $X$. By definition, a
 morphism
$S\rightarrow {\cal H}_X$ is a flat family of subschemes in $X$
parametrized by $S$. The scheme ${\cal H}_X$ is of infinite type
since no bound on "degrees" of subschemes is imposed. The connected
 components of ${\cal H}_X$ are, nevertheless, finite-dimensional
 projective schemes.
\vskip .3cm

\noindent (0.5.2) Any connected component of the scheme ${\cal H}_X$ is
 canonically mapped into the  Chow variety "corresponding"
 to this component. More precisely,
if $K$ is any such connected component then dimensions of
 subschemes from
$K$ are the same and equal, say, $r$. For any scheme $Z\in K$
we define
 the algebraic cycle
$$Cyc(Z) = \sum_{C\i supp(Z) \,\, -\, {\rm irred.},\, dim(C)=r}
 Mult_C(Z) .C \leqno (0.5.3)$$
where $C$ runs over all $r$ -dimensional irreducible components of the
algebraic variety $supp(Z)$ and $Mult_CZ$ is the multiplicity given by
the scheme structure , see [30,42].
\vskip .3cm

\noindent (0.5.4) In the situation of (0.5.2) it follows from
results of [42] that
the cycles $Cyc(Z)$ for all subschemes
$Z\i K$ have the same homology class, say $\delta_K$ and the
 formula (0.5.3)
 defines a regular morphism $K\rightarrow {\cal C}_r(X,\delta_K)$, see [42].
\vskip .3cm

\noindent (0.5.5)
Hilbert schemes have the advantage over Chow varieties in that they
are defined
as objects representing an easily described functor
(that of flat families of subschemes, see (0.5.1)).
In particular, the Zariski tangent space to the scheme
 ${\cal H}_X$ at a point given by a subscheme $Z$ equals (see [47], Proposition
8.1):
$$T_Z {\cal H}_X = H^0(Z, {\cal N}_Z),\quad {\rm where}
\quad {\cal N}_Z = \underline{Hom}(J_Z/J_Z^2, {\cal O}_Z).\leqno(0.5.6)$$
Here $J_Z$ is the sheaf of ideals of the subscheme $Z$.
 The  sheaf
${\cal N}_Z$ is called the normal sheaf of $Z$. It  is
locally free if $Z$ is a locally complete intersection.
 If $Z$ is a smooth variety then ${\cal N}_Z$ is the sheaf
of sections of the normal bundle of $Z$.

\vskip .3cm

\noindent (0.5.6) Consider the situation of (0.1.1) i.e.
an action of an algebraic group $H$ on a projective variety $X$.
 Then for a small open $H$ -invariant set $U\i X$
the orbit closures $\overline{H.x}$ form a flat family. We obtain
an embedding $U/H\hookrightarrow {\cal H}_X$.

\proclaim (0.5.7) Definition. The Hilbert quotient $X///H$ is the closure of
$U/H$ in the Hilbert scheme ${\cal H}_X$.

Thus $X///H$ is a projective algebraic variety compactifying $U/H$.
 "Infinite" points of $X///H$ correspond to subschemes in $X$ which
 are "degenerations" of generic orbit closures.

\vskip .3cm

\noindent (0.5.8) The cycle map (0.5.3) provides a canonical
regular birational morphism
$$\pi: X///H \longrightarrow X//H \leqno (0.5.9)$$
from the Hilbert quotient
 to the Chow quotient. This morphism may be very non-trivial
even in the case when the group $H$ is finite. So $X///H$ provides
 a still finer compactification.
\vskip .3cm

\noindent (0.5.10) In general the Hilbert quotient is rather hard to
 describe. For instance,
in the
 case of torus action on the projective space considered in
section (0.2)  the
 Hilbert quotient is the toric variety corresponding to the
so-called {\it state polytope }
of the toric subvariety $X_A$ introduced by D.Bayer, I.Morrison
 and M.Stillman
[5,6]. However, its exact description depends not only on the
geometry
of the set $A$ (as is the case for the secondary polytope) but also on the
arithmetic nature of relation between elements
of $A$.

We shall see later  that for the torus action on the Grassmannian
Hilbert and Chow quotients coincide thus allowing us to use the
 advantages of both approaches.

\hfill\vfill\eject

\beginsection Chapter 1. GENERALIZED LIE COMPLEXES.

\vskip 1cm
\centerline {\bf (1.1) Lie complexes and the Chow quotient of Grassmannian}
\vskip .5cm

\noindent (1.1.1) Let ${\bf C}^n$ be the coordinate $n$- dimensional
complex vector space
with coordinates $x_1,...,x_n$. By $G(k,n)$ we shall denote
the Grassmannian
of $k$- dimensional linear subspaces in ${\bf C}^n$.
The group $({\bf C}^*)^n$
of diagonal matrices acts on $G(k,n)$.
 Since homotheties act trivially we obtain in fact an action of
 the $n-1$- dimensional algebraic torus $H= ({\bf C}^*)^n/{\bf C}^*$.
Our main object of study in this paper will be the Chow quotient $G(k,n)//H$.

\vskip .3cm

 \noindent (1.1.2) For each subset $I \i \{1,...., n\}$ denote by $L_I$
the coordinate subspace in ${\bf C}^n$ defined by equations
$x_i =0, i\in I$ and by ${\bf C}^I$ the coordinate subspace spanned
by basis vectors from $I$. Thus the codimension of $L_I$ and
the dimension of ${\bf C}^I$ equal to $|I|$, the cardinality of $I$.

 Call a $k$-dimensional subspace
$L\in G(k,n)$ {\it generic} if for any $I\i \{1,..., n\}, |I|=k$
we have $L\cap L_I= 0$. The space $G^0(k,n)$ of all generic subspaces
is an open $H$- invariant subset in $G(k,n)$. It is called the
{\it generic stratum}. It will serve as the open set $U$ from \S0.
\vskip .3cm

 \noindent (1.1.3) The Grassmannian $G(k,n)$ can be seen as
 the variety of $(k-1)$- dimensional projective subspaces in the
projective space
$P^{n-1}$. Using the terminology going back to Pl\"ucker, one usually calls
$(n-1)$- dimensional families of subspaces
in $P^{n-1}$ {\it complexes}.
\vskip .3cm

\proclaim (1.1.4) Definition. By a  Lie complex we shall
 mean an algebraic subvariety in $G(k,n)$ which
is the closure of the $H$- orbit $\overline {H.L}$ of some
generic subspace $L\in G^0(k,n)$.
\vskip .3cm

\proclaim (1.1.5) Proposition. {\rm [19]}  Each Lie complex
 is a $(n-1)$- dimensional variety containing all the $H$- fixed points
 on $G(k,n)$ given by coordinate subspaces ${\bf C}^I,\, |I|=k$.
These $n\choose k$ points are the only singular points of a Lie complex.
Near each of these points a Lie complex looks like the cone over
 $P^{k-1}\times P^{n-k-1}$ in the Segre embedding.
\vskip .3cm

\noindent {\bf (1.1.6) Example.}  Lie complexes in $G(2,4)$
 were extensively
studied in classical literature under the name of
{\it tetrahedral complexes}.
see [2], [27] and references in [20]. Let us describe them in more detail.
 Let $x_1,...,x_4$ be homogeneous coordinates in $P^3$
and $L_i$ be the coordinate plane $\{x_i=0\}$. The configuration
of four planes $L_i$ can be thought of as a tetrahedron.
A line $l\in P^3$ lies in generic stratum $G^0(2,4)$ if and only if
 it does not intersect any of the 6 lines given by the edges
 of our tetrahedron. For such a line the four points of
 intersections $l\cap L_i$ are distinct
 and, as any four distinct points on a projective line,
 posess the cross-ratio
 $r(l\cap L_1,...,l\cap L_4)\in {\bf C}-\{0,1\}$.
Let $\lambda \in {\bf C}-\{0,1\}$ be a fixed number.
The tetrahedral complex
$K_\lambda$ is, by definition, the closure of the set of those
$l\in G^0(2,4)$ for which the cross-ratio $r(l\cap L_1,...,l\cap L_4)$
 equals $\lambda$. Its equation in Pl\"ucker coordinates is
$$p_{12}p_{34} + \lambda p_{13}p_{24} =0.$$
This can be commented as follows.
The classical Pl\"ucker relation gives that three quadratic polynomials
 $p_{12}p_{34}, p_{13}p_{24}$ and $p_{14}p_{23}$ on $G(2,4)$
are linearly dependent i.e. generate a 1-dimensional linear system
(pencil) of hypersurfaces. The tetrahedral complexes are just
hypersurfaces from this pencil.
They are, therefore,  particular cases of quadratic line
complexes.
As was pointed out in [20], the definition of a tetrahedral complex
as the closure of a torus orbit is due to F.Klein and S.Lie.
\vskip .3cm

 \noindent(1.1.7) Clearly all Lie complexes represent the same class  in
$(2n-2)$- dimensional homology of the Grassmannian. Denote this class
$\delta$.
 Let us recall
an explicit formula for $\delta$  found
by A.Klyachko [36]. For any Young diagram
$\alpha =(\alpha_1 \geq ... \geq \alpha_k)$
with no more than $k$ rows and no more than $(n-k)$ columns we shall
 denote by $|\alpha| =\sum \alpha_i$ the number of cells in $\alpha$
 and by
 $\sigma_\alpha$ the Schubert class in $H_{2|\alpha|}(G(k,n))$
 corresponding to $\alpha$ ( see  [24] and \S 3.9 below
 for details on Schubert cycles). These classes form an
integral basis in the homology and
the formula of Klyachko gives a decomposition of
$\delta$ with respect to this basis.
\vskip .3cm

\proclaim (1.1.8) Proposition. {\rm [36]} Let $\alpha$ be a Young diagram
 with $(n-1)$ cells. The coefficient at $\sigma_\alpha$ in the decomposition
 of the fundamental class $\delta$ of a Lie complex with respect to Schubert
 cycles equals
$$\sum_{i=0}^k (-1)^i {n\choose i} dim \Sigma^\alpha ({\bf C}^{k-i}),$$
where  $\Sigma^\alpha ({\bf C}^{k-i})$ is the irreducible representation of
 $GL(k-i)$ with highest weight $\alpha$.
\vskip .3cm

\noindent {\bf (1.1.9) Example.} For a Lie complex in $G(2,n)$ the above
formula
 gives
$$\delta = (n-2)\sigma_{n-2,1} + (n-4)\sigma_{n-3,2} +
(n-6)\sigma_{n-4,3} +...$$
In particular, for  a Lie complex in $G(2,4)$ the formula gives
$\delta =2\sigma_{2,1}$ and $\sigma_{2,1}$ is the class of hyperplane
section of $G(2,4)$ in the Pl\"ucker embedding. This agrees with the
fact that Lie
complexes in $G(2,4)$ are quadratic complexes.

\vskip .3cm

\noindent (1.1.10) The collection of all Lie complexes is naturally
identified with  $G^0(k,n)/H$, the quotient of the generic
stratum (1.1.2).
 We are interested in the Chow quotient $G(k,n)//H$ which is a projective
subvariety in the Chow variety ${\cal C}_{n-1}(G(k,n),\delta)$,
 namely the closure of the set  of all Lie complexes.

Any algebraic cycle from $G(k,n)//H$ will be called a
{\it generalized Lie complex}.  It is our point of view that generalized
 Lie complexes are the "right" generalizations of generic torus orbits
in the Grassmannian.
 We shall see later that each generalized Lie complex can be seen as a
(possibly reducible) algebraic subvariety in $G(k,n)$.
\vskip .3cm

\noindent {\bf (1.1.11) Example.} The Chow quotient $G(2,4)//H$
is isomorphic to the projective line $P^1$. The isomorphism
 is given by the cross-ratio of four points of intersection
$l\cap L_i$  in Example (1.1.6). There are exactly three generalized
 Lie complexes
 in $G(2,4)$ which are not closures of single orbits
(i.e. are not genuine Lie complexes). They are limit positions
 of tetrahedral complexes  corresponding to values $0,1,\infty$
not taken by the cross ratio. Denote by $Z_{ij}\i G(2,4)$
the space of lines intersecting the coordinate line
$x_i=x_j=0$ (the edge of the tetrahedron).
This is a linear section of $G(2,4)$ given by the equation $p_{ij}=0$.
 The three limit complexes are
$$Z_{12}+Z_{34}, Z_{13}+Z_{24}, Z_{14}+Z_{23}.$$

\vskip .5cm
\centerline {\bf (1.2) Chow strata and matroid decompositions of
the hypersimplex.}
\vskip .5cm

\noindent (1.2.1)
Call two $k$- dimensional linear subspaces $L, L'\i {\bf C}^n$ {\it equivalent}
if $dim (L\cap L_I) = dim (L'\cap L_I)$ for any $I\i \{1,...,n\}$.
 Corresponding equivalence classes are called  {\it strata}.
 They are $H$- invariant subsets in $G(k,n)$. A {\it base} of
 a subspace $L\in G(k,n)$ is a $k$-element subset $I\i \{1,...,n\}$
such that $L\cap L_I =0$. It is well- known that two subspaces $L$
and $L'$ lie in the same stratum (i.e. are equivalent) if and anly
if their sets of bases coincide. As a particular case we obtain the
 generic stratum $G^0(k,n)\i G(k,n)$ defined as follows. A space
$L$ lies in $G^0(k,n)$ if and only if each $k$- element subset
is a base for $L$.

This stratification was introduced in [18,19,21]. The set of bases for any
subspace $L\i G(k,n)$ introduces on $\{1,...,n\}$ the
structure of a {\it matroid} of rank $k$. Because of this, this stratification
is often referred to as the {\it matroid stratification} of the Grassmannian.
\vskip .3cm

\noindent (1.2.2) It was remarked in [19], \S 5.1 that
the matroid stratification of the Grassmannian is not a
 stratification in the sense of Whitney.  In particular,
 the closure of a stratum may happen not to be a union of other strata.
\vskip .3cm

\noindent (1.2.3) Let $e_1,...,e_n$ be standard basis vectors in the
 coordinate space
 ${\bf R}^n$. We define the convex polytope $\Delta (k,n)$ called the
$(k,n)$- {\it hypersimplex}
 to be the convex hull of  $n\choose k$ points $e_{i_1}+... +e_{i_k}$
 where $1\leq i_1 < ...< i_k\leq n$.

All these points are vertices of
$\Delta (k,n)$. We shall denote these vertices shortly by
$e_I=\sum_{i\in I}e_i$ where $I\i \{1, ... ,n\},\, |I|=k$.
For any subspace $L\in G(k,n)$ we define its {\it matroid polytope}
$M(L)$ as the convex hull of $e_I$, where $I$ runs over all bases for $L$.
Thus $\Delta (k,n)$ itself is the matroid polytope for a generic subspace.

 The
 hypersimplex  was introduced in [17] and serves as a combinatorial
model both for the Grassmannian with torus action and for any Lie complex.
\vskip .3cm

\proclaim (1.2.4) Proposition. {\rm [19]} Let $L\in G(k,n)$. Then:\hfill\break
a) Any edge of $M(L)$ is parallel to a vector of the form
$e_i- e_j, i\neq j$.\hfill\break
b) The (complex) dimension of the orbit $H.L$ coincides with the real
 dimension of the polytope $M(L)$. The closure $\overline {H.L}$ is a
 projective, normal, toric variety and $M(L)$ is the corresponding polytope
(i.e. the fan of $\overline {H.L}$ is the normal fan of $M(L)$).
In particular:\hfill\break
c) Any Lie complex is a projective toric variety and the
 corresponding polytope is $\Delta (k,n)$. So $p$- dimensional
 $H$-orbits on any Lie complex are in bijection with $p$-dimensional
faces of $\Delta (k,n)$.

The following description of faces of $\Delta (k,n)$  was given in [17].

\proclaim (1.2.5) Proposition. a) Each face of $\Delta (k,n)$
is itself a hypersimplex. \hfill\break
b) Edges of $\Delta (k,n)$ are segments $[e_I, e_J]$ where $J$
differs from $I$ by replacing one element $i\in I$ by another
 $j\notin I$. \hfill\break
c) For $k>1$ there are exactly $2n$ facets (faces of codimension 1)
 of $\Delta (k,n)$. They are
$$\Gamma_i^+ = {\rm Conv} \{e_I, i\in I\} \,\,\,\,{\rm and}
\,\,\,\, \Gamma_i^- = {\rm Conv}\{e_I, i\notin I\}$$
for $1\leq i\leq n$. Each polytope $\Gamma_i^+$ is linearly isomorphic
to the hypersimplex $\Delta (k-1, n-1)$ whereas each $\Gamma_i^-$ is
 isomorphic to $\Delta(k, n-1)$.
\vskip .3cm

\noindent (1.2.6) By a {\it matroid polytope} in $\Delta (k,n)$ we
 shall mean any subpolytope $M\i \Delta(k,n)$ whose vertices are
 among vertices of $\Delta (k,n)$ and edges a have the form described
in part a)  of  Proposition  1.2.4 (i.e. are among edges of $\Delta (k,n)$).
According to Proposition 1.2.5 the polytope $M(L)$ for any subspace
$L\in G(k,n)$ is a matroid polytope. Matroid polytopes of such form
are called {\it realizable}.

The notion of matroid polytope in $\Delta (k,n)$ was introduced in
[19,21]. It was shown in these papers that such polytopes are in
bijection with the structures of rank $k$ matroid on a set $\{1,...,n\}$.

 \vskip .3cm

\noindent (1.2.7) Consider the Pl\"ucker embedding of the Grassmannian
 $G(k,n)$ into the ${n\choose k}-1 $-dimensional projective space
$P(\bigwedge ^k {\bf C}^n)$. The homogeneous coordinates in this
 projective spaces will be denoted $p_I, I\i \{1,..., n\},\, |I|=k$.
The $H$-action on $G(k,n)$ extends to the whole $P(\bigwedge ^k {\bf C}^n)$.
The matroid polytope of a subspace $L\in G(k,n)$ is the image of the orbit
 closure $\overline {H.L}\i G(k,n)\i P(\bigwedge ^k {\bf C}^n) $  under the
momentum map $\mu :P(\bigwedge ^k {\bf C}^n)\rightarrow \Delta (k,n)$
defined as follows [19,21]:
$$\mu (x)= {{\sum_{|I|=k}p_I(x).e_I}\over {\sum_{|I|=k} p_I(x)}}.
\leqno (1.2.8)$$

\noindent (1.2.9) Since $G(k,n)$ is embedded equivariantly into
 $P(\bigwedge ^k {\bf C}^n)$, we obtain the embedding of Chow quotients
$$G(k,n)//H \hookrightarrow P(\bigwedge ^k {\bf C}^n)//H.\leqno (1.2.10)$$
The latter quotient  is, according to Example 0.2, a toric variety of
dimension ${n\choose k}-n$ and the corresponding polytope is the
secondary polytope of the hypersimplex $\Delta (k,n)$. By comparing
 the two Chow quotients
we deduce from [30,31] the following proposition.

\proclaim (1.2.11) Proposition. Let $Z=\sum c_iZ_i$ be a cycle
from $G(k,n)//H$. Then
:\hfill\break
a) Each component $Z_i$ is a closure of some $(n-1)$- dimensional
$H$- orbit $Z_i^0$.\hfill\break
b) Let $M(Z_i)= \mu(Z_i)$ be the matroid polytope of any subspace
$L\in Z_i^0$ (or, what is the same , the image of $Z_i$ under the
 momentum map. Then the polytopes $M(Z_i)$ form a polyhedral
 decomposition of $\Delta (k,n)$.

\noindent {\bf (1.2.12) Example.} Consider the case $k=2, n=4$.
 The hypersimplex
$\Delta (2,4)$ is the 3-dimensional octohedron. Each of the three
generalized Lie complexes
from (1.1.12) gives a decomposition of this octohedron into a union
of two pyramids with a common
quadrangular face.

We have the embedding (1.2.10)
of $G(2,4)//H = P^1$ into  $P(\bigwedge ^2 ({\bf C}^4))//H$ which is a
 toric variety of
dimension 2. This variety is isomorphic to the projective plane $P^2$.
To see this, let us show that the secondary polytope (polygon,
 in our case) $\Sigma$ of the octohedron $\Delta (2,4)$ is in fact
a triangle. Indeed, by definition (0.2.7) vertices of $\Sigma$ are in
bijection with regular triangulations of $\Delta(2,4)$. Each
 triangulation of $\Delta (2,4)$ can be obtained as follows.
Take any  decomposition of $\Delta (2,4)$ into two pyramids as
 above and then decompose each of these pyramids
into two tetrahedra in a compatible way:
\vbox to 4cm{\vfill\vfill (1.2.13)\vfill\vfill}

 Thus there are 3 triangulations which correspond to vertices
 of $\Sigma$ and three pyramidal decompositions which correspond
to edges of $\Sigma$ and hence $\Sigma$ is a triangle.

Since the symmetry group of the octohedron acts on $\Sigma$,
it is a regular triangle. Hence the toric variety corresponding
to $\Sigma$, has $P^2$ as its normalization. The fact that this
 variety is normal can be established by direct computation of
vertices of $\Sigma$ as points of the integer lattice ${\bf Z}^6$
 (according to (0.2.6)) which we leave to the reader.

So  the toric variety $P(\bigwedge ^2 ({\bf C}^4))//H$
is a projective plane. The subvariety $G(2,4)//H = P^1$ is a
 conic in this projective plane inscribed into the coordinate triangle:

\vbox to 4cm{\vfill\vfill  (1.2.14)\vfill\vfill}

\vskip .3cm

\proclaim (1.2.15) Proposition. Let $Z=\sum c_iZ_i$ be a cycle
from $G(k,n)//H$. Then all the multiplicities $c_i$ equal 1 (or 0).

\noindent {\sl Proof.} The recipe for calculation of $c_i$ given
in [30,31] is  the following. We consider the affine {\bf Z}- lattice
$\Xi_i$ generated by the vertices of the polytope $M(Z_i)$ which is
imbedded into the affine lattice $\Xi$ generated by all the vertices
 of $\Delta (k,n)$. Then $c_i =[\Xi :\Xi_i]$. Let us show that in
fact $\Xi =\Xi_i$.  Choose some vertex $e_I$ of $M(Z_i)$. Subtracting
 it from points of $\Xi$ , we  identify $\Xi$ with $\{(a_1,..., a_n)
\in {\bf Z}^n :\sum a_i =0 \}$. Consider now all edges of $M(Z_i)$
 containing $e_I$. By Proposition 1.3 b), they all have the form
$e_j -e_l$, where $j\in I, l\notin I$. Since $M(Z_i)$ has full dimension,
 there are  at least $n-1$ independent edges. However, any $n-1$
independent vectors of the form $e_j-e_l$ generate the lattice
$\Xi$. $\quad\triangleleft$

By the above proposition, generalized Lie complexes (= "infinite"
 points of the Chow quotient $G(k,n)//H$) can be thought of as usual
 reducible subvarieties (instead of cycles) in the Grassmannian, what
further justifies their name. We will therefore denote these complexes
 by $Z=\bigcup Z_i$ to emphasize that they are varieties.

\vskip .3cm

\proclaim (1.2.16) Definition. Two cycles $Z, Z' \in G(k,n)//H$
are called equivalent if the corresponding polyhedral decompositions
of $\Delta (k,n)$ coincide. Equivalence classes under this relation
are called Chow strata.

By considering again the Pl\"ucker embedding we see that our
 stratification of \hfill\break $G(k,n)//H$ is induced from the
stratification of the toric variety $P(\bigwedge ^k {\bf C}^n)//H $
 given by the torus orbits. Each Chow stratum can be specified by a
 finite list of usual strata corresponding to individual matroid
polytopes from the polyhedral decomposition.

\vskip .3cm

\proclaim (1.2.17) Definition. A polyhedral decomposition ${\cal P}$
of the hypersimplex $\Delta (k,n)$ is called a matroid decomposition
 if all the polytopes from ${\cal P}$ are matroid polytopes (1.2.7).
 A matroid decomposition is called realizable if it comes from a
generalized Lie complex (1.2.11).

Thus matroid decompositions of $\Delta (k,n)$ are precisely the labels
by which Chow strata are labelled. The notion of a matroid polytope
being equivalent to that of matroid (1.2.1), a matroid decomposition
represents a new kind of combinatorial structure --- a collection of
 usual matroids with certain properties (that the corresponding
 polytopes form a decomposition of the hypersimplex).

\vskip 1cm

\centerline {\bf (1.3) Example: matroid decompositions of $\Delta (2,n)$.}

\vskip .5cm

In this section we give a complete description of matroid decompositions of
the hypersimplex $\Delta (2,n)$.
The structure involved will turn out to be identical to those in the
 description of stable $n$ -pointed curves of Grothendieck [12]
 and Knudsen [37].

 Recall that vertices of $\Delta (2,n)$ are
of the form $e_{ij} := e_i+e_j, i\neq j, 1\leq i,j\leq n$, where
$e_i$ are the standard basis vectors of ${\bf R}^n$.
\vskip .3cm

\proclaim (1.3.1) Proposition. a) Matroid polytopes in $\Delta (2,n)$
are in bijection with pairs $(J,R)$ where $J\i \{1,...,n\}$ is a
non-empty subset and $R$
is an equivalence relation on $J$ with at least 2 equivalence classes.
The matroid polytope $M(J,R)$ corresponding to $(J,R)$ above has vertices
$e_{ij}$ where $i,j\in J$ are such that $iRj$ does not hold.\hfill\break
b) The dimension of $M(J,R)$ equals $|J|-1$ if $R$ has $\geq 3$
equivalence classes and equals $|J|-2$ if $R$ has exactly 2 equivalence
classes.

\noindent {\sl Proof:} This follows from ([21], Example 1.10 and
Proposition 4,\S 2).
\vskip .3cm

\proclaim (1.3.2) Corollary. Matroid polytopes in $\Delta(2,n)$
 which have full dimension $(n-1)$, are in bijection with equivalence
 relations on $\{1,...,n\}$ with $\geq 3$ equivalence classes.

Thus matroid decompositions of $\Delta (2,n)$ are certain "compatible"
 systems
of equivalence relations on the same set $\{1,...,n\}$. We are going to
describe them.
\vskip .3cm

\noindent (1.3.3) By a graph we mean a finite 1-dimensional simplicial
complex.
So a graph $\Gamma$ is defined by its set of vertices $\Gamma_0$ and
 the set of edges $\Gamma_1$ together with the incidence relation
 connecting these sets.
If $v$ is a vertex of a graph $\Gamma$, the {\it valency} of $v$ is,
 by definition, the number of edges containing $v$.

By a {\it tree} we mean a connected graph $T$ without loops such that
 every vertex of $T$ has the valency either 1 or $\geq 3$. The vertices
 of valency 1 will be called {\it endpoints} of $T$. For any two vertices
$v,w$ of a tree $T$ there is a unique edge path without repetitions
joining these vertices. This path will be denoted $[v,w]$.

Let $A_1,..., A_n$ be formal symbols. By a tree bounding the endpoints
 $A_1,...,A_n$, we shall mean a tree $T$ with exactly $n$ endpoints
which are put into bijection (or just identified) with symbols $A_i$.
Two such trees $T,T'$ are called isomorphic if there
is an isomorphism of graphs $T\rightarrow T'$ preserving $A_i$.

\vskip .3cm

\noindent (1.3.4) Let $T$ be a tree bounding endpoints $A_1,...,A_n$.
Any vertex of $T$ which is not an endpoint will be called interior. Let
  $v\in T_0$ be an interior vertex. We define an equivalence
relation $\cong _v$ on $\{1,...,n\}$ be setting $i\cong _v j$ if the
 edge path
$[A_i,A_j]$ does not contain the vertex $v$.

In other words, the deletion of the vertex $v$ splits the tree into
several
connected components and $i\cong_vj$ if the endpoints $A_i,A_j$ are
 situated
in the same component. The equivalence classes under $\cong_v$ are in bijection
with edges of $T$ containing $v$.
\vskip .3cm

\noindent (1.3.5) Proposition (1.2.5) implies that $\Delta (2,n)$ has
$n$ facets (faces of codimension 1) $\Gamma_i^+ = {\rm Conv}\{e_{ij},
 j\neq i\}$ which are $(n-2)$- dimensional simplices. It is clear that
 in any polyhedral decomposition of $\Delta (2,n)$ each $\Gamma_i^+$
is a facet of exactly one polytope from decomposition.

\proclaim (1.3.6) Theorem. Matroid decompositions of the hypersimplex
$\Delta (2,n)$ are in bijection with isomorphism classes of trees
bounding endpoints $A_1,...,A_n$.\hfill\break
 $\bullet$ Explicitly, if $T$ is such a tree, the corresponding
 decomposition ${\cal P}(T)$ consists of matroid polytopes $M(\cong_v)$
(Proposition (1.3.1)) for all interior vertices $v$ of $T$. \hfill\break
$\bullet$ Conversely, the tree $T$ can be recovered from the corresponding
 matroid decomposition ${\cal P}$ as follows. Internal vertices of
 $T$ are barycenters of  polytopes (of maximal dimension) from
 ${\cal P}$. Endpoints of $T$ are barycenters of facets $\Gamma_i^+$.
 The barycenter of each $\Gamma_i^+$ is joined to the barycenter of
 the unique polytope from ${\cal P}$ containing $\Gamma_i$;
the barycenters of two polytopes from ${\cal P}$ are joined if
 and only if these polytopes have a common facet.

\vskip .3cm

\noindent {\bf (1.3.7) Remark.} Let $v$ be an interior vertex of
the tree $T$. The vertices of the polytope $M(\cong_v)$ are those
 vertices $e_{ij} = e_i+e_j$ of $\Delta (2,n)$ for which the edge
path $[A_i,A_j]$ does contain $v$.

\vskip .3cm

\noindent {\sl Proof of (1.3.6):} Let $T$ be any tree bounding
$A_1,...,A_n$. Let us show that the collection of polytopes
$M(\cong _v)$ forms a polyhedral decomposition of $\Delta (2,n)$.
By definition, this means that the two properties hold:
\vskip .3cm

\noindent {\bf (1.3.8)} Intersection of any two polytopes
 $M(\cong_v), M(\cong_w)$ is a common face of both of them.
\vskip .3cm

\noindent {\bf (1.3.9)} The union of the polytopes $M(\cong_v)$
is the whole hypersimplex $\Delta (2,n)$.
\vskip .3cm

\noindent (1.3.10) \underbar {Proof of {\bf (1.3.8)}.} We shall
prove a stronger statement: that $M(\cong_v)\cap M(\cong_w)$ is
 the convex hull of vertices common to $M(\cong_v), M(\cong_w)$
. By (1.3.7), a vertex $e_{ij}$ is common to $M(\cong_v), M(\cong_w)$
if the edge path $[A_i,A_j]$ contains $[v,w]$ as a sub- path.
Let us subdivide the set $\{1,...,n\}$ into three parts: $X_+, X_-, X_0$.
 We set $i\in X_+$ if the edge path $[A_i,v]$ does not contain points on
 $[v,w]$ other than $v$. We set $i\in X_-$ if the edge path $[A_i,w]$
 does not contain points of $[v,w]$ other than $w$. We set $i\in X_0$
in all other cases:

\vbox to 4cm{\vfill\vfill  (1.3.11)\vfill\vfill}

Recall that $\Delta(2,n)$ lies in ${\bf R}^n$ as the convex hull of
 sums $e_{ij} = e_i+e_j$ of two distinct basis vectors. Consider the
linear function $f$ on ${\bf R}^n$ such that $f(e_i) = +1, -1$ or $0$
 if $i\in X_+, X_-$ or $X_0$ respectively. Then $f$ is non-negative
 on all vertices of $M(\cong_v)$ and non-positive on all vertices of
 $M(\cong_w)$. The only vertices of  $M(\cong_v), M(\cong_w)$ for
which $f=0$ are the vertices common to both these polytopes.
Assertion (1.3.8) is proven.
\vskip .3cm

\noindent (1.3.12) \underbar{Proof of {\bf (1.3.9)}}. It suffices
 to show that any face of codimension 1 of any $M(\cong_v)$ lies
either on the boundary of $\Delta (2,n)$ or is a face of another
 polytope $M(\cong_w)$.
The following description of facets (=faces of codimension 1) of
 matroid polytopes follows from ([21], \S 2, Theorem 5).
\vskip .3cm

\proclaim (1.3.12.1) Proposition. Let $J=\{1,...,n\}$ and
$M = M(J,R)$ be the matroid polytope of full dimension corresponding
 to an equivalence relation $R$ on $J$. Its facets  are the following:
\hfill\break
(1) Facets $\Gamma^+_j(M)$ defined for any $j$ unless $\{j\}$ is an
equivalence class in itself and the total number of classes is 3.
 This facet is the matroid polytope $M(J',R')$, where
$J' = \{1,...,n\}-\{j\}$ and $R'$ is the equivalence
relation induced by $R$. It lies entirely in the boundary of $\Delta (2,n)$.
\hfill\break
(2) Facet $\Gamma^-_K(M)$ defined for any equivalence class
$K\in J/R$. This is
the matroid polytope $M(J, R'')$, where $R''$ is the equivalence
relation with
only two classes of which one is $K$ and the other is formed by
all elements not in $K$.

The notation $\Gamma^\pm$ is compatible with the notation for the
facet of the
full hypersimplex $\Delta(2,n)$ introduced in Proposition 1.2.5.

\vskip .3cm

\proclaim (1.3.12.2) Corollary. Let $T$ be a tree bounding $A_1,...,A_n$
 and $v\in T$ be an interior vertex. The facets of the matroid polytope
 $M(\cong_v)$ not lying in the boundary of $\Delta (2,n)$ are in bijection
 with edges  of $T$ containing $v$ whose second end is also an interior
vertex. The facet corresponding to such an edge $e$ is of the form
$\Gamma^-_K$ where $K$ is the $\cong_v$ -equivalence class corresponding
 to $e$.

Now the assertion (1.3.9) follows from Corollary 1.3.12.2 since
every edge of the tree $T$ joins two vertices and the two matroid
 polytopes corresponding to these vertices have a common facet.
\vskip .3cm

\noindent (1.3.15) We have proven that any tree bounding $A_1,...,A_n$
 gives a matroid decomposition of $\Delta (2,n)$. Conversely, let
 ${\cal P}$ be any matroid decomposition. By taking barycenters of
 polytopes from ${\cal P}$ and joining them as prescribed in Theorem 1.3.6,
 we obtain a certain graph $T$.  Let us show that $T$ is a tree which
 generates the decomposition ${\cal P}$. The fact that $T$ is a tree
 follows from the next lemma

\vskip .3cm
\proclaim (1.3.15.1) Lemma. Let $M\i \Delta (2,n)$ be a matroid polytope
 of full dimension and $\Gamma\i M$ be a facet not lying on the boundary
 of $\Delta (2,n)$. Then $\Gamma$ is equal to the intersection of the
whole $\Delta(2,n)$ with a hyperplane.

\noindent {\sl Proof:} According to Proposition 1.3.12.1, the
polytope $\Gamma$ has the form $M(J,R)$, where $R$ is an equivalence
 relation on $J=\{1,...,n\}$ with only two equivalence classes,
say $A$ and $B$. Define a linear function $g$ on ${\bf R}^n$ whose
value on the basis vector $e_i$ equals 1, if $i\in A$ and equals
 $(-1)$, if $i\in B$. Then $\Gamma$ is the intersection of
$\Delta (2,n)$ with the kernel of $g$. Lemma (1.3.15.1) is proven.

\vskip .3cm

\noindent (1.3.15.2) To finish the proof of Theorem 1.3.6, it
 remains to show that
the tree obtained from the matroid decomposition ${\cal P}$ in (1.3.15),
generates ${\cal P}$. This checking is left to the reader.\vskip .3cm

\noindent {\bf (1.3.16) Example.} There are four matroid decomposition
 of the octohedron $\Delta (2,4)$: one consists of $\Delta (2,4)$ itself
 and each of the others decomposes the octohedron into two pyramids
(Example 1.2.12). These decompositions correspond to the following
 trees bounding endpoints $A_1,...,A_4$:

\vbox to 4cm{\vfill\vfill (1.3.17)\vfill\vfill}

\vskip .3cm

\noindent (1.3.18) We shall show in \S 4 that all matroid decompositions
of $\Delta (2,n)$ are realizable.

\hfill\vfill\eject

\centerline {\bf (1.4) Relation to the secondary variety for
 the product of two simplices.}
\vskip .5cm

In this section we shall compare the Chow quotient $G(k,n)//H$ with a
toric variety of the same dimension. This toric variety will correspond
to the convex polytope which is the secondary polytope for the
 product of two simplices.
\vskip .3cm

\noindent (1.4.1)
Let $P\i {\bf R}^m$ be any convex polytope and $x\in P$ be any vertex.
Denote by $N_xP$ the union of all half-lines drawn from $x$ through all
 the points of $P$. This is an affine cone which we call the {\it normal cone}
 to $P$ at $x$. The base of this cone i.e. a transversal  section of $N_xP$
by an affine hyperplane will be called the {\it vertex figure} of $P$ at $x$.
Thus vertices of the vertex figure correspond to edges of $P$ containing $x$.
 \vskip .3cm

\noindent(1.4.2) Let $A$ be any finite set. Denote by $\Delta ^A$ the
simplex (of dimension $|A|-1$) whose set of vertices is $A$.
By definition, $\Delta ^A$ is the subset in the space ${\bf R}^A$
of functions $A\rightarrow R$ consisting of functions $f(a)$ such
 that $f(a)\geq 0, \forall a$ and $\sum f(a) =1$. To any $a\in A$
there corresponds a vertex $\delta_a$ of $\Delta^A$. This is the
function $A\rightarrow {\bf R}$ taking $a$ to 1 and other elements to 0.

\vskip .3cm

\noindent(1.4.3) Let $e_I$ be a vertex of the hypersimplex $\Delta(k,n)$.
The corresponding vertex figure  is the product of two simplices
 $\Delta ^{k-1}\times \Delta ^{n-k-1}$ or, in the more invariant
notation of (1.4.2), $\Delta^I\times \Delta^{\bar I}$, where $\bar I$
is the complement to $I$.

Indeed, edges of $\Delta(k,n)$ containing $e_I$, are $[e_I, e_I +e_j-e_i]$
where $i\in I,\, j\notin I$.  The required isomorphism takes such an edge
to the vertex $(\delta_i, \delta_j)$ of $\Delta^I\times \Delta^{\bar I}$.

The toric variety associated to the polytope
$\Delta ^{k-1}\times \Delta ^{n-k-1}$ is the product of projective
spaces $P^{k-1}\times P^{n-k-1}$ and the structure of $\Delta (k,n)$
 near a vertex corresponds to the structure of a Lie complex near its
 singular point ( Proposition 1.1.5).

 \vskip .3cm

\proclaim (1.4.4) Proposition-Definition. Let ${\cal P}$ be a matroid
decomposition of $\Delta (k,n)$ and $e_I\in \Delta (k,n)$- a vertex.
Let ${\cal P}_I$ be the induced polyhedral decomposition of the vertex
 figure $\Delta^I\times\Delta^{\bar I}$ . Then all the vertices of
 polytopes constituting ${\cal P}_I$ lie among the vertices of
$\Delta^I\times\Delta^{\bar I}$. If ${\cal P}$ is a realizable
matroid decomposition then ${\cal P}_I$ is a regular polyhedral
 subdivision of $\Delta^I\times\Delta^{\bar I}$.

Recall [22] that a polyhedral subdivision is called regular if it admits
 a strictly convex piecewise- linear function.

\noindent {\sl Proof:} Vertices of polytopes from ${\cal P}_I$ correspond
to edges of polytopes from ${\cal P}$ containing $e_I$. Since all these
polytopes are matroid polytopes, the edges in question correspond to
vertices of $\Delta^I\times\Delta^{\bar I}$.  If ${\cal P}$ is
realizable then it is regular as a polyhedral subdivision of
 $\Delta (k,n)$ and so is ${\cal P}_I$. $\quad\triangleleft$
\vskip .3cm

\noindent (1.4.5) Let $I\i \{1,...,n\}$ be a $k$ -element subset.
 The coordinate subspace ${\bf C}^I\in G(k,n)$ is a fixed point under
 the action of the torus $H \quad$ (1.1.1). Therefore we have the action
 of $H$ on the tangent space
 $T_I\,:=\,T_{{\bf C}^I}G(k,n)$.

The tangent space to $G(k,n)$ at any point $L$ is canonically identified,
see [47], with $Hom (L, {\bf C}^n/L)$. Therefore we have the isomorphism
of $H$ -modules
$$T_I = Hom ({\bf C}^I, \, {\bf C}^{\bar I})\leqno (1.4.6)$$
In other words, $T_I$ is decomposed into $k(n-k)$ one-dimensional
weight subspaces $V_{ij}, i\in I, j\notin I$ such that for any
 $t = (t_1,...,t_n)\in H$ and any $v\in V_{ij}$ one has
$$(t_1,...,t_n).v = (t_i/t_j)v.$$
\vskip .3cm

\noindent (1.4.7)
The character lattice of the torus $H$ is the sublattice in
${\bf Z}^n$ consisting of vectors with the sum of coordinates
 equal to 0. The character corresponding
to the subspace $V_{ij}$ is the vector $e_i-e_j \in {\bf Z}^n$.
 The collection of all vectors $e_i-e_j, i\in I, j\notin I$,
 forms the set of vertices of the simplex
$\Delta^I\times\Delta^{\bar I}$.

\vskip .3cm

\noindent (1.4.8) Call a point $v$ of the tangent space $T_I$
(and the corresponding point of the projectivization $P(T_I)$)
 {\it generic} if all the weight components of $v$ are non-zero.
By a generic $H$ -orbit in $T_I$ or $P(T_I)$ we shall mean the
orbit of a generic point.

\vskip .3cm

\noindent (1.4.9)
For the torus orbits in the Grassmannian $G(k,n)$ we also have a
notion of genericity introduced in (1.1.2). Closures of generic orbits
in $G(k,n)$
were called Lie complexes. Let   $Z=\overline {H.L}, L\in G^0(k,n)$
 be any Lie complex  and
 $TC_IZ\,:= TC_{{\bf C}^I}Z\i T_I$ --
its tangent cone at the point ${\bf C}^I$.
It follows from  (1.1.5) that
$TC_IZ$   is the closure of a generic $H$- orbit in $T_I$.
So we obtain the following proposition.
\vskip .3cm

\proclaim (1.4.10) Proposition. Let us identify the quotient
 $G^0(k,n)/H$ (i.e. the set of generic (1.1.2) $H$- orbits on
 $G(k,n)$) with the set of Lie complexes. Then
 the correspondence $Z\mapsto$ (the projectivization of
$TC_{{\bf C}^I}Z$) defines an open embedding  of $G^0(k,n)/H$
 into the set of generic $H$- orbits in $P(T_I)$.

\vskip .3cm

\noindent (1.4.11) We are going to compare the Chow quotient
$G(k,n)//H$ with $P(T_I)//H$. The latter variety is, according
 to results recalled in \S (0.2), a projective toric variety of
 the same dimension $k(n-k)-n+1$ as $G(k,n)//H$. Due to (1.4.7)
and Theorem 0.2.8, the toric variety $P(T_I)//H$ corresponds to
 the secondary polytope of $\Delta^I\times \Delta^{\bar I}$.
We shall therefore call this variety the {\it secondary variety} of
$\Delta^I\times \Delta^{\bar I}$.
\vskip .3cm

\proclaim (1.4.12) Proposition. The open embedding  from (1.4.10)
 extends to a regular birational  (in particular, surjective) morphism
$$f_I: \, G(k,n)//H\rightarrow P(T_I)//H.\leqno (1.4.13)$$
This morphism takes any generalized Lie complex ${\cal L}$ to the
projectivization of its tangent cone at ${\bf C}^I$ (which is regarded
 as an algebraic cycle with multiplicities 0 or 1 in $P(T_I)$).

\noindent {\sl Proof:} Let $U_I\i G(k,n)$ be the affine chart consisting
of $k$ -dimensional subspace $L$ such that
 $L\oplus {\bf C}^{\bar I} = {\bf C}^n$.
Each such subspace can be regarded as the graph of a linear operator
from ${\bf C}^I$ to ${\bf C}^{\bar I}$. This correspondence identifies
 $U_I$ with the
space of $k\times (n-k)$ -matrices thus introducing coordinates
 $z_{ij}, i\in I, j\notin I$
in $U_I$. The tangent vector space $T_I$ becomes identified with
$U_I$, cf. (1.4.6).
 Consider the action of ${\bf C}^*$ on $U_I$ by homotheties (simultaneous
multiplication of the coordinates $z_{ij}$ by a scalar). Since this action
 is a part of the torus action on $G(k,n)$, we deduce that for any Lie
complex $Z$ the intersection $Z\cap U_I$ will be a conic
(i.e. ${\bf C}^*$ -invariant) subvariety of $U_I$. The same will hold,
by continuity, for any generalized Lie
 complex. The map $f_I$, therefore, takes any generalized Lie complex
$Z$ (which
is a subvariety in $G(k,n)$) into the subvariety in $P(T_I)=P(U_I)$
represented
by the conic subvariety $Z\cap U_I$ in $U_I$. We need to show that $f_I$
 is a regular morphism of algebraic varieties.

Since both $G(k,n)//H$ and $P(U_I)//H$ are projective, it suffices,
 by Serre's GAGA theorem [24], to show that $f_I$ is a holomorphic
map of complex analytic spaces corresponding to these varieties.
However, this follows
from the description
 of the Chow varieties given by D.Barlet [3]. More specifically,
 Barlet gave a condition for a family of $p$-dimensional cycles
$Z(s)\i X$ parametrized by
 a reduced analytic space $S$, to be analytic near a point
 $s_0\in S$. This condition is essentially that for any
codimension $p$ analytic subvariety $Y\i X$ which intersects
$Z(s_0)$ properly, the 0-cycle $Z(s)\cap Y$ depends analytically
on $s$ near $s_0$. To prove that $f_I(Z)$ depends analytically
on $Z\in G(k,n)//H$ we note that any analytic subvariety $Y$ of
 codimension $(n-1)$  in $P(U_I)$ can be lifted to $U_I$
(by setting one of the coordinates to be 1) to a subvariety
$\tilde Y$. If $Y$ intersected some $f_I(Z)$ properly then
 so does $\tilde Y$ with respect to $Z$ due to the conic property.
Analytic dependence of $\tilde Y \cap Z$ on $Z$ implies the analytic
 dependance of $Y\cap f_I(Z)$ which is just
the image of $\tilde Y\cap Z$ in the projectivization. Proposition
(1.4.11) is proven.
 \vskip .2cm

Later we shall make use of
morphisms $f_I$  to construct "coordinate charts" on the Chow quotient.
The following fact is an immediate consequence of (1.4.11).

\vskip .3cm

\proclaim (1.4.14) Corollary. Each regular  polyhedral decomposition
(in particular, triangulation)  of the product of simplices
$\Delta^ I \times \Delta ^ {\bar I}$ has the form ${\cal P}_I$
 for some realizable matroid decomposition of $\Delta (k,n)$.

Thus the problem of classification of all realizable matroid
decompositions of hypersimplex contains the classification problem
 for triangulations of the product of two simplices.

\vskip 1cm
\centerline {\bf (1.5) Relation to the Hilbert quotient.}
\vskip .5cm

\noindent (1.5.1) The Hilbert quotients were defined in n. (0.5).
We want to compare the Chow quotient $G(k,n)//H$ with the Hilbert
 quotient $G(k,n)///H$. Recall  that
there is a  a regular birational
morphism $\pi : G(k,n)///H \rightarrow G(k,n)//H$ to the Chow quotient,
see (0.5.8).
\vskip .3cm

\proclaim (1.5.2) Theorem. The morphism $\pi$ is an isomorphism.

\noindent {\sl Proof:} By definition, points of $G(k,n)///H$ are
 subschemes which are limit positions of Lie complexes.
By Proposition 1.11 every
such subscheme $Z$ is reduced at a generic point of every its component.

\proclaim (1.5.3) Lemma. Any subscheme $Z$ from $G(k,n)///H$ is reduced.

\noindent {\sl Proof:} Consider the intersection of $Z$ with some coordinate
Schubert chart $U_I$. It suffices to prove, for every $I$, that
$Z\cap U_I$ is reduced. The action of $H$ in $U_I$ is a linear one.
This is the action corresponding to the products of simplices and it
follows from the result of B.Sturmfels ([48], Theorem 6.1)
that any limit position of generic $H$ -orbits in $U_I$ is reduced.

Lemma (1.5.3) implies that the morphism $\pi$ is bijective on
${\bf C}$ -points. To show that it is an isomorphism of algebraic
 varieties, it suffices to show that for
any $Z\in G(k,n)//H$ and any non-zero Zariski tangent vector
$\xi$ to $G(k,n)///H$ at $Z$  the vector $d\pi (\xi)$, the
image  of $\xi$ under the differential of $\pi$, is non-zero.
The tangent space $T_Z{\cal H}_G$ to the whole Hilbert scheme
at $Z$ equals $H^0(Z,{\cal N})$ where
 ${\cal N}_Z = \underline{Hom}_{{\cal O}_G}(J_Z/J_Z^2, {\cal O}_Z)$
 is the normal sheaf of $Z$ (see \S 0E).
Let $Z_{reg}$ be the smooth part of $Z$. The restriction of ${\cal N}_Z$
to $Z_{reg}$ is the sheaf of section of the normal bundle of $Z_{reg}$
in the usual sense.
 Hence our vector $\xi\i T_Z{\cal H}_G$ gives a normal vector field
on $Z_{reg}$.
Since $Z$ is reduced, this field is non-zero if $\xi$ is non-zero,
so the assertion follows. Theorem (1.5.2) is proven.

\vskip 1cm
\centerline {\bf (1.6)  The (hyper-) simplicial structure on the
collection of $G(k,n)//H$.}
\vskip .5cm

\noindent (1.6.1) Recall (1.1.2) that $G^0(k,n)$ denotes the
 generic stratum in the Grassmannian $G(k,n)$.
 For any $i\in \{1,...,n\}$ there are the {\it intersection
 and projection maps}
$$G^0(k,n-1)\buildrel B_i\over\longleftarrow G^0(k,n) \buildrel
 A_i\over\longrightarrow G^0(k-1, n-1)\leqno (1.6.2)$$
defined as follows. Let $J_i: {\bf C}^{n-1}\hookrightarrow {\bf C}^n$
 be the embedding taking $(x_1,...,x_{n-1})$ to $(x_1,...,x_{i-1},
0,x_i,...,x_n)$.
 The intersection map $A_i$ sends a $k$ -dimensional subspace
 $L\i {\bf C}^n$
to $J_i^{-1}L$. The projection map $B_i$ is induced by the
projection ${\bf C}^n \rightarrow {\bf C}^{n-1}$ forgetting
the $i$ -th coordinate.

The formal structure of these maps is analogous to that of
 faces of the
hypersimplex $\Delta (k,n)$ (Proposition 1.2.4). The existence
 of such a system of "face" maps was the original reason for
 introducing hypersimplexes
 and then Grassmannians into the problem of combinatorial calculation of
characteristic classes [17].
  More recently, these maps were used by  A.A.Beilinson, R.D.
MacPherson and V.V.Schectman in [7] to give a "constructuble"
approximation to K-theory.
\vskip .3cm

\noindent (1.6.3) As was noted in [7], the maps (1.6.2) descend
to maps of the quotients
$$G^0(k,n-1)/({\bf C}^*)^{n-1}\buildrel b_i\over\longleftarrow
G^0(k,n)/({\bf C}^*)^n \buildrel a_i\over\longrightarrow
G^0(k-1, n-1)/({\bf C}^*)^{n-1}\leqno (1.6.4)$$
of the generic strata by their respective tori.
\vskip .3cm

\noindent (1.6.5) Clearly there is no way to extend maps
(1.6.2) to whole Grassmannians: if the
subspace $L$ is contained in the hyperplane $\{x_i=0\}$
then $A_i (L)$ will have
wrong dimension and similarly for $B_i$. However, it turns
 out that for Chow quotients the situation is different.

\proclaim (1.6.6) Theorem. The maps $a_i,b_i$ in (1.6.4)
can be extended to regular morphisms of projective algebraic
 varieties
$$G(k,n-1)//({\bf C}^*)^{n-1}\buildrel \tilde b_i\over\longleftarrow
G(k,n)//({\bf C}^*)^n \buildrel \tilde a_i\over\longrightarrow G(k-1,
 n-1)//({\bf C}^*)^{n-1}\leqno (1.6.7)$$

Proof of Theorem (1.6.6) will occupy the rest of this section.
\vskip .3cm

\noindent (1.6.8)
Let $e_1,...,e_n\in {\bf C}^n$ be the standard basis vectors and
${\bf C}^{n-1}_i$ be the coordinate hyperplane spanned by $e_j, j\neq i$.
For any $i\i \{1,...,n\}$ we consider the varieties
$$G^+_i = \{L\in G(k,n): e_i\in L\},\,\,\, G^-_i = \{L\in G(k,n):
 L\i {\bf C}^{n-1}_i\}.\leqno (1.6.9)$$
They are the analogs of the family of coordinate hyperplanes in $P^{n-1}$.
The next proposition is immediate.
\vskip .3cm

\proclaim (1.6.10) Proposition. a) As abstract varieties,
 $G_i^+$ are isomorphic to
 $G(k-1,n-1)$ and $G_i^-$ - to $G(k,n-1)$.\hfill\break
b) Both $G_i^+$ and $G_i^-$ are linear sections of $G(k,n)$
 in the Pl\"ucker embedding. More precisely, let $\Pi_i^+$
and  $\Pi_i^-$  be projective subspaces in $P(\bigwedge ^k{\bf C}^n)$
 given by vanishing of the Pl\"ucker coordinates $p_I, i\notin I$
 or, respectively, $p_I, i\in I$. Then $G_i^\pm = G(k,n)\cap \Pi_i^\pm$.
 \hfill\break
c) The image of the subvarieties $G_i^\pm$ under the moment map $\mu$ from
(1.1) is the facet $\Gamma_i^\pm$ of $\Delta (k,n)$ introduced in
(1.2.5). Moreover, we have $G_i^\pm = \mu^{-1}(\Gamma_i^\pm)$.
\vskip .3cm

\noindent (1.6.11)  The isomorphisms in part a) of   (1.6.10) are
as follows. the isomorphism
 $u_i:G(k-1,n-1)\rightarrow G_i^+$   takes a
$(k-1)$ -dimensional
 $\Lambda\i {\bf C}^{n-1}$ to $J_i(\Lambda)\oplus {\bf C}e_i$, where the
embedding $J_i : {\bf C}^{n-1}\hookrightarrow {\bf C}^n$ was defined
 in (1.6.1). The isomorphism $v_i:G(k,n-1)\rightarrow G_i^-$ takes a
 $k$ -dimensional
$M\i {\bf C}^{n-1}$ into $J_i(M)$.
\vskip .3cm

\noindent (1.6.12) Let us turn to the construction of
$\tilde a_i, \tilde b_i$. Note that the
maps $a_i, b_i$ of the quotients of generic strata have the following
 transparent description in terms of Lie complexes (closures
 of generic orbits).
\vskip .3cm

\proclaim (1.6.13) Lemma.  Let $Z=\overline {H.L}$ be a Lie complex
 in $G(k,n)$.
 Then the Lie complex $a_i(Z)$ in $G(k-1,n-1)$ representing the orbit of
 $A_i(L)$, is equal to $u_i^{-1} (Z\cap G_i^+)= u_i^{-1}(Z\cap \Pi_i^+)$.
 Similarly, the Lie complex
$b_i(Z)$ in $G(k,n-1)$, representing the orbit of $B_i(L)$, is equal to
$v_i^{-1} (Z\cap G_i^-) = v_i^{-1}(Z\cap \Pi_i^-)$.\hfill\break

\noindent {\sl Proof:}  Denote by $g_i(t)$ the diagonal matrix
 $(1,...,1,t,1,...,1)\in ({\bf C}^*)^n$ (where $t$ is on the $i$ -
th place). Let $pr_i:{\bf C}^n\rightarrow {\bf C}^{n-1}_i$ be the
coordinate projection.
Let $L\in G^0(k,n)$ be a generic subspace. Then we have
$$(L\cap {\bf C}^{n-1}_i)\oplus {\bf C}e_i = \lim_{t\rightarrow \infty}
g_i(t).L,\,\, \quad pr_i(L)  = \lim_{t\rightarrow 0} g_i(t).L.\leqno
 (1.6.14)$$
This shows that $Z\cap G_i^+$ (resp. $Z\cap G_i^-$) contains
 the orbit of $A_i(L)$ (resp. $B_i(L)$).

On the other hand, the vanishing of Pl\"ucker coordinates $p_I$,
 for which $e_I\notin \Gamma_i^\pm$, forms a system of equations
 for $Z\cap G_i^+$, as follows from the general theory of toric
varieties [49]. Lemma (1.6.13) is proven.

The lemma just proven shows that $a_i$ and $b_i$ are given by
intersecting Lie complexes with projective subspaces $\Pi_i^\pm$.
However, $dim (Z\cap \Pi_i^\pm) = dim (Z) -1$ whereas $\Pi_i^\pm$
 have high codimension. Thus to prove that the intersection gives
a regular morphism of Chow quotients, extra work is needed.
\vskip .3cm

\noindent (1.6.15) \underbar{Proof that $\tilde a_i, \tilde b_i$
 are regular morphisms.} For any face $\Gamma\i \Delta (k,n)$ and
 any Lie complex $Z$ we shall denote by $Z(\Gamma)$ the closure of
 the $H$ -orbit in $Z$ corresponding to $\Gamma$. In particular, the
 codimension 1 faces of $\Delta (k,n)$ are $\Gamma_i^\pm$ from
Proposition 1.8 and the corresponding orbit closures $Z(\Gamma_i^\pm)
 = Z\cap G_i^\pm = Z\cap \Pi_i^\pm$ are the varieties we are studying
. Our aim is to show that the Chow form of $Z(\Gamma_i^\pm)$ can be
 polynomially expressed via that of $Z$.

\proclaim (1.6.15.1) Lemma. Consider any coordinate hyperplane
$\{p_I=0\}$ in $P(\bigwedge ^k {\bf C}^n)$ given by vanishing of
 a Pl\"ucker coordinate $p_I$. Then we have the equality of cycles
$$Z\cap \{p_I=0\} = \sum_{\Gamma: e_I\notin \Gamma} Z(\Gamma),$$
where $\Gamma$ runs over codimension 1 faces of $\Delta(k,n)$
 (i.e. over $\Gamma_i^\pm)$ not containing the vertex $e_I$.

\noindent {\sl Proof:} The lemma says that the order of vanishing
 of $p_I$ on $Z(\Gamma)$ equals 0 if $e_I\in\Gamma$ and 1 if
 $e_I \notin \Gamma$. According to the general rule (valid for
 any toric variety in an equivariant projective embedding [49])
 this order equals the distance from $e_I$ to the affine
hyperplane spanned by $\Gamma$, the distance being measured
 in natural integer units indiced by the lattice. In our case,
 if $e_I\notin \Gamma$, then the said distance equals one.

\proclaim (1.6.15.2) Corollary. Let $R_Z(l_1,...,l_n), l_i\in
\bigwedge^k({\bf C}^n)^*$, be the Chow form of a Lie complex $Z$.
 Let $\pi_I$ be the coordinate projection to the coordinate
hyperplane $Ker (p_I) \i \bigwedge^k({\bf C}^n)$. Then for any
 linear functionals $\lambda_1,...,\lambda_{n-1}\in (Ker p_I)^*$ we have
$$R_Z(p_I, \pi_I^*\lambda_1,...,\pi_I^*\lambda_{n-1}) =
\prod _{\Gamma: e_I\notin\Gamma} R_{Z(\Gamma)} (\lambda_1,...,\lambda_{n-1})$$
where on the right hamd side stand Chow forms of subvarieties
 $Z(\Gamma) \i P(Ker p_I)$.
\vskip .3cm

\noindent (1.6.15.3) \underbar {End of the proof that
$\tilde a_i, \tilde b_i$ are regular}.  Consider some
facet of $\Delta(k,n)$, say, $\Gamma_i^+$  and let
 $e_I\notin \Gamma_i^+$ be some vertex (i.e., $i\in I$).
 Consider the coordinate projection
$Ker p_I \rightarrow \Pi_i^+$
of coordinate subspaces in $\bigwedge^k({\bf C}^n)$.
Here $\Pi_i^+$ is defined, as above, by vanishing of all
$p_I$ with $i\in I$. The projection of $\bigcup_{\Gamma:e_i\notin \Gamma}
 Z(\Gamma)$ to $\Pi_i^+$ is the component $Z(\Gamma_i^+$ we are
 interested in, plus the union of some coordinate $(n-2)$ -
dimensional subspaces (which are images of other components).
 Let $\mu_1,...,\mu_{n-1}\in (\Pi_i^+)^*$ be linear forms.
Denote by $\tilde\mu_i$ the extension of $\mu_i$ to all of
$Ker(p_I)$ by means of the coordinate projection $Ker p_I
\rightarrow \Pi_i^+$. We obtain the equality
$$\prod_{\Gamma:e_I\notin\Gamma} R_{Z(\Gamma)}(\tilde \mu_1,...,
\tilde\mu_{n-1}) = R_{Z(\Gamma_i^+)}(\mu_1,...,\mu_{m-1}) . \
prod_S R_S(\mu_1,...,\mu_{n-1})^{\nu_S},$$
where $S$ runs over coordinate $(k-1)$ -dimensional subspaces
 in $\pi_i^+$, $R_S$ is the Chow form of the subspace $S$ and
$\nu_S$ are some exponents. From this equality we obtain that
$R_{Z(\Gamma_i^+)}$ can be obtained from the right hand side by
 division to a fixed polynomial. Since the left hand side itself
depends polynomially on $R_Z$, Theorem (1.6.6) is completely proven.

\hfill\vfill\eject

\centerline {\bf Chapter 2. PROJECTIVE CONFIGURATIONS AND THE}
\centerline {\bf GELFAND-MACPHERSON ISOMORPHISM.}

\vskip 1cm

I.M.Gelfand and R.W.MacPherson have established in [20] an important
 correspondence between torus orbits in $G(k,n)$ and projective configurations
 i.e.   $GL(k)$- orbits on  $(P^{k-1})^n$.
 In this section shall show that this correspondence extends to an
isomorphism of Chow quotients.

\vskip 1cm

\centerline{\bf (2.1) Projective configurations and their Chow quotient.}
\vskip .5cm

\noindent (2.1.1) Consider the $(k-1)$ -dimensional projective space $P^{k-1}$.
By a {\it configuration} we shall mean an ordered collection
$M=(x_1,...,x_n)$ of $n$ points in $P^{k-1}$.
The general linear group $GL(k)$ acts on $P^{k-1}$ by projective
transformations. This induces an action on the space $(P^{k-1})^n$ of
configurations.

 The study of orbits of this action is a classical problem of
projective geometry. See [14] for investigations from the standpoint
 of Mumford's
geometric invariant theory.
\vskip .3cm

\noindent (2.1.2) The elements of a configuration $M$ can also be
visualized as hyperplanes
(in the dual projective space). This point of view will be useful
later. In this subsection we shall just consider elements of $M$ as points.

\vskip .3cm

\noindent (2.1.3) We will be interested in the Chow quotient
$(P^{k-1})^n//GL(k)$. To apply Definition 0.1.7 it is first
desirable to know which configurations are "generic enough".
The answer, of course, is the following.
\vskip .3cm

\noindent (2.1.4) A configuration
$ M = (x_1,...,x_n)$ of points in $P^{k-1}$ will be said to be
{\it in general position} if any $i$ of these points, $i\leq k$,
span a projective subspace of dimension exactly $i-1$. The set of
all such configurations will be denoted by
$(P^{k-1})^n_{gen}$. Orbits of configurations in general position
 will be referred to as {\it generic orbits} in $(P^{k-1})^n$.
\vskip .3cm

\noindent (2.1.5) Generic  $GL(k)$ - orbits  on $(P^{k-1})^n$ depend
on continuous
parameters only when $n\geq k+2$. We shall assume in the sequel that
 this condition holds. In this case generic orbits have dimension $k^2-1$
since
the stabilizer of a generic configuration consists only of homotheties.
\vskip .3cm

 \noindent (2.1.6) For any $0\leq m\leq k-1$ denote by $[m]$ the $2m$ -
dimensional homology class of $P^{k-1}$ represented by $P^m$. By K\"unneth
formula, the graded homology space of $(P^{k-1})^n$ is the $n$ -fold
tensor power of the graded homology space
$H_*(P^{k-1})$. Therefore,
 the basis for the $2p$ -th homology group $H_{2p}( (P^{k-1})^n )$
 is given by tensor products
$[m_1]\otimes ...\otimes [m_n], \sum m_i =p$.
\vskip .3cm

\proclaim (2.1.7) Proposition. The homology class of the closure of
 any generic
$GL(k)$ -orbit in $(P^{k-1})^n$, $n\geq k+2$, is a variety of dimension
 $k^2-1$ and of homology class
$$\delta =\sum_{m_1+...+m_n=k^2-1\atop m_i\leq k-1} [m_1]\otimes ...
\otimes [m_n].$$
The set of closures of generic orbits is a subvariety in the Chow variety
${\cal C}_{k^2-1}((P^{k-1})^n,\delta)$ isomorphic to the quotient
$(P^{k-1})^n_{gen}/GL(k)$.

\noindent {\sl Proof:} Let $Z= \overline{GL(k).M}$ be the closure
of any $k^2-1$ - dimensional orbit  and $\delta\in H_{2(k^2-1)}
( (P^{k-1})^n,{\bf Z})$ -its homology class. The coefficient in
$\delta$ at $[m_1]\otimes ...\otimes [m_n]$ can be calculated as
 follows. Take generic projective subspaces $L_i\i P^{k-1}$ of
codimension $m_i$.
Our coefficient is just the intersection number of $Z$ with
 $L_1\times ...\times L_n$. In other words, this is the number
 of projective transformations which take each point $x_i$ of
 our configuration $M = (x_1,...,x_n)$ inside $L_i$. The condition
$g(x_i)\i L_i$  is a linear condition on matrix element of a matrix
$g\in GL(k)$ of codimension $m_i$. Taking into account all $L_i$, we
 obtain a system of $k^2-1$ linear equations on matrix elements of $g$.
 By Bertini's theorem applied to $Z$, if $(L_1,...,L_n)$ are generic
 enough, the intersection $Z\cap (L_1\times ...\times L_n)$ consists
 of finitely many points. For our linear system this implies that
 for generic $L_j$ just
one of two cases holds:
\item{(2.1.7.1)} The space  of solutions of the system is 1-dimensional
and consists of multiples of a non-degenerate matrix.
\item{(2.1.7.2)} The space of solutions is contained in the variety of
 degenerate matrices.

\noindent In the first case the coefficient equals one, in the second
 it equals 0.  We need to  show that for $x_1,...,x_n$ in general position,
 the case a) always holds.
\vskip .3cm

\noindent (2.1.7.3) Consider the product of Grassmannians $\Pi =
\prod G(k-m_i,k)$ i.e. the variety of all tuples $(L_1,...,L_n)$
as above. Let $\Pi_Y\i \Pi$ be the subvariety of those tuples for
which $x_i\i L_i$ for all $i$.
\vskip .3cm

\proclaim (2.1.7.4) Lemma. Let $Y=(x_1,...,x_n)$ be any configuration
with $k^2-1$ - dimensional orbit. Then the truth of the case (2.1.7.2)
 above (or, equivalently, the vanishing of the coefficient at
$[m_1]\otimes ...\otimes [m_n]$) is equivalent to the following fact:
For generic tuple of subspaces $(L_1,...,L_n)\i \Pi_Y$ its stabilizer
in $PGL(k-1)$ has positive dimension.

\noindent {\sl Proof:} Case (2.1.7.2) means that the union of $GL(k)$
- orbits of points from $\Pi_Y$ is not dense in $\Pi$.  The codimension
of $\Pi_Y$ in $\Pi$ equals $k^2-1$. Therefore case (2.1.7.2) means that
 for any $Y\in \Pi_Y$ its
 orbit has dimension smaller than $k^2-1$.

Now it is clear that if $x_1,...,x_n$ are in general position then for
any $m_1,...,m_n$ summing to $k^2-1$ it is possible to choose codimension
 $m_i$ subspace $L_i$ through $x_i$ such that the whole collection
$(L_1,...,L_n)$ has trivial stabilizer in $PGL(k)$. Proposition 2.1.7
 is proven.

\noindent {\bf (2.1.8) Example.} Consider the action of $GL(2)$ in
$(P^1)^4$. The closure of the orbit of a 4-tuple of distinct points
is a 3-dimenisional variety. It contains four 2-dimenional orbits
 $W_i$ where $W_i$ is the set of points $(x_1,...,x_4)$
 such that all three $x_j$ with $j\neq i$, coincide with each other
 but differ from $x_i$. The closure of each $W_i$ is isomorphic to
 $P^1\times P^1$. These closure intersect along
 the 1-dimensional orbit which is the set of coinciding tuples.
\vskip .3cm

\noindent (2.1.9) We are interested in the Chow quotient
 $(P^{k-1})^n//GL(k)$. By definitions, points of this quotient
are certain algebraic cycles in $(P^{k-1})^n$ of dimension $k^2-1$
 and homology class $\delta$ given by
 Proposition 2.1.7. Moreover, since the stabilizer of a configuration
can not be a unipotent subgroup in $PGL(k)$, we can apply Theorem 0.3.1
to conclude that  components of any cycle from $(P^{k-1})^n//GL(k)$ are
 closures of $k^2-1$ -dimensional orbits.
\vskip .3cm

\noindent {\bf (2.1.9) Example.}  For 4 distinct points on $P^1$ the only
invariant is the cross-ratio which identifies $(P^1)^4_{gen}/GL(2)$
 with $P^1-\{0,1,\infty\}$. Denote by $Z_\lambda$ the closure of the
 orbit given by 4-tuples with cross-ratio $\lambda$. When
$\lambda\rightarrow 0,1,\infty$, the variety $Z_\lambda$ degenerates
into one of three cycles in $(P^1)^4$. Namely, let $\Delta_{ij}$ to
 be the subset in $(P^1)^4$ given by $\{x_i=x_j\}$. Then the three
cycles in question are
$$\Delta_{12}+\Delta_{34}, \Delta_{13}+\Delta_{24},
 \Delta_{14}+\Delta_{23}.$$
For example, suppose that our four points $x_i$ depends on a
parameter $t$ and degenerate in such a way that $x_1(0)=x_2(0)$
 but $x_3(0)$ and $x_4(0)$ are different from them (see Fig. 2.1.10).
 Let $Z(t), t\neq 0$ be the closure of the orbit of $M(t) =
( (x_i(t)_{i=1,...,4} )$.
Then, of course, the orbit of the limit position $( (x_i(0) )$,
 i.e. $Z_{12}$, will be a part of the cycle $Z(0)=lim_{t\rightarrow 0}Z(t)$,
but not the only part! Indeed, we can perform, for each $t$,
 a projective transformation $g(t)$ which stretches $x_1(t)$ and
$x_2(t)$ back to some fixed distance. This transformation shrinks
 the remaining points $x_3(t)$ and $x_4(t)$ close to each other.
The
limit of the point $g(t)Y$ will lie on the second component $Z_{34}$.

\vbox to 4cm{\vfill\vfill (2.1.10) \vfill\vfill}
\vskip .3cm

\noindent (2.1.11) Similarly,  if we have a degeneration
$M(t)= (x_1(t),...,x_n(t))$ of a family of $n$ points on
$P^1$ such that just two points merge, e.g. $x_1(0)=x_2(0)$
 and all the other $x_i(0)$ remain distinct, then
$\lim_{t\rightarrow 0}\overline {GL(2)M(t)}$ will consist
of two components. The first is the orbit of the limit
configuration $(x_1(0)=x_2(0),x_3(0),...,x_n(0)$. The second
 component is the
set of $(x_1,...,x_n)$ such that $x_3=...=x_n$ and $x_1,x_2$ are
arbitrary.
We shall see later (\S 4) that this phenomenon exactly corresponds
 to the degeneration of $(P^1,x_1(t),...,x_n(t))$ in the Knudsen's
 moduli space $\overline{M_{0,n}}$ of stable $n$ -punctured curves
of genus 0.

\vskip 1cm

\centerline {\bf (2.2) The Gelfand-MacPherson correspondence.}
\vskip .5cm

\noindent (2.2.0) Let us recall the original idea of [20,40] how
 to construct a
configuration from a point in Grassmannian. It will be more
convenient for
us to speak in this section about configurations of {\it hyperplanes}
instead of points.
\vskip .3cm

\noindent (2.2.1) Let $L\i {\bf C}^n$ be a $k$-dimensional subspace
 not lying in any
coordinate hyperplane $H_i =\{x_i=0\}$. Then $(L\cap H_i)$ form a
 configuration of hyperplanes in $L$ i.e. a point in $(P(L^*))^n$.
If a subspace $L'$ is obtained from $L$ by the action of a torus element,
 we shall obtain a projectively isomorphic configuration of hyperplanes
 in $L'$. A class of projective isomorphism of configurations of $n$
 hyperplanes in $(k-1)$- dimensional projective spaces is the same as
 a $GL(k)$- orbit in  the Cartesian power of a fixed projective space
 $(P^{k-1})^n$. Note that not every configuration of hyperplanes can be
obtained, up to an isomorphism, from $L$ as above. To make the assertion
precise, denote by $G_{max}(k,n)\i G(k,n)$ the set of $L$ such that
$dim (H.L)= n-1$. Similarly denote by
$((P^{k-1})^n)_{max}$ the set of configurations $\Pi = (\Pi_1,...,\Pi_n)$
 such that $dim(GL(k).\Pi) = k^2-1$. The Gelfand- MacPherson correspondence
induces the bijection of orbit sets
$$\Lambda :G_{max}(k,n)/H \rightarrow ((P^{k-1})^n)_{max}/GL(k)
 \leqno (2.2.2)$$
see [20].
\vskip .3cm

\noindent (2.2.3) Note that sets in both sides of (2.2.2) are not,
 in general algebraic varieties
since $G_{max}(k,n)$ and $(P^{k-1})^n)_{max}$ contain unstable points.
For comparison of Mumford's quorients of both sides in (2.2.2) see
 section (2.4) below.

The main result of this section is the following theorem

\proclaim (2.2.4) Theorem. The Gelfand-MacPherson correspondence
 (2.2.2) extends to   an isomorphism of Chow quotients
$$\Lambda :G(k,n)//H\longrightarrow (P^{k-1})^n//GL(k).\leqno (2.2.5)$$

This fact permits one to apply the information about behaviour of
$(n-1)$ - dimensional torus orbits (which may be obtained by
 techniques of toric varieties and A- resultants [30]), to the study
 of $(k^2-1)$ dimensional orbits of $GL(k)$ which are at first
glance harder to understand.
\vskip .3cm

\proclaim (2.2.6) Corollary. Every cycle in $(P^{k-1})^n//GL(k)$
is a sum of closures of some $(k^2-1)$ -dimensional orbits with
 multiplicities 0 or 1.

Before starting to prove  Theorem 2.2.5, let us give a simple
 matrix interpretation of the correspondence in question.

\vskip .3cm

\noindent (2.2.7) Let $M(k,n)$ be the vector space of all complex
 $k$ by $n$ matrices, $M_0(k,n)\i M(k,n)$ - the space of matrices
 of rank $k$, and $M'(k,n)\i M(k,n)$- the space of matrices whose
every row is a non-zero vector in ${\bf C}^k$. The group $GL(k)$
acts on $M(k,n)$ form the right, and $({\bf C}^*)^n\i GL(n)$ -
from the left and we have the identifications
$$GL(k)\backslash M_0(k,n) = G(k,n),\,\,\,\,\, M'(k,n)/({\bf C}^*)^n =
 (P^{k-1})^n.\leqno (2.2.8)$$
The Gelfand- MacPherson correspondence comes from consideration of
both types of orbits in (2.2.5) as double $(GL(k), ({\bf C}^*)^n)$
- orbits in $Mat (k,n)$.
\vskip .3cm

\noindent (2.2.9) Let us carry on these considerations for Chow quotients.
 Note that each \hfill\break $(GL(k),  ({\bf C}^*)^n)$ -orbit in the
 vector space $M(k,n)$ is invariant under multiplications by scalars
 (in this vector space) and thus may be identified with a subvariety
 in the projectivization $P(M (k,n))$. Instead of double orbits we
 can speak about left orbits of the product $GL(k)\times
 ({\bf C}^*)^n$. Consider the Chow quotient $P(M(k,n))//GL(k)
\times ({\bf C}^*)^n$. To prove Theorem 2.2.4 it suffices to
 construct isomorphisms

$$G(k,n)//({\bf C}^*)^n\buildrel \alpha\over\rightarrow P(M(k,n))
//GL(k)\times ({\bf C}^*)^n \buildrel \beta\over\leftarrow
(P^{k-1})^n//GL(k).\leqno (2.2.10)$$

\vskip .3cm

\noindent (2.2.11) The existence of these morphisms does not
 present any problem.

The morphism $\alpha$ associates to any cycle $Z=\sum c_iZ_i$
in $G(k,n)//({\bf C}^*)^n$ the cycle $\sum \overline {p^{-1}(Z_i)}$
where $p:P(M_0(k,n))\rightarrow G(k,n)$ is the projection from (2.2.8)
(The multiplicities $c_i$ all are equal to 1 by Proposition 1.2.15).
 Similarly, the morphism $\beta$ associates to any cycle
$W= \sum m_iW_i$ in $(P^{k-1})^n//GL(k)$ the cycle
$\sum m_i\overline {q^{-1}(W_i)}$, where $q: P(M'(k,n))
\rightarrow (P^{k-1})^n$ is the other projection arizing
from (2.3). To show that $\alpha$ and $\beta$ thus defined are regular
maps, it suffices to apply Barlet's criterion of analytic
 dependence of a cycle on a parameter [3]. Since $\alpha$ and
$\beta$ are both given by inverse images in fibrations,
this criterion is trivially applicable.
\vskip .3cm

\noindent (2.2.12) Let us show that $\alpha$ is an isomorphism.
 To do this, note that any
generic $GL(k)\times ({\bf C}^*)^n$- orbit in $P(M(k,n)$ has dimension
 $(n-1)(k^2-1)$. Each component of a cycle $Z$ from
 $P(M(k,n))//GL(k)\times ({\bf C}^*)^n $ is the closure of a single orbit
 which therefore should be the inverse image of an orbit of maximal dimension
in $G(k,n)$. The algebraic cycle $W$ formed by these orbits lies clearly in
 the Chow quotient $G(k,n)//({\bf C}^*)^n$ and this is the unique element of
this Chow quotient such that $\alpha (W) =Z$. This proves that $\alpha$
is bijective on {\bf C}-points. Denote by $\alpha ^{-1}$ the inverse map.
To prove that $\alpha$ is indeed an isomorphism of
algebraic varieties we need to prove that $\alpha^{-1}$ is regular too (which
need not necessarily be the case if the varieties involved are not normal).
However, this again follows from Barlet's criterion similarly to the proof
 of Theorem 1.4.12.

Similarly we prove that $\beta$ is an isomorphism. Theorem 2.2.4 is proven.

\vskip 1cm
{\bf (2.3) Duality (or association).}
\vskip .5cm

It is known classically (since A.B.Coble [11]) that projective
 equivalence classes of
configurations of $n$ ordered points in $P^{k-1}$ are in bijection
with projective equivalence classes of configurations of $n$ points
in $P^{n-k-1}$. This correspondence is known as the {\it association}
 [11,14]  and was used in the context of matroid theory (see [21], \S 2..3).

 The most transparent way to define the association is via
the Gelfand-MacPherson
correspondence.

\vskip .3cm

\noindent (2.3.1) Let us
identify the dual subspace to the coordinate space ${\bf C}^n$
 with ${\bf C}^n$
by means of the standard pairing. By considering orthogonal complements
 to $k$ -dimensional subspaces we obtain an isomorphism
 $G(k,n)\cong G(n-k,n)$.
The torus $H=({\bf C}^*)^n$ acts in both Grassmannians and the
 said isomorphism is $H$ -equivariant. Hence it induces the
isomorphism of coset spaces
$ G(k,n)/H \rightarrow G(n-k,n)/H$.
Taking into account the Gelfand - MacPherson isomorphism (2.2.2),
we obtain the following isomorphism
$$A_{k,n}: (P^{k-1})^n_{max}/GL(k) \rightarrow
 (P^{n-k-1})^n_{max}/GL(n-k) \leqno (2.3.2)$$
where the subscript $"max"$ means the set of points whose
orbits have the maximal dimension.
The isomorphism (2.3.2) will be called the {\it association
 isomorphism}. By construction, this system of isomorphisms
 is involutive i.e. $A_{k,n}\circ A_{n-k,n} = Id$.
\vskip .3cm

\noindent (2.3.3) If $(x_1,...,x_n), (y_1,...,y_n)$ are $n$ -
tuples of points in $P^{k-1}$ and $P^{n-k-1}$ respectively then
we shall say that $(y_i)$ is associated to $(x_i)$ (and vice versa)
 if their orbits under projective transformations have maximal
dimensions and are taken into each other by the association isomorphism.
\vskip .3cm

\noindent (2.3.4) Explicitly the  configuration associated to
 $(x_i)$ can be calculated as follows. Let ${\bf x}_i\in {\bf C}^k$
 be vectors whose projectivizations are $x_i$.
By definition, we have to find a $k$ -dimensional subspace
$L\i {\bf C}^n$ and an isomorphism ${\bf C}^k \rightarrow L^*$
 which takes ${\bf x}_i$ into the
restriction of the $i$ -th coordinate function to $L$. Then the
 vectors
${\bf y}_i\in {\bf C}^n/L$ defined as the projections of the
 standard basis vectors, will represent the associated configuration.
 In other words, we have to find a complete $(n-k)$ -dimensional
 system of linear relations between ${\bf x}_i$, namely
$y_{j1}{\bf x}_1+...+y_{jn}{\bf x}_n=0, j=1,...,n-k$.
 Vectors ${\bf y}_i$ representing the associated configuration
 $(y_i)$ are given by columns of the matrix $||y_{ji}||$.
This gives the following criterion.
\vskip .3cm

\proclaim (2.3.5) Proposition. Let ${\bf x}_i \in {\bf C}^k,
 {\bf y}_i\in {\bf C}^{n-k}$ are $n$ -tuples of vectors such
 that the corresponding configurations of points $x_i\in P^{k-1},
 y_i\in P^{n-k-1}$ have orbits of maximal dimension. Then
$(y_i)$ is associated to $(x_i)$ if and only if there is a unique,
 up to constant, linear relation in ${\bf C}^k\otimes {\bf C}^{n-k}$:
$$\sum_i \lambda_i \,\,({\bf x}_i \otimes {\bf y}_i) = 0$$
which is such that all $\lambda_i \neq 0$.

This can be reformulated as follows.

\vskip .3cm

\noindent (2.3.6) Let $P$ be some projective space and $C\i P$
be some finite subset.
We shall say that $C$ is a {\it circuit} (in the sense of
 matroid theory, see [21]) if $C$ is projectively dependent
 but any its proper subset is projectively independent.
\vskip .3cm

\proclaim (2.3.8) Reformulation. Let $x=(x_1,...,x_n)\in (P^{k-1})^n$
 and $y=(y_1,...,y_n)\in (P^{n-k-1})^n$ be two $n$ -tuples whose orbits
 with respect to projective transformations have maximal dimensions
 (i.e. $k^2-1$ and $(n-k)^2-1$). Consider the Segre embedding
$P^{k-1}\times P^{n-k-1} \hookrightarrow P^{k(n-k)-1}$. Then
 $y$ is associated to $x$
 if and only if the points $(x_i,y_i)\in P^{k(n-k)-1}$  form a circuit.

"Normally" one would expect that points $(x_i,y_i)$ are
projectively independent.
\vskip .3cm

\noindent (2.3.9) For the case of $2k$ points in $P^{k-1}$
the source and the target of the association isomorphism are
 the same, so it is possible to speak about a configuration
being self-associated.
The following characterization of self-associated configurations
due to A.Coble [11] is a corollary of Reformulation 2.3.8.

\proclaim (2.3.10) Corollary. Let $x=(x_1,...,x_{2k})\in
(P^{k-1})^{2k}$ be a $2k$ -tuple. Consider the 2-fold
 Veronese embedding $P^{k-1}\hookrightarrow P(S^2{\bf C}^k)$.
 Then the configuration $x$ is self-associated if and only if
the images of of $x_i$ in the Veronese embedding form a circuit.

\noindent {\bf (2.3.11) Example.}  Let $k=2$. The association
induces an isomorphism between the set of projective equivalence
 classes of $n$ -tuples of distinct points on $P^1$ and the set
of projective equivalence classes of $n$ -tuples of distinct points
in $P^{n-3}$. This correspondence can be seen geometrically as
 follows ([14], Ch.III, \S 2, Proposition 2).

 Given $n$ points $y_1,..., y_n \in P^{n-2}$ in general position,
there is a unique rational normal curve (Veronese curve, for short)
 in $P^{n-3}$ through these points. This curve is isomorphic to $P^1$
and hence $y_i$ represent  on it a configuration of $n$ points
in $P^1$, which is the associated configuration to that of $y_i$
 in $P^{n-3}$. Conversely, given $n$ distinct points on $P^1$,
 we consider the $n-3$ -fold Veronese embedding of $P^1$. It
 identifies $P^1$ with a Veronese curve in $P^{n-3}$. The images
 $y_i$ of $x_i$ in this embedding are in general position as it
 may be seen by calculating the Vandermonde determinant. These
 points represent the configuration associated to that of $x_i$
 on $P^1$.
\vskip .3cm

\noindent {\bf (2.3.12) Example.} Let $x_1,...,x_6$ be a
configuration of 6 points in $P^2$ in general position.
 Corollary 2.3.10 means in this case that the configuration
$(x_i)$ is self associated if and only if the six points $x_i$
lie on a conic. Further examples can be found in [11,14].

\vskip .3cm

\noindent (2.3.13) Theorem 2.2.4 implies that the association
isomorphism extends to the Chow quotients of the spaces of
projective configurations. In other words, we have the following fact.

\proclaim (2.3.14) Corollary. There is an isomorphism of Chow quotients
 $$(P^{k-1})^n//GL(k) \rightarrow (P^{n-k-1})^n//GL(n-k).$$

\vskip 1cm
\centerline {\bf (2.4) Gelfand-MacPherson correspondence
and Mumford's quotients.}
\vskip .5cm

 For completeness sake we include here the comparison of Mumford's quotients
of $G(k,n)$ modulo torus and of $P^{k-1}$ modulo projective transformations.
\vskip .3cm

\noindent (2.4.1) First of all, the theory of Mumford is sensitive not only
 to the structure of orbits but also to the choice of group generating
 these orbits. In order that things behave well, we should consider the
 subgroup $H_1 = \{ (t_1,...,t_n)\in ({\bf C}^*)^n: \prod t_i = 1\}$
acting on $G(k,n)$ and the subgroup $SL(k)\i GL(k)$ acting on $(P^{k-1})^n$.
\vskip .3cm

\noindent (2.4.2) Recall (see n.0.4) that to define  Mumford's
quotient by  any group $G$ acting on any variety $X$ we should
 fix two things: an ample line bundle ${\cal L}$ on $X$ and a
linearization i.e. an extension $\alpha$ of $G$- action to ${\cal L}$.
\vskip .3cm

\noindent (2.4.3) First consider the $H_1$ -action on the
 Grassmannian $G(k,n)$. The Picard group of $G(k,n)$ is generated
by the sheaf ${\cal O}(1)$ in the Pl\"ucker embedding.
so there is essentially no freedom in choosing ${\cal L}$.
We set ${\cal L} = {\cal O}(1)$. For this choice of ${\cal L}$
 a linearization is given by an integral vector $a=(a_1,...,a_n)$
 defined modulo multiples of $(1,...,1)$. Denote
$t^a = t_1^{a_1}...t_n^{a_n}$ the character of $H_1$
corresponding to $a$.  The $H_1$ -action on ${\bf C}^n$
 corresponding to $a$, has the form
$$(t_1,...,t_n) \mapsto diag (t^a.t_1,...,t^a.t_n).$$
This action induces an $H_1$ -action on $\bigwedge ^k {\bf C}^n$
 which is the linearization corresponding to $a$.
\vskip .3cm

\noindent (2.4.4) Denote by $A(k,n)$ the coordinate ring of $G(k,n)$
in the Pl\"ucker embedding.
It is well-known [14] that $A(k,n)$ can be identified with the ring
of polynomials $\Phi(M)$ in entries of an indeterminate $(k\times n)$ -
matrix $M = ||v_{ij}||$ which satisfy the condition
$\Phi (gM) = \Phi(M)$ for any $g\in SL(k)$. In particular, the Pl\"ucker
coordinate $p_I$ corresponds to the polynomial in $v_{ij}$ given by
 the $(k\times k)$ - minor of $M$ on columns from $I$.

The Mumford quotient $(G(k,n)/H_1)_{{\cal O}(1),a}$ is, by definition,
 the projective spectrum $Proj (A(k,n)^{H_1})$ of the invariant
 subring in $A(k,n)$.
\vskip .3cm

\noindent (2.4.5) Consider now the $SL(k)$ -action on $(P^{k-1})^n$.
 For any integral vector $a = (a_1,...,a_k)$ denote by
 ${\cal O}(a) = {\cal O}(a_1,...,a_n)$ the line bundle
 on $(P^{k-1})^n$ whose local sections are functions multihomogeneous
 of degrees $(a_1,...,a_k)$. It is well-known [24,25] that bundles
${\cal O}(a)$ exhaust the Picard group of $(P^{k-1})^n$. For any
 $a\in {\bf Z}^n$ the bundle ${\cal O}(a)$ has exactly one $SL(k)$ -
linearization since the center of $SL(k)$ has dimension 0.
This linearization will be denoted by $\lambda$.
\vskip .3cm

\noindent (2.4.6) The bundle ${\cal O}(a)$ is ample if and only
 if all $a_i>0$. Assuming that this is the case,
let $B(k,n,a)= \bigoplus_d B(k,n,a)_d$ be the homogeneous coordinate
 ring of $(P^{k-1})^n$ in the projective embedding given by
${\cal O}(a)$. The degree $d$ homogeneous component $B(k,n,a)_d$
 of this ring consists of polynomials $F(w_1,...,w_n)$
in coordinates on $n$ vectors $w_i\in {\bf C}^k$ such that
$F(t_1w_1,...,t_nw_n) = t^{da} F(w_1,...,w_n)$ for any $t_i\in {\bf C}^*$.
Writing the vectors in coordinate form as columns $w_i =
(v_{1i},...,v_{ki})^t$, we realize elements of $B(k,n,a)_d$
 as polynomials $F(M)$ in entries of an indeterminate
$(k\times n)$ - matrix $M = ||v_{ij}||$ such that $F(M.t)
 = t^{da}F(M)$. The Mumford's quotient
 $((P^{k-1})^n/SL(k))_{{\cal O}(a),\lambda}$ is,
by definition, the projective spectrum $Proj (B(k,n,a)^{SL(k)})$.
\vskip .3cm

\proclaim (2.4.7) Theorem. Let $a=(a_1,...,a_n) \in {\bf Z}^n$
 be an integer vector. If at least one $a_i\leq 0$ then the
 Grassmannian $G(k,n)$ does not contain $a$ -stable orbits.
 If all $a_i$ are positive then we have an isomorphism
of Mumford quotients
$$(G(k,n)/H_1)_{{\cal O}(1),a)} \cong
((P^{k-1})^n/SL(k))_{{\cal O}(a),\lambda}.$$

\noindent {\sl Proof:} Both varieties are projective spectra
 of the same ring $R=\bigoplus R_d$ where $R_d$ consists of
 polynomials $\Phi(M), M\in Mat (k\times n)$ such that
$$\Phi(M.t) = t^{da} \Phi(M), t\in ({\bf C}^*)^n, \quad\quad
 \Phi (gM) = \Phi (M), g\in SL(k).$$

\noindent {\bf (2.4.9) Remark.} The algebra generated by
$k\times k$ -minors of an indeterminate $k\times n$ matrix
is known as the bracket algebra. It traditionally makes its
 appearance in two seemingly different contexts. The first
 appearance is as the coordinate ring of the Grassmannian
 $G(k,n)$ in the Pl\"ucker embedding. The other is in the
 study in the (semi-) invariants of system of vectors by
symbolic method (see [14]). However, the idea of serious
use of Grassmannians for the study of projective configurations
appeared only fairly recently in the papers of Gelfand and MacPherson.

\hfill\vfill\eject

\centerline{\bf  Chapter 3. VISIBLE CONTOURS OF (GENERALIZED) }
\centerline {\bf  LIE COMPLEXES AND VERONESE VARIETIES.}

\vskip 1cm

The Grassmannian point of view on projective configurations
 (i.e. the Gelfand- MacPherson isomorphism, see \S 3.2)
 simplifies considerably the study of the Chow quotient.
Indeed, instead of working with $(k^2-1)$ -dimensional
subvarieties on $(P^{k-1})^n$ which are closures of $PGL(k)$ -
orbits, we have to deal with
Lie complexes in $G(k,n)$ which are $(n-1)$ -dimensional toric
varieties.
\vskip 1cm

\centerline {\bf (3.1) Visible contours and the logarithmic Gauss map.}
\vskip .5cm

\noindent (3.1.1) There is a classical method (see e.g., [2,24,27])
 to analyze any complex of projective subspaces in $P^{n-1}$ i.e.
 an $(n-1)$ -dimensional subvariety $Z\i G(k,n)$. Namely, take any
 point $p\in P^{n-1}$ and consider the subvariety
$$Z_p = \{L\in Z: p\in L\}$$
of subspaces in $Z$ which contain $p$. This subvariety will be
called the {\it visible contour} of $Z$ at $p$.

Let $G(k-1,n-1)_p\i G(k,n)$ be the variety of all $(k-1)$ -
dimensional projective subspaces containing $p$. It is clear
that $Z_p = Z\cap G(k-1,n-1)_p$.
\vskip .3cm

\noindent (3.1.2) Still another step towards a visualization
 of the complex $Z$ at a point $p$ is done as follows. Let
$P^{n-1}_p$ be the space of lines in $P^{n-1}$ through $p$.
 Then $G(k-1,n-1)_p$ is identified with the variety of all
 $(k-2)$ -dimensional projective subspaces in $P^{n-2}_p$.
 We define the {\it visible sweep} of $Z$ at $p$ to be the
 subvariety $Sw_p(Z)\i P^{n-2}_p$ which is the union of all
 the projective subspaces corresponding to elements of $Z_p$.

\vskip .3cm

\noindent {\bf (3.1.3) Remarks.} a) If $k=2$ then $Z$ consists
 of lines in $P^{n-1}$. The lines belonging to the complex $Z$
 can be thought of as rays of light piercing the space, so $Z_p$
 is the contour which is seen by an observer at a point $p$.
In this case the visible contour is the same as the visible sweep.

b) Although the consideration of the locus $Z_p$ is classical,
 there seems to be no good name in the literature for it. The
term "complex cone" which is used sometimes [27] for the union
 of subspaces from $Z_p$ (i.e. the cone over the visible sweep,
 in our terminology) is obviously unsuitable for modern usage.

c) Dually, one can take any hyperplane $\Pi\i P^{n-1}$ and
consider the locus of subspaces from $Z$ which lie in $\Pi$.

\vskip .3cm

\noindent (3.1.4)  Since $codim G(k-1,n-1)_p = n-k$, we find
 that for any complex $Z\i G(k,n)$ and a generic $p$ the variety
$Z_p$ will have dimension $p-1$. Thus $Z_p$ will be a curve if $Z$
 consists of lines, a surface if $Z$ consists of planes etc.

\vskip .3cm

\noindent (3.1.5) We shall use the approach of visible contours
to study Lie complexes and, more generally, closures of arbitrary
 $(n-1)$ -dimensional torus orbits in $G(k,n)$ (such closures can
 be components of generalized Lie complexes).
Visible sweeps of Lie complexes will be studied in \S (3.6).
\vskip .3cm

\noindent (3.1.6) Let us realize our torus $H =
({\bf C}^*)^n/{\bf C}^*$ as an open subset in $P^{n-1}$
 consisting of points with all homogeneous coordinates non-zero.
The point $e=(1:...:1)\in P^{n-1}$ becomes the unit in $H$.
 Denote by {\bf h} the Lie algebra of $H$. It is identified
with the tangent space to $P^{n-1}$ at $e$. Explicitly, ${\bf h}
 = {\bf C}^n/\{(a,...,a)\}$.
For any $x\in H$ let $\mu_x:H\rightarrow H$ be the operator of
 multiplication by $X$.

Any subvariety $L\i P^{n-1}$ not lying inside a coordinate
hyperplane gives a subvariety $L\cap H$ in the algebraic  group $H$.

\proclaim (3.1.7) Definition. Let $X\i H$ be a $p$ -dimensional
 algebraic subvariety. The logarithmic Gauss map of $X$ is the
(rational) map $\gamma_X:X\rightarrow G(p,{\bf h})$ which takes
 a smooth point $x\in X$ to the $p$ -dimensional subspace
$d(\mu_x^{-1})(T_xX)\i T_eH = {\bf h}$  --- the translation
 to the unity of the tangent space $T_xX\i T_xH$.

The name "logarithmic" comes from the fact that explicit formula
 for $\gamma_X$ involves logarithmic derivatives (see below).
\vskip .3cm

\noindent (3.1.8) Let  $L\i P^{n-1}$ be a $(k-1)$ -dimensional
 projective subspace not lying in a coordinate hyperplane.
 The orbit closure $\overline{H.L}$ has dimension $n-1$ i.e.
 this is a complex. Since this complex is $H$ - invariant,
its visible contour $(\overline{H.L})_p$ at any point
$p\in P^{n-1}$ with all coordinates non-zero, will be isomorphic
 to the visible contour at the point $e=(1:...:1)$.

Before stating the next proposition let us note that the
 Grassmannian $G(k-1,n-1)_e$, where visible contours lie,
 is canonically identified with $G(k-1,{\bf h})$. (Correspondingly,
 the space $P^{n-2}_e$ where visible sweeps lie, is $P({\bf h})$.)
\vskip .3cm

\proclaim (3.1.9) Proposition. a) If the subspace $L$ does not lie
 in a coordinate hyperplane then the visible contour
$(\overline{H.L})_e$ coincides with the closure of the
 image of $L\cap H$ under the logarithmic Gauss map.
In particular, this visible contour is a rational variety. \hfill\break
b) The intersection of $\overline{H.L}$ with the
 sub-Grassmannian $G(k-1,n-1)_e$ is proper and transversal at its
 generic point.

\noindent {\sl Proof:} a) Neither the complex $\overline{H.L}$
 nor the image of $L\cap H$ under the logarithmic Gauss map will
 change if we translate $L$ by the $H$ -action. So we can (and will)
 assume that $L$ contains the point $e=(1:...:1)$. For $h\in H$ the
 translated subspace $h^{-1}L$ contains $e$ if and only if $h\in L$.
 Thus the variety $(H.L)_e = (H.L)\cap G(k-1,n-1)_e$ consists of
subspaces $h^{-1}.L, h\in L$. In other words, $(H.L_e$ is the image
 of the map $L\cap H \rightarrow G(k-1,n-1)_e$ taking $h\in L\cap H$
 to the subspace $h^{-1}.L$. This map clearly coincides with the
logarithmic Gauss map.

b) If $Z$ is any complex in $G(k,n)$ then the assertion will be
 true for the
intersection of $Z$ with $G(k-1,n-1)_p$, where $p\in P^{n-1}$ is
a generic point. In our case, due to the invariance under the torus
 action, the situation at any $(p=(p_1:...:p_n)$ with all $p_i\neq 0$
is the same as at $e$.
\vskip .3cm

\proclaim (3.1.10) Theorem. Suppose that $L\i P^{n-1}$ is a $(k-1)$ -
dimensional subvariety belonging to the generic stratum $G^0(k,n)$.
 Then the logarithmic Gauss map $\gamma_{L\cap H}$ extends to a regular
 embedding $L\hookrightarrow G(k-1, {\bf h}) = G(k-1,n-1)_e$.
In other words, the visible contour $(\overline{H.L})_e$ is
identified with $L$ itself.

\noindent {\sl Proof:} First let us show that the logarithmic
Gauss map $\gamma$
extends to all of $L$ as a regular map. This will be done by
 calculation in coordinates which we shall also use on other occasions.
\vskip .3cm

 \noindent (3.1.11) Let  $x_1,...,x_n$ be homogeneous coordinates
in $P^{n-1}$. Let $p= (y_1:...:y_n)\in L$ be any point. Since $L$
 lies in the generic stratum, there are $1\leq i_1<...<i_{k-1}\leq n$
 such that $y_j\neq 0$ for $j\notin \{i_1,...,i_{k-1}\}$. After
 renumbering variables we can (and will) assume that $y_k,...,y_n$
 are non-zero. Consider the affine space $L-\{x_n=0\}$ which
 contains our point $p$. Introduce in this space affine
coordinates $z_1,...,z_{k-1}$ where $z_i = x_i/x_n$.
We can set $x_n$ to be 1 on $L-\{x_n=0\}$ and express all
the other coordinates as affine-linear functions in $z_i$ i.e.
$$x_i = z_i, i=1,...,k-1,\quad\quad x_i = f_i(z) =
 \sum_{\nu=1}^{k-1} a_{i\nu}z_\nu + a_k, i=k,...,n-1.$$
We also set $f_i(z) = z_i$ for $i=1,...,k-1$.
\vskip .3cm

\noindent (3.1.12) We identify the torus
 $H=({\bf C}^*)^n/{\bf C}^*\i P^{n-1}$ with the set
  $\{(t_1,...,t_{n-1},1)\in ({\bf C}^*)^n\}$ i.e. with
 $({\bf C}^*)^{n-1}$. Its Lie algebra is therefore identified
 with ${\bf C}^{n-1}$. In this notation the map $\gamma$ takes
a point $z=(z_1,...,z_{k-1})$ to the $(k-1)$ -dimensional
 subspace in ${\bf C}^{n-1}$ spanned by the rows of $(k-1)$
 by $(n-1)$ -matrix
$||\partial {\rm log} f_i /\partial z_j||, i=1,...,n-1, \,
 j=1,...,k-1$. We can multiply the $j$ -th row by $z_j$
without changing this subspace. After this the matrix takes the form
$$\pmatrix{ 1&0&...&0&{a_{k1}z_1\over f_k(z)}&...&...&
{a_{n-1,1}z_1\over f_{n-1}(z)}\cr
0&1&...&0&{a_{k2}z_2\over f_k(z)}&...&...&{a_{n-1,2}z_2\over f_{n-1}(z)}\cr
...&...&...&...&...&...&...&...\cr
...&...&...&...&...&...&...&...\cr
0&0&...&1&{a_{k,k-1}z_{k-1}\over f_k(z)}&...&...&
{a_{n-1,k-1}z_{k-1}\over f_{n-1}(z)}}\leqno (3.1.13)$$
This matrix is clearly regular near our point $p$ since
 $f_k(p),...,f_{n-1}(p)$ are non-zero. The rank of the
 matrix (3.1.13) being equal $k$, we deduce that $\gamma$
is regular at $p$. We have proven that $\gamma$ extends to
 a regular morphism $L\rightarrow G(k-1,n-1)$.
\vskip .3cm

\noindent (3.1.14) Let us finish the proof of Theorem 3.1.10
 by showing that the logarithmic Gauss map $\gamma$ is an
 embedding. Consider the set of all $(k-1)$ by $(n-1)$ matrices
of which the first $(k-1)$ columns form the unit $(k-1)$ by
 $(k-1)$ matrix. The entries of the remaining $(n-k)$ columns
 are independent affine coordinates in the open Schubert cell
 ${\bf C}^{(k-1)(n-k)}\i G(k-1,n-1)$. Let us show that entries
 of any given column of (3.1.13) whose number is greater that $k$,
 alone suffice to separate all points of $L$. Indeed, consider, say,
 the $p$ -th column, $p>k$ and regard its entries as defining a
transformation
$(z_1,...,z_{k-1})\mapsto (s_1,...,s_{k-1})$ where $s_i =
 {a_{pi}z_i\over \sum_{\nu=1}^{k-1}a_{p\nu}z_\nu}$. This
 is a projective transformation corresponding to the $k$ by $k$ matrix
$$ T_p = \pmatrix{a_{p1}&0&0&...&0\cr
0&a_{p2}&0&...&0\cr
...&...&...&...\cr
a_{p1}&a_{p2}&a_{p3}&...&a_{pk}} \leqno (3.1.15)$$
Since our subspace $L\i P^{n-1}$ belongs to the generic stratum,
 every entry of the matrix $a_{ij}, i=1,...,k-1, j=1,...,n-k$,
 is non-zero. Hence the matrix (3.1.15) defines a non-degenerate
 projective transformation and separates the points (as well as
 the tangent vectors). Theorem 3.1.10 is proven.

\vskip 1cm

\centerline {\bf (3.2) Bundles of logarithmic forms on $P^{k-1}$
 and visible contours.}
\vskip .5cm

\noindent (3.2.1) It is well known [24] that maps from any
 projective variety $X$ to Grassmannians of the form $G(r,V),
 dim (V)>r$, are in
  correspondence with rank $r$ vector bundles on $X$. More precisely,
 given such a bundle $E$ we consider the vector space $V= H^0(X,E)^*$.
Suppose that $E$ is generated by global sections and let $N$ be the
 dimension of $V$. Define a map $\phi_E:X\rightarrow G(r,V) =
G(N-r,H^0(X,E))$  as follows. For a point $x\in X$ the value
 $\phi_E(x)$ is the codimension $r$ subspace in $H^0(X,E)$
consisting of all the section vanishing at $x$.  Conversely,
 suppose given a map $\phi: X\rightarrow G(r,V)$. Let $S$ be
the tautological rank $r$ bundle on $G(r,V)$ (whose fiber at
a subspace $L$ is $L$ itself). Associate to $\phi$ the bundle
 $\phi^*S^*$ on $X$.

The bundles on $P^{k-1}$ corresponding to visible contours of
 Lie complexes have the following description.
\vskip .3cm

\noindent (3.2.2) Let $M=(M_1,...,M_n)$ be a configuration of $n$
 hyperplanes in $P^{k-1}$ which are in general position. Then $M$
 is a divisor with normal crossings and we can define the sheaf
 $\Omega^1_{P^{k-1}}({\rm log} M)$ of differential 1-forms on
 $P^{k-1}$ with logarithmic poles along $M$, see [13].
 By definition, the space of sections of this sheaf near a
 point $x\in P^{k-1}$ is generated (over ${\cal O}_{P^{k-1},x}$)
by 1-forms regular at $x$ and also by forms ${\rm dlog} f_i$ where
 $f_i$ are local equations of hyperplanes from $M$ containing $x$.

An important property of the sheaf $\Omega^1_{P^{k-1}}({\rm log} M)$
is that it is locally free i.e. can be seen as a rank $(k-1)$ vector
 bundle over $P^{k-1}$.

\proclaim (3.2.3)
Proposition. Let $M=(M_1,...,M_n)$ be a configuration of $n$
hyperplanes in the projective space $L = P^{k-1}$ which are
 in general position.
 and let $f_i$ be a linear form defining $M_i$. Then:\hfill\break
a) The space $W = H^0(L, \Omega^1_L({\rm log} M))$ has dimension
 $n-1$ and consists of forms
$$\sum_i \alpha_i {\rm dlog} f_i = {\rm dlog} \prod_i f_i^{\alpha_i},
 \alpha_i\in {\bf C}, \sum \alpha_i =0.$$
Higher cohomology groups of $\Omega^1_L({\rm log} M)$ vanish.
 \hfill\break
b) The vector bundle $E= \Omega^1_L({\rm log} M)$ defines a
 regular embedding $\phi_E:L\hookrightarrow G(k-1,W^*)$. \hfill\break
c) Suppose that  $L$ is realized as a subspace in the coordinate
 $P^{n-1}$ so that $M_i$ is given by the vanishing of the $i$ -th
coordinate. Then $\phi_E$ coincides with the (extension of) the
logarithmic Gauss map $\gamma_{L\cap H}$, and the image $\phi_E(L)$
 coincides with the visible contour of the Lie complex
$\overline{H.L}$.
\vskip .3cm

\noindent (3.2.4) \underbar {\sl Proof of (3.2.3),a):} The sheaf
$\Omega^1_L$
 of regular 1-forms is obviously a subsheaf of $\Omega^1_L({\rm log} M)$.
 To describe the quotient, we shall, following P.Deligne [13],
 denote by $\tilde M$ the disjoint union of hyperplanes in $M$
and let $\epsilon:\tilde M\rightarrow L$ be the natural map.
Then we have the exact sequence
$$0\rightarrow \Omega^1_L \rightarrow \Omega^1_L({\rm log} M)
\buildrel Res\over\longrightarrow \epsilon_*{\cal O}_{\tilde M}
\rightarrow 0 \leqno (3.2.5)$$
where $Res$ is the Poincar\'e residue morphism, see [13].
Consider the corresponding long exact sequence of cohomology.
 The equality $H^0(L,\epsilon_*{\cal O}_{\tilde M}) = {\bf C}^n$
 means that the residue of a global logarithmic form along each
$M_i$ is constant. The sum of the residues is given by the boundary map
$H^0(L,\epsilon_*{\cal O}_{\tilde M}) \rightarrow H^1(L,\Omega^1)
 = {\bf C}$ should be zero. Since the forms exhibited in the
formulation are indeed global sections of our sheaf, we obtain
 the statement about $H^0$. The vanishing of higher
$H^i$ follows from known information about the cohomology of
the sheaf ${\cal O}$ on $P^{k-2}$ and $\Omega^1$ on $P^{k-1}$.
\vskip .3cm

\noindent (3.2.6) \underbar {\sl Proof of (3.2.3) b) and c):}
 We can assume that $M$ is given by the intersection of an
embedded $L=P^{k-1}\i P^{N-1}$ with coordinate hyperplanes
$\{x_i=0\}$. Then, by n.a), the basis of
$H^0(L,\Omega^1_L({\rm log} M))$ is given by 1-forms
 ${\rm dlog} (x_1/x_n), i=1,..., n-1$. We identify the space
 of section with ${\bf C}^{n-1}$ by using this basis. Now
looking at explicit formula (3.1.12), we find that the map
 $\phi_E:L\rightarrow G(k-1,n-1)$ is defiend by the formula
identical to that of the logarithmic Gauss map.
\vskip .3cm

\proclaim (3.2.7) Proposition. The Chern classes of
$E= \Omega^1_{P^{k-1}}({\rm log} (M_1+...+m_n))$  have the form
$$c_i(E) = {n-k+i\choose i} \in H^{2i} (P^{k-1}, {\bf Z}) = {\bf Z}.$$
In particular, the determinant (= top exterior power) of $E$ is
 isomorphic to ${\cal O}_{P^{k-1}}(n-k)$.

\noindent {\sl Proof:} This follows at once from the exact sequence (3.2.5).
\vskip .3cm

\noindent {\bf (3.2.8) Example.} Consider a Lie complex in $G(2,n)$,
 the Grassmannian of lines in $P^{n-1}$. Let this complex have the
 form $Z=\overline{H.l}$, where $l$ is a line belonging to the
generic stratum. The visible contour $Z_e$ lies
in the projective space $P^{n-2}_e$ of all lines in $P^{n-1}$
 through the point $e$.
Proposition 3.2.3 means that $Z_e$ is a rational normal curve
(Veronese curve, for short) in $P^{n-2}_e$. More precisely, it
 is the embedding of $l=P^1$ defined by the invertible sheaf
 $\Omega^1_l({\rm log} (m_1+...+m_n))\cong {\cal O}_l(n-2)$.
 Here $m_i\in l$ is the point of intersection of $l$ with the coordinate
 hyperplane $\{x_i=0\}$.

\vskip 1cm

\centerline {\bf (3.3)  Visible contour as a Veronese variety
 in the Grassmannian.}
\vskip .5cm

\noindent (3.3.1) Recall [24,25] the $d$ -fold {\it Veronese embedding}
$$P^{k-1} = P({\bf C}^k) \hookrightarrow P(S^d{\bf C}^k),
\quad\quad x\mapsto x^d\leqno (3.3.2)$$
of $P^{k-1}$ into the projectivization of $(S^d{\bf C}^k)$,
 the space of homogeneous degree $d$ polynomials in $k$ variables.
 This is the embedding corresponding to the line bundle ${\cal O}(d)$.
We shall say that a $(k-1)$ -dimensional subvariety $X\i P^N$ is a
 $d$ -fold {\it Veronese variety} if there is a projective equivalence
 $P^N\cong P(S^d{\bf C}^k)$ taking
$X$ into the image of (3.3.2). A Veronese curve in $P^N$ is the same
as a rational normal curve of degree $N$.
\vskip .3cm

\noindent (3.3.3) The dimension of of the projective space $P(S^d{\bf C}^k)$
of the $d$ -fold Veronese embedding (3.3.2) equals $N= {d+k-1\choose d-1}-1$.
 Note that the same dimension is attained by the projective space of the
 Pl\"ucker embedding of the Grassmannian $G(k-1,d+k-1)$. Therefore it
makes sense to look for Veronese subvarieties in
$P^{{d+k-1\choose d-1}-1}= P(\bigwedge ^{k-1}{\bf C}^{d+k-1})$
 which lie on the Grassmannian.
\vskip .3cm

\noindent (3.3.4) We shall say that a $(k-1)$ -dimensional subvariety
$X\i G(k-1,d+k-1)$ is a $d$ -fold Veronese variety if it becomes such
after the Pl\"ucker embedding of
$G(k-1,d+k-1)$.
\vskip .3cm

\proclaim (3.3.5) Proposition. Let $M=(M_1,...,M_n)$ be a
configuration of hyperplanes in $P^{k-1}$ in general position,
 $E= \Omega^1_{P^{k-1}}({\rm log} M)$ -the corresponding
logarithmic bundle and $\phi_E:P^{k-1} \rightarrow G(k-1,n-1)$ -
 the embedding corresponding to $E$. Then $\phi_E(P^{k-1})$ is
an $(n-k)$ -fold Veronese variety in $G(k-1,n-1)$.

\noindent {\sl Proof:} Let us construct an isomorphism of linear spaces
$\bigwedge^{k-1}{\bf C}^{n-1}\rightarrow S^{n-k}{\bf C}^k$
taking $\phi_E(P^{k-1})$ into the standard Veronese variety.
 Let $S$ be the tautological rank $(k-1)$ bundle on $G(k-1,n-1)$.
 Then $E=\phi_E^*(S^*)$ and hence
$$\phi_E^*(\bigwedge^{k-1}S^*) = \bigwedge^{k-1}E =
 {\cal O}_{P^{k-1}}(n-k)$$
 by Proposition 3.2.7. Thus we obtain a linear map of restriction
$$\bigwedge^{k-1}{\bf C}^{n-1}= H^0(G(k-1,n-1),
 \bigwedge^{k-1}S^*)\buildrel r\over\longrightarrow H^0(P^{k-1},
\phi_E^*(\bigwedge^{k-1}S^*)) \cong \leqno (3.3.6)$$
$$ \cong H^0(P^{k-1}, {\cal O}(n-k)) = S^{n-k}{\bf C}^k.$$
\vskip .3cm

\noindent (3.3.7)
Let us show that the restriction  map $r$ in (3.3.6) is an isomorphism.

 Since spaces in both sides have the same dimension, it suffices to show
that $r$ is injective i.e. that the variety $X = \phi_E(P^{k-1})$
does not lie in any hyperplane in $P(\bigwedge^{k-1}{\bf C}^{n-1})$.
Take an affine chart in $P^{k-1}$ in which the last hyperplane $M_n$
is the infinite one. All the other hyperplanes $M_i$ are then defined
 by vanishing of affine-linear functions
$f_i, i=1,...,n-1$ on ${\bf C}^{k-1}= P^{k-1}-M_n$.
 The fact that the variety $X$ lies in a hyperplane means that
 there is a collection of numbers $a_{i_1,...,i_{k-1}}$, not all of them zero,
 such that the meromorphic $(k-1)$ -form
$$\Omega = \sum_{1\leq i_1<...<i_{k-1}\leq n-1}
a_{i_1,...,i_{k-1}}{\rm dlog}f_{i_1} \wedge ...
\wedge {\rm dlog} f_{i_{k-1}}$$
on ${\bf C}^{k-1}$ vanishes identically. However, the coefficient
 $a_{i_1,...,i_{k-1}}$ can be read off $\Omega$ as the residue at
the intersection point $M_{i_1}\cap ...\cap M_{i_{k-1}}$ so all
the coefficients should be zero. Proposition 3.3.5. is proven.

\vskip .3cm

\noindent (3.3.8) Let ${\bf h} $ be the quotient of ${\bf C}^n$
by the subspace of $(a,...,a), a\in {\bf C}$. Note that this subspace
is canonically identified with the Lie algebra of the torus $H$ and with
the tangent space to $P^{n-1}$ at the point $e=(1,...,1)$. We shall denote,
therefore, by $G(k-1,n-1)_e$ the Grassmannian
of $(k-1)$ -dimensional subspaces in ${\bf h}$.
\vskip .3cm

\proclaim (3.3.9) Definition. By a special Veronese subvariety in
 $G(k-1,n-1)_e$ we shall mean a subvariety of the form $\phi_E(P^{k-1})$,
 where:\hfill\break
a) $E=\Omega^1_{P^{k-1}}({\rm log}M)$ is the logarithmic bundle
corresponding to some configuration $M=(M_1,...,M_n)$ of hyperplanes
in general position;\hfill\break
b) The space $H^0(E)$ is identified with $\{(a_1,...,a_n)\in
 {\bf C}^n: \sum a_i=0\}$ as in Proposition 3.2.3, and its dual -
 with ${\bf h}$.

Thus the notion of a special Veronese variety makes an
 explicit appeal to a choice of coordinate system.
\vskip .3cm

\noindent (3.3.10)  Note that by Proposition 3.2.3 special
 Veronese varieties are precisely the visible contours of Lie
 complexes in $G(k,n)$. In particular, these variety define,
around a generic point of $G(k-1,n-1)_e$, a foliation with $k-1$
- dimensional fibers which is just the intersection of
$G(k-1,n-1)_e$ with the
 foliation given by the orbits of $H$. Let us note also the
 following corollary.
\vskip .3cm

\proclaim (3.3.11) Corollary. The set $G^0(k,n)/H =
(P^{k-1})^n_{gen}/GL(k)$ of projective equivalence \hfill\break
classes of configuration of $n$ hyperplanes in $P^{k-1}$
in general position is in one-to-one correspondence with
the set of special Veronese varieties in $G(k-1,n-1)_e$.
This correspondence takes a configuration $M=(M_1,...,M_n)$
into the subvariety $\phi_E(P^{k-1})$, where $E= \Omega^1({\rm log}M)$.

\vskip .3cm

\noindent (3.3.12) Clearly all special Veronese varieties in $G(k-1,n-1)_e$
 represent the same homology class $\Delta \in H_{2k-2}(G(k-1,n-1),{\bf Z})$.
 A precise determination of $\Delta$ will be given in \S
3.9 below. By Corollary 3.3.11 we obtain an embedding of $G^0(k,n)/H$
into the Chow variety ${\cal C}_{k-1}(G(k-1,n-1),\Delta)$.
Denote by $V$ the closure of $G^0(k,n)/H$ in this Chow variety.
So it is the variety of cycles in ${\cal C}_{k-1}(G(k-1,n-1),\Delta)$
 which are limit positions of special Veronese varieties.
\vskip .3cm

\noindent (3.3.13) Similarly, all special Veronese varieties
 in $G(k-1,n-1)_e$ form a flat family. Let ${\cal H}$ be the
Hilbert scheme parametrizing all subschemes
in $G(k-1,n-1)$, cf. (0.5.1). Define the variety $W$ to be
 the closure of $G^0(k,n)/H$
 in the Hilbert scheme ${\cal H}$. So it is the variety of
subschemes in ${\cal C}_{k-1}(G(k-1,n-1),\Delta)$ which are
 limit positions of special Veronese varieties.

Our next result shows that all the information about the Chow
quotient $G(k,n)//H$ is contained in visible contours.
\vskip .3cm

\proclaim (3.3.14) Theorem. The correspondence $Z\mapsto Z_e$
extends to an isomorphism of the variety $G(k,n)//H$ with $V$ and $W$.

\noindent {\sl Proof:} Since, by Proposition 3.1.9 b), every orbit
 closure $\overline{H.L}$ intersects the variety $G(k-1,n-1)_e$
 properly, we can conclude, by using the result of Barlet [4]
 that the map $Z\mapsto Z_e = Z\cap G(k-1,n-1)_e$ defines a
regular morphism $\psi:G(k,n)//H\rightarrow V$. Proposition
3.1.9 implies that $\psi$ is set-theoretically a bijection.
 To show that this is an isomorphism of algebraic varieties
 it suffices to apply once again the reasoning with normal
 vector fields on a generalized Lie complex $Z$ used in the
 proof of Theorem 1.5.2. Similarly for the Hilbert scheme
compactification.
\vskip .3cm

\noindent {\bf (3.3.15) Remark.} A natural problem would be
to study all Veronese subvarieties in Grassmannians. In general,
 not every such variety is projectively equivalent to a special one.
 This is because there are rank $(k-1)$ vector bundles on $P^{k-1}$
 which have the same Chern classes as $\Omega^1({\rm log}M)$ but
not having this form. A study of bundles $\Omega^1({\rm log}M)$ from
 the point of view of stable vector bundles on projective spaces
 will be undertaken in a subsequent paper of I.Dolgachev and the author.

\vskip 1cm

\centerline {\bf (3.4) Properties of special Veronese varieties.}
\vskip .5cm

As has been recalled in \S 2, Lie complexes in $G(k,n)$ correspond
to projective equivalence classes of configurations of $n$
 hyperplanes in $P^{k-1}$ in general position. We have seen
 in the previous subsection that the space $P^{k-1}$ can be
recovered from the corresponding Lie complex $Z$ as its
visible contour $Z_e$. Let us recover the configuration too.
\vskip .3cm

\noindent (3.4.1) For any points $x_1,...,x_m\in
 P^{n-1}$ let $<x_1,...,x_m>$ denote their projective span.
 We shall define also by $G_{<x_1,...,x_m>}$ the subvariety
in $G(k,n)$ formed by $P^{k-1}$'s containing $<x_1,...,x_m>$.
 As an abstract variety, it is isomorphic to $G(k-p, n-p)$,
where $p= dim<x_1,...,x_m> + 1$. Let $e_i\in P^{n-1}$ be the
 images of the standard basis vectors of ${\bf C}^n$.

\proclaim (3.4.2) Proposition. Let $M=(M_1,...,M_n)$ be a
 configuration of $n$ hyperplanes in $P^{k-1}$ in general
position, $E=\Omega^1_{P^{k-1}}({\rm log} M)$, and
$X = \phi_E(P^{k-1})\i G(k-1,n-1)_e$ -the corresponding
special Veronese variety (i.e. the visible contour of the Lie
 complex corresponding to $M$). Then $\phi_E(M_j) = X \cap G_{<e,e_j>}$.

\noindent {\sl Proof:} Let us give a coordinate description of
$\phi_E$ which, unlike the description given in (3.1.11), is
symmetric with respect to permutation of hyperplanes.
\vskip .3cm

\noindent (3.4.3) Let $z_1,...,z_k$ be homogeneous coordinates
in $P^{k-1}$ and $g_j(z) = \sum_i a_{ij}z_j$ -the linear equations
 of $M_j, j=1,...,n$. Let ${\bf h}$ denote, as before, the quotient
${\bf C}^n/\{(a,...,a)\}$. The map $\phi_E: P^{k-1}\rightarrow
 G(k-1, {\bf h}) = G(k-1,n-1)_e$ is defined as follows.
\vskip .3cm

\noindent (3.4.5) Let $z=(z_1:...:z_k)\in P^{k-1}$ be generic.
 Consider the Jacobian $(k\times n)$ -matrix
$N(z) = ||\partial {\rm log} f_i/\partial z_j||$. Due to the
identity $\sum z_j {\partial {\rm log} f_i\over\partial z_j} =
1, \,\, \forall i$, the $k$ -dimensional subspace spanned by
rows of this matrix contains the vector $(1,...,1)$ and hence
 defines a $(k-1)$ -dimensional subspace in ${\bf h}$ which is
 precisely $\phi_E(z)$. For any subset $I\i \{1,...,n\}, |I|=k$,
denote by $p_I(N(z))$ the $k\times k$ -minor of $N(z)$ on columns
 from $I$.
\vskip .3cm

\noindent (3.4.6) We can (and will) assume, by renumbering the
coordinates, that the number $j$ in the formulation of Proposition
 3.4.2  equals 1. The subspace generated by rows of $N(z)$ lies in
 the sub-Grassmannian $G_{<e,e_1>}$ if and only if all minors
 $p_I(N(z)), 1\notin I$, vanish. To eveluate the limit of this
subspace for $x\rightarrow M_i$, multiply $N(z)$ by $diag
(f_1(z),...,f_k(z))$ from the left. Then any minor $p_I$ of
$diag (f_1(z),...,f_k(z)).N(z)$ with $1\notin I$, will contain
a row vanishing on $M_i$ and hence will vanish on $M_i$ itself.
 On the other hand, the minor $p_{1,2,...,k}$ of
 $diag (f_1(z),...,f_k(z)).N(z)$ is constant.
Since we can write instead of $(f_1,...,f_k)$ any
 $(f_1,f_{i_2},...,f_{i_k}), 1<i_2<...<i_k\leq n$,
this proves that $\phi_E(M_1) = \phi_E(P^{k-1})\cap G_{<e,e_1>}$.
 Proposition 3.4.2 is proven.
\vskip .3cm

\proclaim (3.4.7) Corollary. Any special Veronese variety
contains the $n\choose k-1$ points
$$<e,e_{i_1},...,e_{i_{k-1}}> \in G(k-1,n-1)_e,
 1\leq i_1<...<i_{k-1}\leq n.$$
\vskip .3cm

\proclaim (3.4.8) Proposition. The intersection
$\phi_E(P^{k-1})\cap G_{<e,e_j>}$ is itself a special
 Veronese variety corresponding to the projective space
$M_j$ and the configuration of hyperplanes $(M_i \cap M_j), i\neq j$.

\noindent {\sl Proof:} Straightforward. Left to the reader.
\vskip .3cm

\noindent {\bf (3.4.9) Example.} Consider the case $k=2$ when
 $G(k,n)$ consists of lines in $P^{n-1}$. The variety of Lie
complexes in $G(2,n)$ is, by \S 2, the same as the quotient
 $((P^1)^n - \bigcup \{x_i=x_j\})/GL(2)$ i.e. the set of
 projective equivalence classes of $n$ -tuples of distinct
 points on $P^1$. As we have seen in Example 3.2.8 that the
 visible contour of any  Lie complex in $G(2,n)$ is a
Veronese curve in $P^{n-2}_e$, the variety of lines in
$P^{n-1}$ throught $e$. Corollary 3.4.7 means that every
special Veronese curve in $P^{n-2}_e$ contains $n$ points
 $<e,e_i>$ which are in general position.

 It is a classical fact that for any points $p_1,...,p_n\in P^{n-2}$
in general position the set $V_0(p_1,...,p_n)$ of all Veronese curves
 through $p_i$ is in bijection with $(P^1)^n_{gen}/GL(2)$ (see [14],
 Ch.III, \S 2, Proposition 3).
\vskip .3cm

\noindent {\bf (3.4.10) Example (continued).} As a transparent
 particular case, consider the case of 4 points $p_1,...,p_4$ in
 $P^2$. Veronese curves in $P^2$ are just smooth conics. Conics
 through $p_1,...,p_4$ form a 1-dimensional pencil ${\cal L}=P^1$.
 There are exactly three degenerate conics in this pencil namely
 unions of lines
$$<p_1,p_2>\cup <p_3,p_4>,\,\, <p_1,p_3>\cup <p_2,p_4>,\,\,
 <p_1,p_4>\cup <p_2,p_3>.$$
The set of cross-ratios of $p_1,...,p_4$ regarded on conics from
 ${\cal L}$ is in bijection with the set of non-degenerate conics
 from ${\cal L}$ i.e. with $P^1$ minus 3 points.

\vskip 1cm

\centerline {\bf (3.5) Steiner constructions of Veronese varieties
in Grassmannians.}
\vskip .5cm

Veronese curves in projective spaces posess a lot of remarkable properties,
see [2,14,46]. Most of these properties do not generalize to higher-dimensional
 Veronese varieties. It is our opinion that the "right" class of ambient
 spaces for $p$ -dimensional Veronese varieties is formed not by projective
spaces but by Grassmannians of the form $G(p,V)$. In this section we shall
show that (special) Veronese varieties admit a "synthetic" construction  in
the spirit of Steiner.
\vskip .3cm

\noindent (3.5.1) Consider some projective space $P^{m-1}$.
  Let $L\i P^{m-1}$ be a projective subspace of codimension $d$.
 Denote by $]L[$ the space of all hyperplanes in $P^{m-1}$ containing $L$.
We shall call it the {\it star} of $L$. This is a projective space of dimension
$d - 1$.
\vskip .3cm

\noindent (3.5.2) We may like to have a {\it parametrization} of the star $]L[$
i.e. an identification
 $f: P^{d-1}\rightarrow ]L[$
of $]L[$ with the standard $P^{d-1}$. Such an identification is the same as
 a linear operator $E:{\bf C}^{m}\rightarrow {\bf C}^d$ whose kernel
is  ${\bf L}$, the linear subspace corresponding to $L$.
 Indeed, given such an $E$, we obtain a bijection $\Pi\mapsto
 E^{-1}(\Pi)$ between hyperplanes in ${\bf C}^d$ and
 hyperplanes in ${\bf C}^m$
containing ${\bf L}$ i.e. hyperplanes from $]L[$.

In coordinate notation,  we write $E$ as a row of linear functions
$g_i:{\bf C}^m\rightarrow {\bf C}$.
Then
to any $(\lambda_1:...:\lambda_d)\in P^{d-1}$ we associate the hypleplane
$Ker (\sum \lambda_i g_i) \in ]L[$.
\vskip .3cm

\noindent (3.5.3)
It will be convenient for us to view a parametrization
$f$ above  as a linear form
$\sum \lambda_i g_i$
on ${\bf C}^m$ whose entries are linear forms in $\lambda_1,...,\lambda_k$.
This is tantamount to viewing a linear operator
 $E:{\bf C}^{m}\rightarrow {\bf C}^d$ as
 an element of
${\bf C}^d \otimes ({\bf C}^m)^*$.

\vskip .3cm

\noindent (3.5.4) Recall  Steiner's construction of Veronese curves
in  $P^m$ [24]. Take $m$ projective subspaces of codimension 2,
 $L_1,...,L_m\i P^m$. The star $]L_i[$ of each $L_i$ is just a
pencil of hyperplanes i.e. it is isomorphic to the projective
line $P^1$. Let us identify these pencils with each other, e.g.
 by choosing projective equivalence $f_i:P^1\rightarrow ]L_i[$.
Consider the curve in $P^m$ which is the image of $P^1$ under the map
$$t \mapsto (f_1(t)\cap ...\cap f_m(t)).$$
This is a Veronese curve. It depends on the choice of subspaces $L_i$ and
of identifications $f_i$.

In classical terminology, one would say that a Veronese curve can
be obtained as the locus of intersections of corresponding
hyperplanes from $m$ pencils in correspondence.

\vskip .3cm

\proclaim (3.5.5) Construction. ({\rm The Grassmannian Steiner
 construction.})
 Take $n-k$ projective subspaces in $P^{n-2}$, say,
 $L_1,...,L_{n-k}$, of codimension $k$. Put the stars
 $]L_i[$ into 1-1 correspondence with each other, e.g.,
 by choosing projective isomorphisms $f_i:P^{k-1}\rightarrow ]L_i[$.
Then consider the subvariety in $G(k-1,n-1)$ given by the parametrization
$$t \mapsto (f_1(t)\cap ...\cap f_{n-k}(t)),\,\, t\in P^{k-1}.$$

This is a direct generalization of the construction in (3.5.4).
 Using the fact (3.5.3)
that parametrized stars are the same as linear forms with
coefficentls linearly depending on parameters, we get the
following reformulation of the construction.
\vskip .3cm

\proclaim (3.5.6) Reformulation. Take a linear operator
 $A:{\bf C}^k\rightarrow Hom ({\bf C}^{n-1}, {\bf C}^{n-k})$
 such that for any non-zero $z\in {\bf C}^k$
the operator $A(z):{\bf C}^{n-1}\rightarrow {\bf C}^{n-k}$ is surjective.
The Grassmannian Steiner construction is the subvariety in $G(k-1, n-1)$
consisting of points $Ker A(z), z\in {\bf C}^k-\{0\}$.

 \proclaim (3.5.7) Theorem. Any special $(k-1)$ -dimensional
Veronese variety in
$G(k-1,n-1)$ can be obtained by  the Grassmannian Steiner  construction.

\noindent {\sl Proof:} Let $X$ be a special Veronese variety
 coming from a configuration $(M_1,...,M_n)$ of hyperplanes
in $P^{k-1}$ in general position.
As in (3.1.11) we can assume that $M_n$ is the infinite hyperplane
and choose affine coordinates $z_1,...,z_{k-1}$ in ${\bf C}^{k-1}
= P^{k-1} -M_n$ such that
$M_i$ is given by the equation $z_i=0$ for $i=1,...,k-1$.
\vskip .3cm

\noindent (3.5.8) Consider the  coordinate space ${\bf C}^{n-1}$
 with coordinates $y_1,...,y_n$.
and basis vectors $e_1,...,e_{n-1}$
Decompose it
into the direct sum ${\bf C}^{k-1} \oplus {\bf C}^{n-k}$
where ${\bf C}^{k-1}$ is spanned by $e_1,...,e_{k-1}$ and ${\bf C}^{n-k}$ -
by $e_k,...,e_{n-1}$.
\vskip .3cm

\noindent (3.5.9) By definition, the variety $X\i G(k-1,n-1)$ has the rational
parametrization $z\mapsto \gamma (z)$, where $z\i P^{k-1}$ and $\gamma$
 is the logarithmic Gauss map. Explicit formula (3.1.13) gives
 that for generic $z\in {\bf C}^{k-1}$ the subspace $\gamma (z)$
 is the graph of the linear operator
${\bf C}^{k-1}\rightarrow {\bf C}^{n-k}$ given by the matrix
$$  B(z) = \pmatrix{ {a_{k1}z_1\over f_k(z)}&{a_{k2}z_2\over f_k(z)}
&...&{a_{k,k-1}z_{k-1}\over f_k(z)}\cr
{a_{k+1,1}z_1\over f_k(z)}&{a_{k+1,2}z_2\over f_k(z)}&...&
{a_{k+1,,k-1}z_{k-1}\over f_k(z)}\cr
...&...&...&...\cr
{a_{n-1,1}z_1\over f_{n-1}(z)}&{a_{n-1,2}z_2\over f_{n-1}(z)}&...&
{a_{n-1,k-1}z_{k-1}\over f_{n-1}(z)}}\leqno (3.5.10)$$
where $f_j(z) = \sum_{\nu =1}^{k-1} a_{j\nu}z_\nu + a_{jk}$ is the
 equation of the hyperplane $M_j, j=k,...,n-1$.
In other words, the subspace $\gamma (z)$ is spanned by the $(k-1)$ vectors
$e_i + \sum_{j=k}^{n-1} {a_{ji}z_i\over f_j(z)}e_j$.
It is immediate to see that $\gamma(z)$ is the intersection of $(n-k)$
hyperplanes given, in the standard coordinates $y_1,...,y_{n-1}$, by
 linear equations
$$f_j(z) y_j - (a_{j1}z_1) y_1 - ... - (a_{j,k-1}z_{k-1})y_{k-1} =0,
 \,\, j=1,...,n-k.\leqno (3.5.11)$$
The linear functions in (3.5.11), considered together, define a linear
 operator
$a(z): {\bf C}^{n-1}\rightarrow {\bf C}^{n-k}$ whose matrix elements
 are affine functions of $z_1,...,z_{k-1}$.

Let us complete the affine coordinates $z_1,...,z_{k-1}$ in
 ${\bf C}^{k-1}$ to homogeneous coordinates $z_1,...,z_k$ in
 $P^{k-1}$ so that the vanishing of $z_k$ defines the infinite
 hyperplane. Then affine-linear functions $f_j(z_1,...,z_{k-1})$
 will become linear functions $F_j(z) = F_j(z_1,...,z_k) =
\sum _{\nu =1}^k a_{j\nu} z_\nu$.
The $(n-k)$ linear functions in (3.5.11) give rise to a family
of linear operators
$A(z_1,...,z_k): {\bf C}^{n-1} \rightarrow {\bf C}^{n-k}$ given
 by the matrix
$$ \pmatrix{ -a_{k1} z_1& -a_{k2} z_2&...& -a_{k,k-1}
z_{k-1}&F_k(z)&0&...&0\cr
-a_{k+1,1}z_1& -a_{k+1,2}z_2&...&-a_{k+1, k-1}z_{k-1}&0
&F_{k+1}(z)&...&0 \cr
...&...&...&...&...&...&...&...\cr
-a_{n-1,1}z_1 & -a_{n-1,2}z_2 & ...& -a_{n-1, k-1}z_{k-1}&
0&0&...&F_{n-1}(z)}\leqno (3.5.12)$$
whose entries are linear functions in $z_1,...,z_k$. Theorem 3.5.7 is proven.
\vskip .3cm

\noindent {\bf (3.5.13) Remark.} It is immediate to extract from
the formula (3.5.11) the $(n-k)$ subspaces $L_i\i {\bf C}^{n-1}$
whose stars $]L_i[$ are identified
(in the synthetic version (3.5.5) of the Grassmannian
Steiner construction). Namely, $L_i, i=k,k+1,...,n-1$
is the span of the vectors $e_1,...,e_{k-1}, e_i$.

We observe that the position of these subspaces is rather special.
 The identifications of the stars are also very special. The extremely
interesting question
of possibility of Steiner construction of more general Veronese varieties
in Grassmannians will  be treated elsewhere.
\vskip .3cm

\noindent (3.5.14) Varieties in \underbar {projective spaces} defined
 by various generalizations of the Steiner's construction (3.5.4) were
 studied in detail in the book [46] by T.G.Room. To obtain such
a generalization, one takes $r$ subspaces $L_1,...,L_r\i P^{m-1}$
of codimension $d$, identifies all the stars $]L_i[$ with each other
and considers the codimension $r$ subspaces in $P^{m-1}$ which are the
intersections of corresponding hyperplanes from these stars. If $d<r$
then the union of these subspaces is a proper subvariety in $P^{m-1}$
which will be called a {\it projectively generated variety} [46].
The fundamental remark of Room is that any projectively generated
variety in $P^{m-1}$ can be given by a system of equations of which
each equation has the form of a determinant with entries - linear forms
on $P^{m-1}$.

We shall use this idea in the next section to get a better hold of Veronese
 varieties in Grassmannians.
\vskip 1cm

\centerline {\bf (3.6) The sweep of a Veronese variety in Grassmannian.}
\vskip .5cm

\noindent (3.6.1) Let $X$ be any subvariety in the Grassmannian
$G(k-1,n-1)$. So $X$ is a family of $(k-2)$ -dimensional projective subspaces
in $P^{n-2}$. The {\it sweep} of $X$ is, by definition, the subvariety
$Sw(X)\i P^{n-2}$ defined as the union of the subspaces from $X$.
\vskip .3cm

\noindent (3.6.2) We shall be mostly interested in the case when
 $X \i G(k-1,n-1)$ be a ($k-1$ -dimensional) special
Veronese variety, see (3.3.9). In other words (3.3.10), $X$
is the visible contour of a Lie complex $Z=\overline{H.L}$ i.e.
 the locus of $P^{k-1}$'s in $P^{n-1}$ which belong to the complex
 $Z$ and contain the chosen
point $e=(1,...,1)$. The sweep of $X$ is just what we called in
(3.1.2) the visible sweep of the complex $Z$ at $e$. So it is
 the projectivization of the cone in $P^{n-1}$ with vertex $e$
given by the union of all $P^{k-1}$'s from the complex $Z$ which
 contain $e$.

\vskip .3cm

\noindent (3.6.3) Let  ${\bf h} = {\bf C}^n / \{(a,...,a)\}$
be the Lie algebra of the maximal torus $H\i PGL(n)$. Recall
 (3.1.8) that the Grassmannian $G(k-1,n-1)$ in which the visible
 contours (and hence special Veronese varieties) lie,
is in fact $G(k-1,{\bf h})$. Therefore the sweep of any special Veronese
variety lies naturally in the projective space $P({\bf h})$
\vskip .3cm

\noindent (3.6.4) Let $t_1,...,t_n$ be standard coordinate functions
on ${\bf C}^n$. A linear form
$\sum c_i t_i$ descends to a linear form on {\bf h} if $\sum c_i=0$.
 In particular, the {\it roots} i.e. the linear forms $t_i-t_j$ are
 forms on {\bf h}.

\vskip .3cm

\noindent (3.6.5) The possibility of defining $X$ by the Grassmannian
Steiner construction (3.5.5) implies that $Sw(X)$ is always projectively
generated variety in the sense of (3.5.14).
\vskip .3cm

\proclaim (3.6.6) Theorem. Suppose that $k\leq n-k$. Let $z_1,...,z_{k}$
be homogeneous coordinates in $P^{k-1}$.
Suppose that a configuration $M=(M_1,...,M_n)$ of hyperplanes in
$P^{k-1}$ consists of $k$ coordinate hyperplanes $M_i = \{z_i=0\},
 i=1,...,k$ and $(n-k)$ other hyperplanes
$M_j = \{ \sum_{i=1}^k a_{ji}z_j =0\}$.
Then the sweep of the Veronese variety in $P^{n-2} = P({\bf h})$
corresponding to $M$ is given by vanishing of all $k\times k$ -minors
 of the following $k\times (n-k)$ -matrix
of linear forms on {\bf h}:
$$ A^\dagger (t_1,...,t_n) = ||a_{ji}(t_j-t_i)||,
i=1,...,k, \,\, j= k+1,...,n.
\leqno (3.6.7)$$

\noindent {\sl Proof:}
Let $A:{\bf C}^k \rightarrow Hom( {\bf C}^{n-1}, {\bf C}^{n-k})$
 be the linear system of linear operators such that $X$
consists of kernels of $A(z), z\in {\bf C}^k-\{0\}$.
An explicit formula for $A$ is given in (3.5.12).
Using partial dualization, let us associate to $A$ a linear operator
$$A^\dagger: {\bf C}^{n-1} \rightarrow Hom ({\bf C}^k, {\bf C}^{n-k})$$
A point $t\in {\bf C}^{n-1}$ lies in the kernel of $A(z)$ for some non-zero
$z\in {\bf C}^k$ if and only if  the linear operator $A^\dagger (t):
{\bf C}^k\rightarrow {\bf C}^{n-k}$ has non-trivial kernel i.e. the rank of
$A^\dagger (t)$ is less than $k$. Thus the sweep $Sw(X)$ is
 defined by vanishing of all $k\times k$ minors of the matrix
 $A^\dagger (t)$ of linear forms on $P^{n-2}$.

To see  that $A^\dagger$ has the claimed form, we use the formula
 (3.5.12).
This formula was written with respect to the non-symmetric  system
 of coordinates $y_1,...,y_{n-1}$
in {\bf h}. In the language of (3.6.4) we have $y_i = t_i-t_n$. Substituting
this to (3.5.12) and transposing the matrix, we arrive at the formula (3.6.7).
 Theorem 3.6.6 is proven.

\vskip .3cm

\noindent (3.6.8) Let us describe a more geometric construction for
the sweep of the Veronese variety corresponding to a projective configuration.

Let $Mat(k,n-k)$ be the vector space of all $k$ by $(n-k)$ matrices
 and   $P(Mat(k,n-k))$ be the projectivization of this space.
 The projectivization of the set of matrices of rank 1 is just
 the Segre embedding $P^{k-1}\times P^{n-k-1} \i P(Mat (k,n-k))$.
Let $\nabla \i P(Mat(k,n-k))$ be the projectivization of the space
 of matrices of rank $<k$. This is an algebraic subvariety of codimension
$n-2k+1$.

\vskip .3cm

\noindent (3.6.9) Let now $(x_1,...,x_n)$ be a configuration of
 \underbar {points} in $P^{k-1}$ in general position. (Recall that
 modulo projective isomorphism,  configurations of points give the same
 orbit space as configurations of hyperplanes.) Let $(y_1,...,y_n)$ be
 the  configuration of points in $P^{n-k-1}$
associated to $x_1,...,x_n$ (see \S 2.3 about association).
 By Reformulation (2.3.8), the points
$z_i=(x_i,y_i) \i P^{k-1} \times P^{n-k-1} \i P(Mat(k,n))$
 form a circuit i.e.  span a projective space, say, $L$ whose dimension
 is $n-2$  and are in general position as points of $L$. The space
 $P^{n-2}_e = P({\bf h})$ also comes with a circuit given by points
$\bar e_i$,-- projectivizations of images of the basis vectors
$e_i\in {\bf C}^n$ in ${\bf h} = {\bf C}^n/{\bf C}$. Hence there is
 a unique  projective transformation $\phi : L\rightarrow P({\bf h})$
 taking $z_i$ to $\bar e_i$. We shall be interested in the intersection
 $\nabla\cap L\i L$ where $\nabla$ is the determinantal variety in (3.6.8).
\vskip .3cm

\proclaim (3.6.10) Proposition. The map $\phi$ identifies
the subvariety $\nabla\cap L\i L$ with the sweep
 $Sw(X (x_1,...,x_n)) \i P({\bf h})$ of the Veronese
variety corresponding to the configuration $(x_1,...,x_n)$.

\noindent {\sl Proof:} Denote the sweep $Sw (X (x_1,...,x_n))$
 shortly by $S$. Let $M_1,...,M_n$ be the hyperplanes in the dual
 $P^{k-1}$ corresponding to $x_i$. After choosing suitable homogeneous
 coordinates we can apply Theorem
(3.6.6) which gives a representation of $S$ as the inverse image
 of $\nabla$
under the linear embedding $A^\dagger : {\bf h} \rightarrow Mat (k,n-k)$.
We regard $A^\dagger$ as a map ${\bf C}^n \rightarrow Mat (k,n-k)$ using
the isomorphism ${\bf h} = {\bf C}^n /{\bf C}$. Let $e_i\in {\bf C}^n$
 be the standard basis vectors. Proposition (3.6.10) is a consequence
 of the following statements:
\vskip .3cm

\noindent (3.6.11) The matrices $A^\dagger (e_i)$ lie in the Segre
 embedding $P^{k-1}\times P^{n-k-1} \i P( Mat(k,n-k))$ i.e.,
 $rank A^\dagger (a_i) =1$.

\vskip .3cm

\noindent (3.6.12) The configuration of hyperplanes
 $Ker (A^\dagger (e_i)) \i P^{k-1}$ is projectively isomorphic
 to $(M_1,...,M_n)$ and the configuration of points
$Im (A^\dagger (e_i)) \i P^{n-k-1}$ is associated to $(M_1,...,M_n)$.

\vskip .3cm

Both these statements are immediate from the explicit form (3.6.7)
of the matrix $A^\dagger$.
\vskip .3cm

\proclaim (3.6.14) Corollary. Let $n=2k$. Then Veronese varieties
 corresponding to a configuration $(M_1,..., M_{2k})\i P^{k-1}$ and
to the associated configuration, have the same sweep.

\vskip .3cm

\noindent (3.6.15) Any determinantal variety i.e. variety defined by
vanishing of monors of a matrix of linear forms, bears two canonical
 families of projective subspaces, so-called $\alpha$- and $\beta$ -
families [46]. Let us recall their construction and explain their
 relevance to our situation.

Let $k\leq n-k$ and $\nabla\i P(Mat (k,n-k))$ denote, as before,
 the projectivization of the space of matrices of rank $<k$. For
any 1-dimensional subspace
$\lambda \i {\bf C}^k$ set
$$\Pi_\alpha (\lambda) = P(\{M:{\bf C}^k \rightarrow {\bf C}^{n-k}
 : M(\lambda) =0\})\i\nabla. \leqno (3.6.16)$$
This is a projective subspace in $\nabla$ of codimension $k-1$. Thus
 we get a family of projective subspaces in $\nabla$ (called $\alpha$
 -subspaces) of codimension $k-1$, parametrized by $P^{k-1} = P({\bf C}^k)$.

Similarly, for any hyperplane $\Lambda \i {\bf C}^{n-k}$ set
$$\Pi_\beta (\Lambda) = P(\{M:{\bf C}^k \rightarrow {\bf C}^{n-k}:
 Im(M)\i \Lambda \})\i\nabla .\leqno (3.6.17)$$
This is a projective subspace in $\nabla$ of codimension $n-k-1$.
We get a family of projective subspaces in $\nabla$ (called $\beta$
 -subspaces) of codimension $n-k-1$ parametrized by $P^{n-k-1} =
 P(({\bf C}^{n-k})^*)$.
\vskip .3cm

\noindent (3.6.18) Let $L\i P(Mat(k,n-k))$ be a projective subspace
 of dimension $n-2$. Consider the variety $S=L\cap\nabla$. It contains
projective subspaces
$\Pi_\alpha (\lambda) \cap L$ whose dimension is at least
$n-k-1$ (they will be called the $\alpha$ -subspaces in $S$)
 and subspaces $\Pi_\beta (\Lambda)\cap L$ whose dimension is
at least $k-1$ (they will be called $\beta$ -subspaces in $S$).
 The role of this subspaces in our situation is as follows.

\vskip .3cm

\proclaim (3.6.19) Proposition. Let $(x_1,...,x_n)$ be a configuration
of points in $P^{k-1}$ in general position and $(y_1,...,y_n)$ -
the associated configuration in $P^{n-k-1}$. Let $X\i G(k-1,n-1),
 X ' \i G(n-k,n-1)$ be the Veronese varieties corresponding to
$(x_i)$ and $(y_i)$ (their dimensions equal, respectively, $k-1$
 and $n-k-1$). Let $S,S'\i P^{n-2}_e = P({\bf h})$ be the sweeps
 of these varieties. Then $S=S'$ and the subspaces from $X$ (resp.
 from $X '$) lying on $S=S'$ are precisely the $\beta$- (resp.
$\alpha$-) subspaces on $S$ defined in (3.6.18).

\noindent {\sl Proof:} This is a reformulation of Theorem 3.5.7
about Steiner construction of $X$.

\beginsection (3.7) An example: Visible contours and sweeps of Lie
complexes in $G(3,6)$.

\vskip .5cm

In this section we study in detail the construction of \S 3.6 in the
 particular case corresponding to configuration of 6 points on $P^2$.
In other words, we consider the case $k=3, n=6$.
\vskip .3cm

\noindent (3.7.1) To any sextuple $(x_1,...,x_6)\in (P^2)^6_{gen}$ a
Veronese surface  $X(x_1,...,x_6)\i G(2,5)$ is associated. Its sweep,
denoted
$S(x_1,...,x_6) \i P^4$ is a cubic hypersurface since by Theorem 3.6.6
it is given by vanishing of the determinant of a 3 by 3 matrix of
linear forms.
We shall study such hypersurfaces.
\vskip .3cm

\noindent (3.7.2) We are interested in configurations modulo
 projective isomorphism. So we can consider equally well the
sixtuple of lines $M_i\i \check P^2$ dual to $x_i$. This sextuple
 represents the same element of $(P^2)^6/GL(3)$.

\vskip .3cm

\noindent (3.7.3) We can always assume that lines $M_i$ have
the particular form considered in Theorem 3.6.6: for i=1,2,3
 the line $M_i$ is given by the equation $z_i=0$ and $M_j$ for
 $j=4,5,6$ is given by the equation $\sum a_{ij}z_j=0$.
 The $3\times 3$ matrix $||a_{ij}||$ is defined by a projective
 isomorphism class of $(M_1,...,M_6)$ not uniquely but only up to
 multiplication of rows and columns by non-zero scalars. Generic
position of lines $L_\nu$ implies that all $a_{ij}\neq 0$. Hence
by multiplication of rows and columns by scalars we can take
 $||a_{ij}||$ into a unique matrix of the form
$$\pmatrix{ 1&1&1\cr
1&a&b\cr
1&c&d}. \leqno (3.7.4)$$
In this way the quotient $(P^2)^6_{gen}/GL(3)$ becomes identified
with the space of $(a,b,c,d)$ such that all the minors of the matrix
(3.7.4) are non-zero.
\vskip .3cm

\proclaim (3.7.5) Proposition. The points $x_1,...,x_6$ lie on a
conic (or, equivalently, the lines $M_1,...,M_6$ are tangent to a
 conic) if and only if the matrix elements $a,b,c,d$ of the matrix
(3.7.4) satisfy the equation $\Psi(a,b,c,d)=0$, where
$$\Psi(a,b,c,d)= det \pmatrix{a(1-c)&b(1-d)\cr c(1-a)&d(1-b)} =
ad-bc+abc+bcd-acd-abd.\leqno (3.7.6)$$

\noindent {\sl Proof:} This is Proposition 2.13.1 of [41].

\vskip .3cm

\noindent (3.7.7) The Veronese variety $X = X(x_1,...,x_6)$ lies
in $G(2,5)$, the space of lines in $P^4= P({\bf h}), {\bf h} =
 {\bf C}^6/{\bf C}$. Let
$p_i\in P^4, i=1,...,6$ be the point corresponding to the standard
basis vector $e_i\in {\bf C}^6$ (In the realization of $P^4$ as the
 space of lines in $P^5$ through $e=(1,...,1)$, the point $p_i$
 corresponds to the line $<e,e_i>$. By corollary 3.4.7, the variety
$X$ contains all the lines $<p_i,p_j>$. This implies that points $p_i$
 are singular points of the sweep $S=S(x_1,...,x_6)$ of $X$. Indeed,
the tangent directions at any $p_i$ to lines $<p_i,p_i>, j\neq i$,
span the whole tangent space $T_{p_i}P^4$ which will therefore
coincide with the Zariski tangent space of $S$ at $p_i$.

A more detailed information about the singularities of $S$ is given
 in the next proposition which is the main result of this section.
\vskip .3cm

\proclaim (3.7.8) Proposition. a) If $x_1,...,x_6\in P^2$ are in
general position and do not  lie on a conic then the sweep
$S(x_1,...,x_6)$ has only 6 singular points namely $p_i$ and
these points are simple quadratic singularities.\hfill\break
b) If $x_1,...,x_6$ are in general position and lie on a conic
$K\i P^2$ then $S(x_1,...,x_6)$ has a curve $C$ of singular points
 which is a Veronese curve in $P^4$ containing $p_1,...,p_6$.
 In this case the configuration of $x_i$ on $K\cong P^1$ is
isomorphic to that of $p_i$ on $C$. The variety $S(x_1,...,x_6)$
is the union of all chords of $C$.

\noindent {\sl Proof:} Consider the varieties
$$P^2\times P^2\i \nabla \i P^8 = P(Mat(3,3))\leqno (3.7.9)$$
where $\nabla$ consists of degenerate matrices and $P^2\times P^2$ -
of matrices of rank 1. It is well-known that $P^2\times P^2 = Sing (\nabla)$.

Assume that the configuration $(M_1,...,M_6)$ of lines dual to $x_i$
 has the form spacified in Theorem 3.6.6 with the $3\times 3$ matrix
$||a_{ij}||$ given in the normal form (3.7.4). Let $A^\dagger :
 P^4=P({\bf h}) \rightarrow P(Mat(3,3))$ be the embedding given
 by formula (3.6.7). In other words (taking into account the normal
 form of $||a_{ij}||$) we have
$$A^\dagger (t_1,...,t_6) = \pmatrix{t_1-t_4&t_1-t_5&t_1-t_6\cr
t_2-t_4&a(t_2-t_5)&b(t_2-t_6)\cr
t_3-t_4&c(t_3-t_5)&d(t_3-t_6)}\leqno (3.7.10)$$
Theorem 3.6.6 implies that our sweep $S$ equals $(A^\dagger)^{-1}(\nabla)$.
\vskip .3cm

\noindent (3.7.11) Let $L\i P(Mat(3,3))$ be the image of $A^\dagger$.
It is immediate to check that the degree of Segre variety $P^2\times
 P^2\i P^8$ equals 6.
 On the other hand, the subspace $L\i P^4$ already intersects
$P^2\times P^2$ in 6 points $q_i=A^\dagger (e_i), i=1,...,6$.
Hence there remains one of two possibilities:
\vskip .2cm
\noindent \underbar {Case 1.} $L$ intersects $P^2\times P^2$
transversally in 6 points $q_i = A^\dagger (e_i)$.

\vskip .2cm

\noindent \underbar {Case 2.} The intersection
$L\cap (P^2\times P^2)$ contains a component of positive dimension.

Part a) of Proposition 3.7.8 will follow from the next two lemmas.
\vskip .3cm

\proclaim (3.7.12) Lemma. If $x_1,...,x_6$ are in general position
 and do not lie on a conic then for $L=A^\dagger (P({\bf h}))$ the
Case 1 holds.

\vskip .3cm

\proclaim (3.7.13) Lemma. If for $L$ the Case 1 holds then $L\cap\nabla$
has $q_i$ as the only singular points and the singularities at
 $q_i$ are simple quadratic.
\vskip .3cm

\noindent (3.7.14) \underbar{\sl Proof of Lemma 3.7.13:} The only
 possibility which we have to exclude is that there is a point
$q\in L\cap\nabla$ which is smooth on $\nabla$ (i.e. $rank (q) =2$)
 and  such that $T_qL\i T_q\nabla$. Do rule out this possibility,
let $\lambda\in P^2$ be the point corresponding to $Ker(q)$ and let
 $\Pi = \Pi_\alpha (\lambda)$ be the corresponding $\alpha$- subspace
i.e. the projectivization of the space of all $3\times 3$ matrices
annihilating $\lambda$. Then $dim (\Pi) =5$. Since $L$ is connected
in the embedded tangent space to $\nabla$ at $q$ (which is
7-dimensional), we have $dim(L\cap\Pi) \geq 2$. However, by
Proposition 3.6.19, the intersection of $L$ with all $\alpha$
- subspaces should be 1-dimensional.
\vskip .3cm

\noindent (3.7.15) \underbar {\sl Proof of Lemma 3.7.12:} It
suffices to show that each $q_i$ is a isolated singular point
 of the intersection $L\cap (P^2\times P^2)$. To this end, we
prove that the tangent spaces $T_{q_i}L, T_{q_i}(P^2\times P^2)
 \i T_{q_i}P^8$ intersect only in 0. Since the roles of $q_i$
are symmetric, it is enough to consider $i=1$. The point $q_1
 = A^\dagger (e_1)$ is given, in virtue of the normal form (3.6.27)
 of $A^\dagger$, by the matrix
$$\pmatrix{1&1&1\cr 0&0&0\cr0&0&0}.$$
The tangent space at $q_1$ to the locus of rank 1 matrices is easily
 seen to consist of matrices of the form
$$\pmatrix{\lambda_1&\lambda_2&\lambda_3\cr
\lambda_4&\lambda_4&\lambda_4\cr
\lambda_5&\lambda_5&\lambda_5}.$$
Therefore the intersection $T_{q_1}L\cap T_{q_1}(P^2\times P^2)$ is obtained
from the space of solutions $(t_1,...,t_6)$ of the linear system
$$t_2-t_4 = a(t_2-t_5) = b(t_2-t_6),\quad t_3-t_4 = c(t_3-t_5) =
d(t_3-t_6) \leqno (3.7.16)$$
by factorization by the 2-dimensional subspace spanned by $e_1=
(1,0,...,0)$ and $(1,1,...,1)$. Hence we are reduced to the
following statement:
\vskip .3cm

\proclaim (3.7.18) Lemma. If $x_i$ do not lie on a conic (i.e.
if the polynomial $\Psi (a,b,c,d)$ given by (3.7.6) does not vanish)
then the linear system (3.7.16) has 2-dimensional space of solutions.

\noindent {\sl Proof:} This is a system of four equations on 6
 variables whose matrix of coefficients has the form
$$\pmatrix { 0&1-a&0&-1&a&0\cr
0&1-b&0&-1&0&b\cr
0&0&1-c&-1&c&0\cr
0&0&1-d&-1&0&d}.$$
Let us disregard the first column, then move the column with
$(-1)$'s to the left and then subtract the first row from all
 the other rows. We obtain the matrix
$$\pmatrix{ -1&1-a&0&a&0\cr
0&a-b&0&-a&b\cr
0&-1+a&1-c&c-a&0\cr
0&-1+a&1-d&-a&d}.$$
It is immediate to see that all 4 by 4 minors of this matrix  have
 the form $\pm \Psi (a,b,c,d)$.

We have proven part a) of Proposition 3.7.8.

\vskip .3cm

\noindent (3.7.19) Let us prove part b) of Proposition 3.7.8. So assume
that $x_1,...,x_6\in P^2$ are points in general position lying on a
 conic $K$. Consider the 2-fold Veronese embedding
$$v_2:P^2\hookrightarrow P^5 = P(S^2({\bf C}^3)) \i P(Mat (3,3))
\leqno (3.7.20)$$
where $P(S^2({\bf C}^3))$ is embedded into $P(Mat(3,3))$ as the
 space of symmetric matrices. Let $L\i P(S^2({\bf C}^3))$ be the
 projective envelope of $v_2(x_i)$. Let also $\nabla_{sym} \i
 P(S^2({\bf C}^3))$ be the space of degererate quadratic forms
i.e. $\nabla_{sym} = \nabla \cap  P(S^2({\bf C}^3))$. Note
 that $P^2 = Sing (\nabla_{sym})$ has codimension 2 in $\nabla_{sym}$.

Since $x_i$ lie on a conic, their configuration is self-associated
 (Example 2.3.12). Now the interpretation of $S(x_1,...,x_6)$ given
 in (3.6.9), (3.6.10) implies that $S(x_1,...,x_6) = L\cap \nabla_{sym}$.

The conic $K$ is equal to $v_2^{-1}(L\cap v_2(P^2))$. Since $x_i$
are in general position, $K$ is smooth. Now since $P^2 = Sing
 (\nabla_{sym})$, we find that $C = v_2(K) = Sing (L\cap\nabla_{sym})$
 is the singular curve of $L\cap\nabla_{sym} = S(x_1,...,x_6)$. This
 is clearly a Veronese curve in $L=P^4$. The Veronese embedding $v_2$
 indetifies $k$ with $C$ and points $x_i\in K$ with our distinguished
 points $p_i\in C$. Finally, $\nabla_{sym}$ is the union of chords of
 $v_2(P^2)$ (every degenerate quadratic form in 3 variables is a sum
 of two quadratic forms of rank 1). Hence $S(x_1,...,x_6)$ contains
 the union of chords of $C$. Since the latter is also a 3-dimensional
 variety, the two varieties in question
coincide.

Proposition 3.7.8 is completely proven.
\vskip .3cm

\noindent (3.7.21) Fix 6 points $p_1,...,p_6\in P^4$ in general
 position. Any two choices of $p_i$ can be taken to each other
 by a unique projective isomorphism.
Let ${\cal L}\i P(S^3({\bf C}^5))$ be the linear system of all
 cubic hypersurfaces in $P^4$ which contain $p_i$ as singular
points. It is clear by dimension count that ${\cal L}$ has
dimension 4. On the other hand, taking $p_i$ to be the standard
 points (images of $e_i$ in $P({\bf h})$ we have constructed a
 4-dimensional family of sweeps $S(x_1,...,x_6)$ which all belong
 to ${\cal L}$. Hence we have the following corollary.
\vskip .3cm

\proclaim (3.7.22) Corollary. Let $p_1,...,p_6\in P^4$ be points
in general position. A generic cubic hypersurface $S\i P^4$ for
 which $p_i$ are singular points, is projectively equivalent to
the visible sweep of some Lie complex in $G(3,6)$. In particular,
 $S$ can be realized as $P^4\cap \nabla$ for a suitable embedding
 $P^4\i P(Mat(3,3))$. The variety $S$ contains two families of
 lines ($\alpha$ and $\beta$- families, in the determinantal
realization), whose parameter spaces $P,P'$ are isomorphic to
$P^2$. These families give rise to two Veronese surfaces in
$G(2,5)$ which correspond to a pair of associated configurations
$(M_1,...,M_6)\i P, (M'_1,...,M'_6)\i P'$ of lines. Explicitly,
 $M_i\i P$ (resp. $M'_i\i P'$) is the locus of all lines from
the first (resp. second) family which contain the point $p_i$.

\vskip .3cm

\noindent {\bf (3.7.23) Remarks.} a) Any generic intersection
$P^4\cap \nabla \i P(Mat(3,3))$ intersects the Segre variety
$P^2\times P^2 = Sing (\nabla)$ in 6 points (since
$6= deg (P^2\times P^2)$) and hence is a cubic hypersurface
 of the form studied in the above corollary.

b) The correspondence
$$(x_1,...,x_6) \quad mod \quad PGL(3) \quad \longmapsto\quad
 S(x_1,...,x_6)$$
is two-to-one. Hence it defines a two-sheeted covering of fourfolds
$\pi: (P^2)^6_{gen}/GL(3) \rightarrow {\cal L}$. This covering
is well-known classically, see [14,41]. It extends to a map of
the Mumford quotient
$((P^2)^6/GL(3))_{Mumf} \rightarrow {\cal L}$ which is a double
 cover ramified along a hypersurface $W\i {\cal L}$ of degree 2.
This hypersurface is called the modular variety of level 2 (see
 [14]) The projective dual $\check W \i \check {\cal L}$ is a
so-called Segre cubic threefold i.e. a cubic hypersurface with 10
ordinary singular points. (It is known that all such threefolds
 are projectively isomorphic.)

\vskip 1cm

\centerline {\bf (3.8) Chordal varieties of Veronese curves.}
\vskip 1cm

\noindent (3.8.1) Let $C\i P^r$ be a Veronese curve. An $s-1$ -
dimensional projective subspace $L\i P^r$ is called {\it chordal}
 to $C$ if it intersects $C$ in $s$ points (counted with multiplicities).
 Denote by $Ch_{s-1}(C)$ the variety of all chordal $(s-1)$ -dimensional
 subspaces of $C$. Clearly $Ch_{s-1}(C)$ is isomorphic to the $s$ -fold
 symmetric power of $C\cong P^1$ i.e. to the projective space $P^s$.
 We obtain therefore a special class of embeddings $P^s\i G(s, r+1)$.
 It turns out that these embeddings give particular cases of Veronese
 varieties in Grassmannians considered in \S 3.3 above.
\vskip .3cm

\proclaim (3.8.2) Proposition. The chordal variety of any Veronese
curve is a Veronese variety in the Grassmannian.

This proposition is classical and due to L.M.Brown [9]. In modern
language it is a consequence (or, rather, a reformulation) of the
following fact.
\vskip .3cm

\proclaim (3.8.3) Proposition.  There is a unique isomorphism of
$GL(2)$- modules
$$\xi: \bigwedge^k(S^n({\bf C}^2)) \quad \longrightarrow \quad
S^{n-k-1}(S^k({\bf C}^2)) \otimes (\bigwedge^2({\bf C}^2))^{\otimes
 k(k-1)/2}$$
such that for any $l_1,...,l_k\in {\bf C}^2$,
$$\xi (l_1^n\wedge ...\wedge l_k^n) \quad = \quad (l_1...l_k)^{n-k-1}
 \otimes \prod_{i<j} (l_i\wedge l_j). \quad\triangleleft$$

\vskip .3cm

\noindent (3.8.4) Put now $s=k-1, r=n-2$ so that we obtain Veronese
 varieties in $G(k-1,n-1)$ i.e., are in the setting of sections (3.3)
 and (3.4). Recall that special Veronese varieties in $G(k-1,n-1)_e$
 (i.e. in the space of $P^{k-1}$'s in $P^{n-1}$ through $e=(1,...,1)$)
 are in bijection with the set of projective equivalence classes of
 $n$ -tuples of points in $P^{k-1}$ in general position. Let us
 clarify the place in this picture of chordal varieties of Veronese curves.

Suppose that $n$ distinct points $x_1,...,x_n\in P^{k-1}$ happen
to lie on a Veronese curve $D$ (of degree $k-1$) in $P^{k-1}$.
 Then they are in general position in $P^{k-1}$, as it follows
 from the calculation of the Vandermonde determinant. Thus they
represent an element of $(P^{k-1})^n_{gen}/GL(k)$. Such element,
 by Corollary 3.3.11, is represented by a unique special Veronese
 variety $X(x_1,...,x_n)\i G(k-1,n-1)_e$ of dimension $k-1$.
 On the other hand, the curve $D$ being isomorphic to $P^1$,
the points $x_i$ represent an element from $(P^1)^n_{gen}/GL(2)$.
 The latter set, as we have seen in Example 3.4.9, is identified
 with the set of Veronese curves in $P^{n-2}_e = G(1,n-1)_e$
 through $n$ points $<e,e_i>$. Let $C(x_1,...,x_n)$ be the
special Veronese curve representing the configuration of $x_i$
 on $D$.
\vskip .3cm

\proclaim (3.8.5) Theorem.  The special Veronese variety
 $X(x_1,...,x_n)$ coincides with the chordal variety of the
 Veronese curve $C(x_1,...,x_n)$.

Using the language of hyperplane configurations, this can be
reformulated as follows.
\vskip .3cm

\proclaim (3.8.6) Reformulation. Let $D$ be a Veronese curve
in $P^{k-1}$, $M=(M_1,...,M_n)\i P^{k-1}$ -a configuration of
 hyperplanes which are osculating to $D$
(i.e. each $M_i$ intersects $D$ in just one point $x_i$ with
muptiplicity $(k-1)$). Then the embedding $\phi_E:P^{k-1}
\hookrightarrow G(k-1,n-1)$ defined by the vector bundle
$E=\Omega^1_{P^{k-1}}({\rm log} M)$ maps $P^{k-1}$
isomorphically to the chordal variety of some other Veronese
curve $C\i P^{n-1}$. The curve $C$ is the image of $D$ in the
 projective embedding defined by the line bundle $\Omega^1_D
({\rm log} (x_1+...+x_n))$.

\noindent {\sl Proof:} An easy calculation in coordinates shows
 that the restriction of 1-forms defines an isomorphism
$$H^0(P^{k-1}, \Omega^1_{P^{k-1}}({\rm log} M)) \longrightarrow
 H^0(D, \Omega^1_D({\rm log} (x_1+...+x_n)).\eqno (3.8.7)$$
Denote for short the bundle $\Omega^1_{P^{k-1}}({\rm log} M)$
 on $P^{k-1}$  by $E$ and the bundle $\Omega^1_D({\rm log}
(x_1+...+x_n)$ on $D$ - by $F$. Let $\phi_E$ and $\phi_F$ be
the corresponding maps to the Grassmannian and the proejctive
space respectively.
We shall show that under the identification (3.8.7) the image
of $\phi_E$ is the chordal variety of the image of $\phi_F$.
 By definition of $\phi_E, \phi_F$ (see section (3.2) above)
 this is equivalent to part a) the following statement:
\vskip .3cm

\proclaim (3.8.8) Lemma. Let $p\in P^{k-1}$ be a generic point.
 Then there are $(k-1)$ points $y_1(p), ..., y_{k-1}(p) \in D$
 such that a form $\omega\in H^0(P^{k-1},E)$ vanishes at $p$
(as a section of $E$) if and only if the restriction of $\omega$
 to $D$ (as a 1-form) vanishes at all $y_i$. \hfill\break
b) Explicitly,  points $y_i(p)$
are precisely the points of osculation of the $k-1$ osculating
 hyperplanes to $D$ passing through $p$

 Note that through any generic point $p\in P^{k-1}$ there pass
exactly $k-1$ osculating hyperplanes to $D$. (Since osculating
 hyperplanes to $D$ form a Veronese curve $\hat D$ in the dual
 projective space, this just means that the degree of $\hat D$
 is also $k-1$).

 We shall verify Lemma 3.8.8 in coordinates.
\vskip .3cm

\noindent (3.8.9) Let us regard the affine space ${\bf C}^k$ as
the space of polynomials $f(t) = \sum _{i=0}^{k-1} a_it^i$ of
degree $\leq k-1$ in one variable. The hyperplanes
$M_i$  have the form $M_i = \{f: f(x_i)=0\}$, where $x_i\in
{\bf C}$ are distinct numbers. The equation of $M_i$ is,
 therefore, the evaluation map $f\mapsto f(x_i)$. Any point
$p\in P^{k-1}$ is represented  as a polynomial $f(t)$. If $\alpha_i$
 are the roots of $f(t)$ then the points $y_i(p)$ are the polynomials
 $(t-\alpha_i)^{k-1}$. We can normalize any polynomial $f$ to have
 the form $f(t) = \prod (t-\alpha_i)$.

A section of the bundle $E$ is given as the logarithmic differential
of the function
$$f\mapsto \prod_{j=1}^n f(x_j)^{\lambda_j} = \prod_{i=1}^{k-1}
\prod_{j=1}^n (x_j - \alpha_i)^{\lambda_j},\quad \, \sum \lambda_i=0.$$
Our assertion means that this function has a critical point at a
 given $f$ if and only if the function $\prod (t-x_j)^{\lambda_j}$
 has a critical value at each $t=\alpha_i$. But this follows from
the equality
$${\partial\over\partial t} {\rm log} \prod_j
(t-x_j)^{\lambda_j}|_{t=\alpha_i} =  - {\partial\over\partial \alpha_r}
{\rm log} \prod_{j} (x_j-\alpha_i)^{\lambda_j}
 = \sum_j {\lambda_j\over \alpha_i -t_j}. $$
Lemma 3.8.8 and hence Theorem 3.8.5 are proven.

\vskip .3cm

\noindent {\bf (3.8.10) Examples.} a) Consider the case of 5 points
 in $P^2$ i.e. $k=3, n=5$. Such configurations lead to Veronese
 surfaces in $G(2,4)\i P^5$. Since every 5 points in $P^2$ lie on
a unique conic, any special Veronese surface in $G(2,4)$ will be
 the chordal variety of a Veronese curve (twisted cubic) in $P^3$.

b) Similarly, the case of $n$ points in $P^{n-3}$ leads to Veronese
 $(n-3)$ -folds in $G(n-3,n-1)$ which are chordal varieties of
Veronese curves in $P^{n-2}$. Note that $G^0(k,n)/H = G^0(n-k,n)/H$
 by duality and hence the case of $n$ points in $P^{n-3}$ is equivalent
 to that of $n$ points on $P^1$ (see  section (2.3)). Any $n$ points
in $P^{n-3}$ in general position lie on a unique Veronese curve which
 provides the dual configuration.

c) For the case $k=3,n=6$ (six points in $P^2$) we associate to
 sextuples $(x_1,...,x_6)\in (P^2)^6_{gen}$ Veronese surfaces
 $X(x_1,...,x_6)$ in $G(2,5)$, the space of lines in $P^4$.
 When $x_1,...,x_6$ lie on a conic, the surface $X(x_1,...,x_6)$
 is the chordal surface of a Veronese curve in $P^4$. We have
seen this is proposition 3.7.8.

\vskip .3cm

\noindent {\bf (3.8.11) Remark.} It is a remarkable fact that
chordal varieties of Veronese curves (regarded as subvarieties
 in Grassmannians) posess deformations which do not come from
 Veronese curves at all (and represent general projective configurations).

\vskip 1cm

\centerline {\bf (3.9) The homology class of a (special)
Veronese variety in Grassmannian.}

\vskip .6cm

We have associated to each isomorphism class of  generic
configurations of $n$ points in $P^{k-1}$ a certain embedding
of $P^{k-1}$ into the Grassmannian $G(k-1,n-1)$ --- the
Veronese variety. For instance, configurations of points on
$P^1$ correspond to Veronese curves in $P^{n-2}$ through a
 fixed set of $n$ generic points. In this section we calculate
 the homology class $\Delta$ represented by these Veronese varieties
in $G(k-1,n-1)$. It turns out that the coefficients of the expansion
 of $\Delta$ in the basis of Schubert cycles are exactly the dimensions
of irreducible representations of the group $GL(n-k)$.

\vskip .3cm

\noindent (3.9.1) Let us review the homology  theory of the
 Grassmannian $G(p, q)= G(p, {\bf C}^q)$ (see [24] for more
details). Let $\alpha = (\alpha_1 \geq ...\geq \alpha_p\geq 0),
 \alpha_i\leq q-p$, be a decreasing sequence of non-negative
integers. We visualize $\alpha$ as a Young diagram in which
 $\alpha_i$ are the lengths of rows. Because of inequalities
 $\alpha_i\leq q-p$ the diagram $\alpha$ lies inside the
rectangle $[0,q-p]\times [0,p]\i {\bf R}^2$:

\vbox to 4cm{\vfill\vfill (3.9.2) \vfill\vfill}

We define $|\alpha| = \sum \alpha_i$ to be the number of cells in $\alpha$.

\vskip .3cm

\noindent (3.9.3) To each Young diagram $\alpha$ as above we
 associate the lattice path $\Lambda (\alpha) \i [0,q-p]\times
 [0,p]$ going from the points $(0,0)$ to the point $(q-p,p)$.
This path is just the boundary of $\alpha$, see Fig. (3.9.2).
 It consists of exactly $q$ edges $E_1,...,E_q$ which are
horizontal or vertical segments of the lattice. We write $E_i$
 in such order that $E_1$ begins at $(0,0)$ and $E_q$ ends at
 $(q-p,p)$. The number of vertical edges is $p$ and a path is
 completely determined by specifying which of the $E_i$ are
vertical. So the numbers of possible lattice paths in $[0,q-p]
\times [0,p]$ equals $q\choose p$. The same will be the number
 of all Young diagrams in $[0,q-p]\times [0,p]$,
if we count also the "empty" diagram $(0,...,0)$.
\vskip .3cm

\noindent (3.9.4) Let $\alpha \i [0,q-p]\times [0,p]$ be a
Young diagram and $\Lambda (\alpha) = (E_1,...,E_q)$ be the
corresponding lattice path. Associate to $\alpha$ the sequence
 $ht(\alpha) = (0\leq ht_1(\alpha) \leq ...\leq ht_q(\alpha) =p$
by setting $ht_i(\alpha)$ to be equal the ordinate (height) of
 the end of the edge $E_i$.
\vskip .3cm

\noindent (3.9.5) Let $V_\cdot = (V_1\i V_2\i ...\i V_{q-1} \i
V_q = {\bf C}^q)$ be a complete flag of linear subspaces in
 ${\bf C}^q$ so that $dim (V_i) =i$.
Let also $\alpha \i [0,q-p]\times [0,p]$ be a Young diagram.
We define the {\it Schubert variety} $S_\alpha (V_\cdot) \i
 G(p,q)$ to be the locus of subspaces $L \i {\bf C}^q, dim (L) =p$
 such that $dim (L\cap V_i) \geq ht_i(\alpha)$ where $ht_i(\alpha)$
 was defined in (3.9.4).

It is well-known [24] that $S_\alpha (V_\cdot)$ is an irreducible
 variety
of complex dimension $|\alpha|$. The homology class in $H_{2|\alpha|}
 (G(p,q), {\bf Z})$ represented by $S_\alpha (V_\cdot)$ is independent
 on $V_\cdot$. This class is denoted by $\sigma_\alpha$ and called the
 {\it Schubert cycle}.

It is known that the homology group $H_{2r}( (G(p,q), {\bf Z})$ is freely
generated by cycles $\sigma_\alpha$ where $\alpha$ runs over all
 Young diagrams with $r$ cells contained in the retangle
$[0,q-p]\times [0,p]$.

\vskip .3cm

\noindent (3.9.6) Any Young diagram $\alpha$ (not necessarily
 contained in a given rectangle) defines the {\it Schur functor}
$\Sigma ^\alpha$ on the category of vector spaces [38]. By definition,
 for a vector space $V$ the space $\Sigma^\alpha (V)$ is the space of
 irreducible representation of the group $GL(V)$ with the highest
 weight $\alpha$. It can be defined, e.g.,  as the image of the
Young symmetrizer $h_\alpha$ in the tensor space $V^{\otimes |\alpha|}$.
 In particular, for $\alpha = (m,0....,0)$ (a horizontal strip of
length $m$) the functor $\Sigma ^\alpha$ is the symmetric power $S^m$.
 For a vertical strip $\alpha = (1^m) = (1,...,1,0,...,0)$ ( $m$ units)
 the functor $\Sigma^\alpha$ is the exterior power $\bigwedge ^m$.
\vskip .3cm

\noindent (3.9.7) For a Young diagram $\alpha$ we denote by
 $\alpha ^*$ the dual (or transposed) Young diagram defined
 by $\alpha^*_i = Card \{j: \alpha_j\geq i\}$. The rows on
 $\alpha^*$ correspond to columns of $\alpha$ and vice versa.

Now we can formulate the main result of this section.
\vskip .3cm

\proclaim (3.9.8) Theorem. Let $X\i G(k-1,n-1)$ be a special
$(k-1)$ -dimensional Veronese variety and $\Delta \in H_{2k-2}(G(k-1,n-1),
{\bf Z})$ be the homology class of $X$. Then the decomposition of
 $\Delta$ in the basis of Scubert cycles has the form
$$\Delta = \sum_{|\alpha| = k-1} m_\alpha.\sigma_\alpha, \quad\quad
 {\rm where} \quad\quad m_\alpha = dim \,(\Sigma ^{\alpha ^*}
({\bf C}^{n-k})\,\,).$$

\vskip .3cm

\noindent (3.9.9) To prove Theorem 3.9.8, we first take $X$ to
 be the variety of $(k-2)$ -dimensional chords of a Veronese curve
$C\i P^{n-2}$, which is a particular case of Veronese varieties ,
 see section (3.8). Then we degenerate $C$ into a union of lines.

More precisely, let $e_1,...,e_{n-1}\in P^{n-2}$ be the points
 corresponding to standard basis vectors of ${\bf C}^{n-1}$.
Consider the reducible curve $D= D_1\cup ...\cup D_{n-2}$ where
 $D_i$ is the line $<e_i, e_{i+1}>$. Since $D$ can be obtained
as a limit position of Veronese curves, we find that  $\Delta$
 is equal to the homology class of the chordal variety $Ch_{k-2}(D)$,
 see (3.8.1) for the notation.
\vskip .3cm

\noindent (3.9.10) The variety $Ch_{k-2}(D)$ is reducible
and splits into $n-2\choose k-1$ components $X_{i_1,...,i_{k-1}}$
 which correspond to sequences $1\leq i_1< ... < i_{k-1} \leq n-2$.
 The component $X_{i_1,...,i_{k-1}}$ is the locus of chordal subspaces
 $<x_{i_1}, ... , x_{i_{k-1}}>$ where $x_{i_\nu}$ lies on the line
 $<e_{i_\nu}, e_{i_\nu +1}> \i D$. Therefore our homology class $\Delta$
is the sum of homology classes $[X_{i_1,...,i_{k-1}}]$ of all the
 components
of $Ch_{k-2}(D)$.

\vskip .3cm

\noindent (3.9.11) Let us introduce a different, more suitable for
 our purposes, combinatorial labeling of components of $Ch_{k-2}(D)$.

Denote by $W(n-k,k-1)$ the set of all (not necessarily decreasing)
 sequences
 $\lambda = (\lambda_1,...,\lambda_{n-k})\in {\bf Z}_+^{n-k}$  of
 non-negative integers such that $\sum \lambda_i = k-1$. We shall
 call elements of $W(n-k,k-1)$   {\it weights} (in the sense of
 representation theory) in $n-k$ variables and of degree $k-1$.

To any sequence $1\leq i_1< ... < i_{k-1} \leq n-2$ we associate a weight
$\lambda (i_1,...,i_{k-1}) = (\lambda_1,...,\lambda_{n-k}) \in W(n-k,k-1)$
as follows. Let $j_1,...,j_{n-k-1}$ be all elements of the set
$\{1,...,n-2\} - \{i_1,...,i_{k-1}\}$, written in the increasing
order. Set also
$j_{n-k} = n-1$. Now define
$$\lambda(i_1,...,i_{k-1}) = (\lambda_1,...,\lambda_{n-k}),\quad
 {\rm where} \quad \lambda_\nu = j_\nu - j_{\nu -1} - 1. \leqno (3.9.12)$$
The numbers $\lambda_\nu$ are just the lengths of arithmetic
progressions with increment one into which the sequence
$(i_1,...,i_{k-1})$ splits.

The correspondence $(i_1,...,i_{k-1}) \mapsto \lambda(i_1,...,i_{k-1})$
establishes a bijection between the set of all $(k-1)$ -
element subsets in $\{1,...,n-2\}$ and the set $W(n-k, k-1)$.
 This bijection is the labeling we need.

We shall denote by $X(\lambda), \lambda\in W(n-k,k-1)$, the component
$X(i_1,...,i_{k-1}) \i Ch_{k-2}(D)$ where $\lambda(i_1,...,i_{k-1})
 = \lambda$.

Now Theorem 3.9.8 will be a consequence of the following fact.
\vskip .3cm

\proclaim (3.9.13) Theorem. Let $\alpha, \quad |\alpha| =
k-1$, be a Young diagram in the rectangle $[0,n-k-1]\times
[0,k-1]$ and let $\lambda \in W(n-k, k-1)$ be any weight.
Then the homology class of the component $X(\lambda) \i G(k-1, n-1)$
has the form
$$[X(\lambda)] \quad = \quad \sum_{|\alpha| = k-1} K_{\lambda,
\alpha^*} \cdot \sigma_{\alpha},$$
where $K_{\lambda, \alpha^*}$ is the multiplicity of weight
$\lambda$ in the irreducible representation
$\Sigma^{\alpha^*}({\bf C}^{n-k})$ (the Kostka number).

We shall concentrate on the proof of Theorem 3.9.13.

\vskip .3cm

\proclaim (3.9.14) Proposition. Let $\lambda =
 (\lambda_1,...,\lambda_{n-k}) \in W(n-k,k-1)$ be
 a weight. The component $X(\lambda)$ is isomorphic
to the product of projective spaces $\prod P^{\lambda_j}$.
 It is embedded into the Grassmannian as the image of the direct sum map
$$\oplus : \prod P^{\lambda_j} = \prod G(\lambda_j, \lambda_j +1)
 \hookrightarrow G(\sum \lambda_j, \sum (\lambda_j +1)) = G(k-1, n-1).$$

\noindent {\sl Proof:} Let $1\leq i_1<...,i_{k-1}\leq n-2$
 be sequence of integers to which $\lambda$ is associated,
 see (3.9.10). The component $X(\lambda)$ consists of chords
 which join points of lines $<e_{i_\nu}, e_{i_{\nu} +1}>$.
 Let us split the sequence $(i_1,...,i_{k-1})$ into segments
which are arithmetic progressions with increment 1. Then
$\lambda_\nu$ are precisely lengths of these segments.
 Now let $i, i+1,...,i+ \lambda_\nu$ be any such segment.
 The $(\lambda_\nu-1)$ -dimensional chords of the subcurve
 $<e_i, e_{i+1}>\cup <e_{i+1}, e_{i+2}>\cup ... \cup <e_{i+\lambda_\nu},
 e_{i+\lambda_\nu+1}>$ are just arbitrary hyperplanes in the projective
 subspace $<e_i,e_{i+1},...,e_{i+\lambda_\nu+1}>$. Any $k-2$ -
dimensional chord from our component $X_{i_1,..., i_{k-1}} =
X(\lambda)$ is therefore the projective span of hyperplanes in
 the independent projective subspaces $P^{\lambda_j} =
P({\bf C}^{\lambda_\nu +1})$, as required.
\vskip .3cm

\proclaim (3.9.15) Proposition. The coefficient at the
Schubert cycle $\sigma_\alpha$ in the decomposition of the
class $[X(\lambda)]$ equals the multiplicity
of the irreducible representation $\Sigma ^\alpha {\bf C}^{n-k}$
in the tensor product $\bigwedge ^{\lambda_1} ({\bf C}^{n-k})
\otimes ... \otimes \bigwedge ^{\lambda _{n-k}}({\bf C}^{n-k})$
of exterior powers.

The proof is based on the following (known) fact.
\vskip .3cm

\proclaim (3.9.16) Lemma. For three Young diagrams $\alpha, \beta,
 \gamma$ such that $|\gamma| = |\alpha| + |\beta|$ let $c_{\alpha
 \beta} ^\gamma$ be the multiplicity of $\Sigma ^\gamma$ in
 $\Sigma ^\alpha \otimes \Sigma ^\beta $ (the Littlewood
- Richardson number). Then:\hfill\break
a) The image of the cycle $\sigma_\alpha \otimes \sigma_\beta$
 under the direct sum map
$$ \phi: G(p_1, V_1) \times G(p_2, V_2) \rightarrow G(p_1 +p_2,
 V_1\oplus V_2)$$
equals $\sum _\gamma c_{\alpha \beta} ^\gamma \sigma_\gamma$.
\hfill\break
b) If $A$ and $B$ are finite- dimensional vector spaces then
for any Young diagram $\gamma$ we have the isomorphism of $GL(A)
\times GL(B)$ -modules
 $$\Sigma^\gamma (A\oplus B) \cong \bigoplus _{ |\alpha| + |\beta|
 = |\gamma|}
c_{\alpha \beta} ^\gamma \Sigma^\alpha (A) \otimes \Sigma ^\beta (B).$$

\noindent {\sl Proof of (3.9.16):} In part a) it suffices to treat
 the "stable" case when $V_i$ have infinite dimension. We shall
assume that it is so.

Let $H^\cdot (G(p,\infty), {\bf Z})$ be the cohomology ring of
$G(p,\infty)$ and let $Rep(GL(p))$ be the Grothendieck ring of
 polynomial representations of $GL(p)$. Let also
$\Lambda_p= {\bf Z}[x_1,...,x_p]^{S_p}$ denote the ring of
symmetric polynomials in $p$ variables $x_1,...,x_p$. There
are isomorphisms of rings
$$\Lambda_p \cong H^\cdot (G(p,\infty), {\bf Z}) \cong Rep(GL(p)),$$
which take the elementary symmetric function $e_j\in \Lambda_p$ into
the $j$ -th Chern class of the tautological bundle on $G(p, \infty)$
and into the representation $\bigwedge ^j({\bf C}^p) \in Rep(GL(p))$.

 For any Young diagram $\alpha$ denote by $s_\alpha (x_1,...,x_p)
 \in \Lambda_p$ the  {\it Schur polynomial}  (see [38]).
 It corresponds to the following elements of the two above rings:
\item{$\bullet$} The cocycle $\sigma^\alpha\in H^{2|\alpha|}
(G(p,\infty), {\bf Z}))$ dual to the Schubert cycle $\sigma_\alpha$
(i.e. $(\sigma^\alpha, \sigma_\beta) = \delta_{\alpha \beta}$).
\item{$\bullet$} The irreducible representation $\Sigma^\alpha
{\bf C}^p$ of which $s_\alpha$ is the character.

Consider the tensor product $\Lambda_{p_1}\otimes \Lambda_{p_2}$.
 It can be regarded as a ring of polynomials $f(x_1,...,x_{p_1},
y_1,...,y_{p_2})$ symmetric with respect to $x_i$ and with respect
 to $y_i$. Therefore we have an embedding
$$\delta : \Lambda_{p_1+p_2} \rightarrow \Lambda_{p_1} \otimes
\Lambda_{p_2}.$$
(it is a part of Hopf algebra structure on the limit $\Lambda =
\lim \Lambda_p$, see [38]).

The homology space of $G(p,\infty)$ is dual to $\Lambda_p$.
 The map $\phi_*$, induced on homology by the direct sum map
 $\phi$ from part a) of the lemma, is known to be dual to $\delta$.

Similarly, if we assume in part b) of the lemma that $dim(A) =
 p_1, \, dim(B) = p_2$ then the restriction map
$$Rep(GL(p_1+p_2)) \rightarrow Rep(GL(p_1)) \otimes Rep(GL(p_2))$$
is identified with $\delta$.

Thus in part a) we have to find the matrix elements of the dual map
$\delta^*:\Lambda^* \otimes \Lambda^* \rightarrow \Lambda^*$ in the
basis dual to that of $s_\alpha$. In part b) we have to find matrix
 elements of $\delta$ in the basis of $s_\alpha$. So both parts
follow from the equality
$$s_\gamma (x_1,...,x_{p_1}, y_1,...,y_{p_2})  = \sum_{ |\alpha|
 + |\beta| = |\gamma|} c_{\alpha \beta} ^\gamma
s_\alpha (x_1,...,x_{p_1})  s_\beta(y_1,...,y_{p_2}),$$
which is a reformulation of ([38], formula 5.9. Ch.I).

\vskip .3cm

\noindent (3.9.17) \underbar {\sl Proof of Proposition 3.9.15:}
  Note  that the fundamental class of the projective space
$P^{\lambda_i}$ considered as $G(\lambda_i, {\bf C}^{\lambda_i +1})$
 is the Schubert cycle corresponding to the Young diagram $(1^{\lambda_i})$
 i.e. to the vertical strip of lenght $\lambda_i$. The Schur functor
 corresponding to this diagram is the exterior power. Now the result
 follows from  Lemma 3.9.14.
\vskip .3cm

\noindent (3.9.18) \underbar {\sl Proof of Theorem 3.9.13:}
By Lemma 3.9.15, t suffices to show that the weight multiplicity
$K_{\lambda, \alpha^*}$ equals the multiplicity  of $\Sigma^\alpha$
 in $\bigwedge ^{\lambda_1} ({\bf C}^{n-k})\otimes ... \otimes
\bigwedge ^{\lambda _{n-k}}({\bf C}^{n-k})$. By Young duality this
 is equivalent to the saying that $K_{\lambda, \alpha}$, for any
 $\alpha$, is equal to the multiplicity of $\Sigma^\alpha$ in the
product of symmetric powers $S^{\lambda^1} ({\bf C}^{n-k})\otimes
 ...\otimes S^{\lambda_{n-k}}({\bf C}^{n-k})$.

To see the truth of this latter statement, decompose ${\bf C}^{n-k}$
into a sum of 1-dimensional subspaces
$L_1\oplus ...\oplus L_{n-k}$. Then decomposition of $\Sigma^\alpha
({\bf C}^{n-k})$ as a $GL(L_1)\times ... \times GL(L_{n-k})$ -module
 is just the weight decomposition. On the other hand, applying
 repeatedly Lemma 3.9.16 b) we find that
$K_{\lambda, \alpha}$ i.e. the multiplicity of $S^{\lambda_1}(L_1)
\otimes S^{\lambda _{n-k}}(L_{n-k})$ in $\Sigma ^\alpha (L_1\oplus
... \oplus L_{n-k})$ equals the multiplicity of $\Sigma ^\alpha$ in
$S^{\lambda_1} \otimes ... \otimes S^{\lambda _{n-k}}$.

Theorems 3.9.13 and 3.9.8 are completely proven.
\vskip .3cm

\noindent {\bf (3.9.19) Remark.} A different expression for the
coefficients $m_\alpha$ in Theorem 3.9.8 can be obtained from
Klyachko's formula (Proposition 1.1.8) for the homology class
of the whole Lie complex $Z\i G(k,n)$. Our Veronese variety is
 just the visible contour of $Z$ i.e., the intersection $Z\cap
 G(k-1, n-1)_p$, see (3.1.1).
The intersection map
$$H_r(G(k,n), {\bf Z} \rightarrow H_{r-n+k} (G(k-1, n-1)_p, {\bf Z})$$
is easy to describe. It takes  a Schubert cell
$\sigma_{\beta_1,...,\beta_k}$ to $\sigma_{\beta_2,...,\beta_k}$
 if $\beta_1 = n-k$ and to 0 otherwise. This leads to the formula
$$m_\alpha = \sum _{i=0}^k (-1)^i {n\choose i} dim
 (\Sigma^{n-k,\alpha_1,...,\alpha_{k-1}}({\bf C}^{k-i})).$$
According to Theorem 3.9.8 this expression equals just
$dim \Sigma ^{\alpha ^*} ({\bf C}^{n-k})$ but we do not
 know a straightforward proof of this fact.
\vskip .3cm

\noindent {\bf (3.9.20) Remark.}
 Let $\alpha = (\alpha_1,...,\alpha_{k-1}), |\alpha| = k-1$,
 be a Young diagram in the rectangle $[0,n-k]\times [0,k-1]$.
 Denote by $\bar\alpha = (n-k-\alpha_{k-1},..., n-k-\alpha_1)$
the diagram complementary to $\alpha$ in this rectangle. It is
known [24] that for any Young diagram $\beta$ with $|\beta| =
(k-1)(n-k-1)$ the intersection index $\sigma_\alpha \cdot
\sigma_\beta$ equals 1 if $\beta = \bar\alpha$ and to 0
 otherwise. Hence the coefficient $m_\alpha$ in the expansion,
by Shubert cycles, of the cycle represented by the Veronese variety
 $S\i G(k-1, n-1)$, equals $S.\sigma_{\bar\alpha}$.

Let us realize  $\sigma_{\bar{\alpha}}$ as the class of the Schubert
 variety  $S_{\bar {\alpha}} (V_. )$ for a generic flag $V_.$.
Let us take the Veronese variety $S$ to be the chordal variety
of  Veronese curve. We obtain the
 following restatement of Theorem 3.9.8:
\vskip .2cm

\item{} {\sl The number of chordal $(k-1)$ - dimensional
 subspaces of a Veronese curve in $P^{n-2}$ satisfying any
 given Schubert condition, equals the dimension of some
irreducible representation of $GL(n-k)$! The dimension
 of any representation can be  realized in this way.}

\vskip .2cm

\noindent It  would be interesting to find a conceptual
 explanation of this fact e.g., define a $GL(n-k)$ -action
on the vector space freely generated by points from
$S\cap S_{\bar\alpha} (V_.)$. Let us also point out to a
series of papers of A.N. Kirillov and N.Yu. Reshetikhin
 (see [33,34] and references therein) on new combinatorial
 formulas for weight multiplicities $K_{\lambda,\alpha}$.
 Their construction is based on an interpretation of any
 $K_{\lambda,\alpha}$ as the number of solutions of some
 special system of algebraic equations (the equations of
 Bethe - Ansatz). This interpretation seems to be connected
 with the one given above.

\vskip .3cm

\noindent {\bf (3.9.21) Example.} It is well known that
the number of nodes of a plane rational curve of degree
$d$ equals $(d-1)(d-2)/2$. We can obtain this as a
 particular case of Theorem 3.9.8. Let $C\i P^d$ be
a Veronese curve, $L\i P^d$ - a projective subspace
of dimension $d-3$ and $\pi:P^{d} -L \rightarrow P^2$ -
the projection with center $L$. Nodes of the plane curve
 $\pi(C)$ correspond to 1-dimensional chords of $C$
 intersecting $L$. Let $X = Ch_1(C) \i G(2,d+1)$ be
the surface of chordal lines of $C$. By Theorem 3.9.8,
its homology class has the form
$$[X] = dim (S^2 {\bf C}^{d-1}) \cdot \sigma_{1,1} +
 dim (\bigwedge ^2 {\bf C}^{d-1}) \cdot \sigma_{2,0}.$$
The coefficient at $\sigma_{2,0}$ equals the intersection
index of $S$ with the Schubert cycle\hfill\break
 $\sigma_{d-1,...,d-1,d-3}$. The corresponding Schubert variety
is the locus of all lines intersecting a given $(d-3)$ -dimensional
 subspace in $P^d$, for example, $L$. So we find the number of nodes
of $\pi(C)$ to be $dim (\bigwedge ^2 {\bf C}^{d-1}) = (d-1)(d-2)/2$.

\vskip .3cm

\noindent {\bf (3.9.22) Example.} The number of 4-secant lines of a
 spatial rational curve $X\i P^3$ of degree $d$ can be found by
reasoning similar to the above example. This number equals
$$dim (\Sigma^{2,2} {\bf C}^{d-3}) = (d-2) (d-3)^2 (d-4) /12.\leqno
(3.9.23)$$
The right hand side of (3.9.23) is a well-known formula for the number
of quadrisecants, see e.g. [24, Ch.2, \S 5].
\vskip .3cm

\noindent {\bf (3.9.24) Example.} The number of trisecant lines of a
rational curve of degree $d$ in $P^4$ equals
$$dim (\bigwedge ^3 {\bf C}^{n-2}) = (n-2)(n-3)(n-4)/6. \leqno (3.9.25)$$

\hfill\vfill\eject

\centerline {\bf  Chapter 4. CHOW QUOTIENT OF  $G(2,n)$ }

\centerline {\bf AND GROTHENDIECK - KNUDSEN MODULI SPACE $\overline{M_{0,n}}$.}

\vskip 1cm

In this section we study in detail the Chow quotient $G(2,n)//H$
of the Grassmannian $G(2,n)$ of lines in $P^{n-1}$.
We establish the isomorphism of this Chow quotients with the moduli space
$\overline{M_{0,n}}$ of stable $n$ -punctured curves of genus 0
introduced by
A.Grothendieck [12]
and later by F.Knudsen [37]. In particular, $G(2,n)//H$ is a smooth
 variety and the complement to the open stratum is a divisor with
normal crossing.
The relation of  the space $\overline{M_{0,n}}$ to the Grassmannian
 permits us to represent this space as an iterated blow-up of the projective
 space $P^{N-3}$.

\vskip 1cm

\centerline {\bf (4.1) The space  $G(2,n)//H$ and stable curves.}
\vskip .6cm

\noindent (4.1.1) According to Theorem 2.2.4, we have an isomorphism
$$G(2,n)//H = (P^1)^n//GL(2).$$
In other words, our Chow quotient compactifies the space
$$M_{0,n} = ((P^1)^n - \bigcup \{x_i=x_j\})/GL(2)$$
of projective equivalence classes of $n$ -tuples of distinct points on $P^1$.
 The space $M_{0,n}$ can be considered as the moduli space of systems
$(C,x_1,...,x_n)$ where $C$ is a smooth curve of genus 0 and $x_i$
are distinct points on $C$.

\vskip .3cm
\noindent (4.1.2) There is a well-known compactification of
$M_{0,n}$ by means of so-called stable $n$ - pointed curves
of genus 0 introduced by Grothendieck and  Knudsen [12][37].
Let us recall the definitions.
\vskip .3cm

\proclaim (4.1.3) Definition. A stable $n$ -pointed curve of genus 0
  is a connected (but possibly reducible) curve $C$
over $k$ together
 with $n$ smooth  distinct points $x_1,...,x_n\in C$,
 satisfying the following conditions:
\item{(1)} $C$ has only ordinary double points and every irreducible
component of $C$ is isomorphic to the projective line $P^1$.
\item{(2)} The arithmetic genus of $C$ is equal to $0$.
\item{(3)} On each component of $C$ there are at least three points
 which are either marked or double. \hfill\break
Points of $C$ which are either marked or double will be called  special.

The condition (2) is equivalent to saying that the graph formed
by components of $C$ is a tree. We shall prefer the following
"dual" point of view on this tree.
\vskip .3cm

\proclaim (4.1.4) Definition. Let $(C,x_1,...,x_n)$ be a stable
$n$ -pointed curve of genus 0. Its tree ${\cal T}(C,x_1,...,x_n)$
 has the following vertices:
\item{(1)} Endpoints (1-valent vertices) $A_1,...,A_n$ corresponding
 to $x_0,...,x_n$.
\item {(2)} Vertices corresponding to all the components of $C$.
\hfill\break
Two vertices of type (2) are joined by an edge if the corresponding
components intersect. An endpoint $A_i$ is joined by a new edge to
 the vertex of type (2)
corresponding to the unique component containing the point $x_i$.

Definition 4.1.4 is illustrated on Fig. 4.1.5.

\vbox to 4cm{\vfill\vfill (4.1.5) \vfill\vfill}

Thus edges of ${\cal T}(C,x_1,...,x_n)$ correspond to special points of $C$.
\vskip .3cm

\noindent (4.1.6) F.Knudsen has constructed in [37] the moduli
space $\overline{M_{0,n}}$ of stable $n$-pointed curves and
proved that it is a smooth compact algebraic variety.
 To formulate Knudsen's result more precisely, let us
introduce a notion of a stable $n$-pointed curve over an arbitrary
 base scheme $S$. By definition, it is a flat proper morphism
$\pi:C\rightarrow S$ together with $n$ distinguished sections
$s_1,...,s_n:S\rightarrow C$ such that for any geometric point
 $s\in S$ the fiber $C_s=\pi^{-1}(s)$ is a reduced (i.e.
 without nilpotents) algebraic curve and $(C_s, s_1(s),...,s_n(s))$
 is a stable $n$-pointed curve of genus 0. An  isomorphism between
 two such objects $(\pi:C\rightarrow S, s_1,...,s_n)$ and
 $(\pi':C'\rightarrow S, s'_1,...,s'_n)$ over the same base
$S$ is just an isomorphism $f:C\rightarrow C'$ commuting with
 projections and taking $s_i$ to $s'_i$.
\vskip .3cm

\proclaim (4.1.7) Theorem. {\rm [37]} There exists a smooth
 projective complex algebraic variety $\overline{M_{0,n}}$
 such that for any scheme $S$ over {\bf C} the set of isomorphism
 classes of stable $n$-pointed curves of genus 0 over $S$ is
 naturally identified with $Hom(S,\overline{M_{0,n}})$.

An open subset $M_{0,n}$ in $\overline{M_{0,n}}$ is formed by $n$ -
pointed curves $(C,x_1,...,x_n)$ such that $C$ is smooth i.e. $C\cong P^1$.

Now we can formulate the complete description of the Chow quotient
 of $G(2,n)$.

\proclaim (4.1.8) Theorem. The Chow quotients $G(2,n)//H$ and
$(P^1)^n//GL(2)$ are isomorphic to the moduli space $\overline{M_{0,n}}$.

To prove Theorem 4.1.8 note that Theorem 3.3.14 together with
 Corollary 3.4.9 implies the following description of
$G(2,n)//H$.
\vskip .3cm

\proclaim (4.1.9) Corollary. Take $n$ points $p_1,...,p_n$ in
$P^{n-2}$ in general position. Let $V_0(p_1,...,p_n)$ be the
 space of all Veronese curves in $P^{n-2}$ through $p_i$.
 Denote by $V(p_1,...,p_n)$ the closure of $V_0(p_1,...,p_n)$
in the Chow variety and by $W(p_1,...,p_n)$ -the closure in the
Hilbert scheme. Then $V(p_1,...,p_n)\cong W(p_1,...,p_n)\cong G(2,n)//H$.

It was proven in [29] (Theorem 0.1) that $V(p_1,...,p_n)$ and
$W(p_1,...,p_n)$
are isomorphic to $\overline {M_{0,n}}$. More precisely, any
subscheme from $W(p_1,...,p_n)$ is in fact reduced and being
 regarded together with $p_i$, is a stable $n$ -pointed curve of genus 0.
Theorem 4.1.8 is proven.
\vskip .3cm

\noindent {\bf (4.1.10) Remark.} It seems to be difficult to prove
directly that the Chow quotient $(P^1)^n//GL(2)$ coincides with
 $\overline {M_{0,n}}$. However, the Grassmannian picture (i.e.
 the Gelfand- MacPherson isomorphism, see  \S 2) leads to stable
 curves very naturally: these curves are just visible contours of
 generalized Lie complexes.
\vskip .3cm

\noindent (4.1.11) We can now give a translation to the language
of stable curves of general constructions of \S 1.

The combinatorial invariant of a stable $n$ -pointed curve
 $(C,x_1,...,x_n) \in \overline{M_{0,n}}$ is its tree ${\cal T}(C)$
 (Definition 4.2.)
 For each tree ${\cal T}$ bounding the endpoints $1,...,n$ we
define the stratum $M({\cal T})\i \overline {M_{0,n}}$ consisting
of stable $C$ curves with ${\cal T}(C)={\cal T}$. In particular,
 to a 1-vertex tree corresponds the open stratum
 $M_{0,n}\i \overline{M_{0,n}}$.

The combinatorial invariant of a generalized Lie complex
$Z\i G(k,n)//H$ is the corresponding matroid decomposition
of the hypersimplex $\Delta(k,n)$ (Proposition 1.9).
The Chow strata in $G(k,n)//H$ were defined (Definition 1.2.16)
as the loci of $Z$ for which the corresponding matroid decomposition
 is fixed.

It was proven in section 1.3 that matroid decomposition of the
 hypersimplex $\Delta (2,n)$ correspond exactly to trees bounding
 $n$ endpoints $A_1,...,A_n$, see Theorem 1.3.6. Thus we obtain the
 following corollary.
\vskip .3cm

\proclaim (4.1.12) Corollary. All matroid decompositions of the
 hypersimplex
$\Delta (2,n)$ are realizable i.e. come from non-empty Chow strata
in $G(2,n)//H$. These Chow strata have the form $M({\cal T})\i
 \overline {M_{0,n}}$. The stratum $M({\cal T})$ is isomorpic to
the product $\prod M_{0,e(v)}$, where $e$ runs over   $v$ of
${\cal T}$ and $e(v)$ is the number of edges containing $v$.

\vskip .3cm

\noindent (4.1.13) Forgetting $i$-th point on any stable
$n$-pointed curve $(C,x_1 ,..., x_n)\in \overline{M_{0,n}}$
 gives a new $n$-pointed curve. This curve might be unstable
 i.e. the condition d) of Definition 4.1.3 might be violated.
 Clearly, this happens in  the case when the component of $C$
 containing $x_i$ contains only two other double or marked points).
Blowing down this component defines a stable curve $\pi_i(C)$ pointed
with images of $x_j,j\neq i$, see [37]. It was shown in [37] that
$\pi_i$ defines a morphism
  $\overline{M_{0,n}}\rightarrow \overline{M_{0,n-1}}$ which
identifies $\overline{M_{0,n}}$ with the universal family of
 curves over $\overline{M_{0,n-1}}$.

It was shown in [29] that $\pi_i$ corresponds to geometric projection
of (degenerate) Veronese curves from $V(p_1,...,p_n)$ (Corollary 4.5)
 from the point $p_i$.
In terms of generalized Lie complexes (= points of $G(2,n)//H$) the
 projection $\pi_i$ is described as follows.
\vskip .3cm

\proclaim (4.1.14) Proposition. Let $Z\i G(2,n)$ be a generalized
 Lie complex. Let $G(2,n-1)^i$ be the space of lines in $P^{n-1}$
which in fact lie in the $(n-2)$- dimensional projective subspace
 spanned by basis vectors $e_j, j\neq i$. Then $Z\cap (G,2,n-1)^i$
is a generalized Lie complex in $G(2,n-1)_i$. The operation of
intersection with $G(2,n-1)^i$ corresponds, under the identification
of Theorem 4.1.8, to the projection $\pi_i:\overline{M_{0,n}}
\rightarrow \overline{M_{0,n-1}}$.

\vskip 1cm

\centerline {\bf (4.2) The birational maps $\sigma_i: \overline
{M_{0,n}} \rightarrow P^{n-3}$.}

\vskip 1cm

\noindent (4.2.1) The Grothendieck - Knudsen space
$\overline {M_{0,n}}$ can be seen as a "high-brow"
 compactification of the space $M_{0,n}$ of projective
equivalence classes of $n$ -tuples of distinct points
 on $P^1$, see (4.1.1). On the other hand, every three
distinct points on $P^1$ can be brought to the points
 $0,1,\infty$ by a unique projective transformation.
Doing this with the first three points of any $n$ -tuple,
 we find that
$$M_{0,n} \cong \{(x_4,...,x_n) \in {\bf C}^{n-3} : x_i \neq 0,1,
 \quad \forall i \quad {\rm and} \quad x_i\neq x_j, \quad \forall
i\neq j\}.$$
So  $M_{0,n}$ as an open subset in ${\bf C}^{n-3}$. This suggests a
"naive" compactification of $M_{0,n}$ which is just the projective
 space $P^{n-3}$ compactifying ${\bf C}^{n-3}$. One expects then
that $\overline{M_{0,n}}$, being the finer compactification, maps
 to $P^{n-3}$ by means of a regular birational map.

\vskip .3cm

\noindent (4.2.2) As shown in [29], the regular map
$\overline {M_{0,n}} \rightarrow P^{n-3}$ can be constructed
as follows. Realize $\overline {M_{0,n}}$ as the space $V(p_1,...,p_n)$
 of limit position of Veronese curves in $P^{n-2}$ containing given
 generic points $p_1,...,p_n$. For any curve $C\in V(p_1,...,p_n)$
all $p_i$ are smooth points of $C$. Fix some $i$ and consider the
 projective space $P^{n-3}_i$ of all lines in $P^{n-2}$ through $p_i$.
By associating to any curve $C\in V(p_1,...,p_n)$ its embedded tangent
 line $T_{p_i}C$ one gets a regular map
$$\sigma_i : \overline {M_{0,n}} \rightarrow P^{n-3}_i.$$
 It was demonstrated in [29] that $P^{n-3}_i$ is exactly the "naive"
compactification of $M_{0,n}$ mentioned in (4.2.1). It dependence on
$i$ is easy to explain: we need to specify which
point of an $n$ -tuple is set to be $\infty$. So the construction of
(4.2.1) corresponds to $i=3$.

Here we are going to study the maps $\sigma_i$ in more detail.
\vskip .3cm

\noindent (4.2.3) Let $L_i, i=1,...,n$, be the line bundle on
 $\overline{M_{0,n}}$ whose fiber at a pointed curve
 $(C,x_1, ..., x_n)$ is $T^*_{x_i}C$, the cotangent space to $C$ at $x_i$.
 Clearly $L_i \cong \sigma_i^* ({\cal O}_{P^{n-3}_i}(1))$
 The following fact was proven in [29].
\vskip .3cm

\proclaim (4.2.4) Proposition. For any $i\in \{1,...,n\}$
 the space $H^0(\overline{M_{0,n}}, L_i)$ has dimension $n-2$.
 The corresponding morphism $\gamma_{L_i}$ is everywhere regular,
birational and, moreover, one-to one outside the subvariety in
 $P^{n-3}_i$ formed by lines which lie an a hyperplane spanned
by $p_i$. In the Veronese picture the space
$P(H^0(\overline{M_{0,n}},L_i)^*)$ is identified with
$P^{n-3}_i$ and $\gamma_{L_i}$ is identified with $\sigma_i$.
\vskip .3cm

\noindent (4.2.5) Let us give a description of maps $\sigma_i$
on the language of generalized Lie complexes. So we start with
 the standard coordinatized projective space $P^{n-1}= P({\bf C}^n)$.
   For any point $x\in P^{n-1}$  we shall denote by $P^{n-2}(x)\i G(2,n)$
 the space of lines in $P^{n-1}$ meeting $x$. Let $e_i\in P^{n-1}$
 be the $i$-th basis vector.
\vskip .3cm

\proclaim (4.2.6) Proposition. Any Lie complex
 (and hence any generalized Lie complex) in $G(2,n)$ contains
each projective space $P^{n-2}(e_i)$.

\noindent {\sl Proof:} Let $Z$ be a Lie complex. The intersection
$Z\cap P^{n-2}_{e_i}$  contains the closure of a generic torus orbit
in $P^{n-2}(e_i)$. Since this generic orbit is dense in $P^{n-2}(e_i)$,
the assertion follows.
\vskip .3cm

\noindent (4.2.7) Note that the dimension of a (generalized) Lie complex
 $Z$ is just by one greater then that of $P^{n-2}(e_i)$.
 Hence at a generic point $l$ of $P^{n-2}(e_i)$ the tangent
 space $T_lZ$ represents a line in the normal space
$N_l(Z/G(2,n))= T_lG(2,n)/T_lZ$. In  the construction
 of the Veronese curve corresponding to $Z$ (as
 the visible contour, see  \S (3.1)) we considered
 the set of all lines in $Z$ meeting a given point
$u$ or, in other words, the intersection $Z\cap P^{n-2}(u)$.
 More precisely, we specialized to $u=e=(1,...,1)$.

\proclaim (4.2.8) Proposition. The space $P^{n-3}_i$ is
 naturally identified with the projectivization of the
 normal space to $P^{n-2}(e_i)$ in $G(2,n)$ at the point
 $p_i=<e_i,u>$. The map $\sigma_i$ is identified to the map
taking any generalized Lie complex $Z$ to the line in the
 above normal space given by the subspace $T_lZ$.

\noindent {\sl Proof:} The subvarieties $P^{n-2}(e_i)$
and $P^{n-2}(u)$ in $G(2,n)$ are of middle dimension and
intersect transversely in the point $p_i=<e_i,u>$. Therefore
the normal space in question is naturally identified with the
 tangent space at $p_i$ to $P^{n-2}(u)$. Since $\sigma_i$ is
 defined by considering the tangent line to $Z\cap P^{n-2}(u)$
 at $p_i$, the assertion follows.
\vskip .3cm

\noindent (4.2.9) The advantage of the description  in (4.1.20)
is that it clearly states the dependence on the choice of a point
 $u$. It also shows how to obtain a more invariant description
 of $\sigma_i$. To do this, one should  by consider all the
 1-dimensional subspaces in the normal spaces to $P^{n-2}(e_i)$
 in $G(2,n)$ at all the generic points or, in other words, the
 corresponding subbundle in the normal bundle. Let us describe
these bundles.
\vskip .3cm

\proclaim (4.2.10) Proposition. Let $x\in P^{n-1}$ be any point.
 Then the normal bundle of the subvariety $P^{n-2}(x)\i G(2,n)$
is naturally isomorphic to the twisted tangent bundle
 $T_{P^{n-2}(x)}\otimes {\cal O}_{P^{n-2}(x)}(-1)$.

\noindent {\sl Proof:} Let $l\in P^{n-2}_x$ be any line in
$P^{n-1}$ containing $x$ and let $ N =  N_{l/P^{n-1}}$ be
 the normal bundle of $l$. The tangent space $T_lG(2,n)$
is identified (by Kodaira-Spencer) with the space $H^0(L, N)$
i.e. with the space of normal vector fields on $l$. The subspace
$T_lP^{n-2}_x \i T_l G(2,n)$ consists of those fields $v$ which
vanish at $x$ i.e. $v(x)=0$. Hence we have a linear map
$$(N_{P^{n-2}(x)/G(2,n)})_l \mapsto (N_{l/P^{n-1}})_x,\quad v\mapsto v(x).$$
It is immediate to check that this map is in fact an isomorphism.
Let now $E$ be the vector bundle on $P^{n-2}(x)$ whose fiber over
 a line $l$ is the normal space $(N_{l/P^{n-1}})_x$. We have proven that
 $N_{P^{n-2}(x)/G(2,n)}$  is isomorphic to $E$. Let us regard
$P^{n-2}_x$ as the projectivisation of the vector space
$W = T_x P^{n-1}$. Then we have the following description
of the bundle $E$ on $P(W)$: the fiber of $E$ at a 1-dimensional
 linear subspace $\Lambda \i W$ is $W/\Lambda$. This is the
 standard description of $T_{P(W)} \otimes {\cal O}_{P(W)}(-1)$,
 the so-called Euler sequence, see [45].
\vskip .3cm

\noindent (4.2.11) Proposition 4.2.10 implies that the projectivization
 of the normal bundle \hfill\break
$N_{P^{n-2}(x)/G(2,n)}$ is the same as that of the tangent bundle
$T_{P^{n-2}(x)}$. So  any (generalized) Lie complex defines a 1-
 dimensional subbundle in the tangent bundle of $P^{n-2}(e_i)$
(this subbundle can be defined only over generic points;
over some special points of $P^{n-2}(e_i)$ it may have
singularities i.e. become a non- locally free coherent sheaf).
In other words, we have a field of directions in $P^{n-2}(e_i)$.
Let us denote this field by $\Sigma_i(Z)$ in $P^{n-2}(e_i)$.
Let $H_i$ be the coordinate hyperplane in $P^{n-1}$ opposite to $e_i$.
The projection from $e_i$ identifies $H_i$ with $P^{n-2}(e_i)$
 so we can consider the direction field $\Sigma_i(Z)$
as being defined on $H_i$. It can be regarded as the choice
- free materialization of $\sigma_i (C)$ where $C$ is the stable
 curve corresponding to $Z$.
\vskip .3cm

\noindent (4.2.12) Let us describe the  direction field $\Sigma_i(Z)$
 geometrically. Note that since $Z$ is $H$- invariant the visible
 contour $Z_p$ does not change, up to isomorphism, for points $p$
lying in one torus orbit. We shall  look how does $Z_p$ split when
$p$ goes to a point on $H_i$ not lying on other $H_j$.

\vskip .3cm

\proclaim (4.2.13) Proposition. The field of directions $\Sigma_i(Z)$
is well defined at any point of $H_i$ not lying on coordinate hyperplanes.
 For such a point $x$ the visible cone  $Z_x = Z\cap P^{n-2}(x)$
splits into a stable curve (family of lines) in $P^{n-2}(x)$, all
 whose lines lie in $H_i$ and a plane  pencil of lines containing
$<x,e_i>$. The direction $\Sigma_i(Z)(x)$ is given by the unique
line of this pencil lying in $H_i$.

\noindent {\sl Proof:} Let $C$ be the visible contour of $Z$ at
 the point $e=(1,...,1)$. Let $K$ be the union of all lines from $C$.
 This is a cone with vertex $e$ containing all the lines $<e,e_i>$.
 All points of $<e,e_i>$ except $e$ are smooth points of $K$.
Let $L$ be the embedded tangent (2-) plane to $K$ along
$<e,e_i>$. The $L$ lies in $P^{n-2}(e_i)$ and is, by definition,
 equal to $\sigma_i(C)$. To prove the proposition, it suffices to
 consider the case when
$x\in H_i$ is the barycenter of the coordinate simplex i.e. all
homogeneous coordinates of $x$ are equal to 1 except the $i$ -th
 which is equal to 0. Consider the transformation $\gamma (t) =
 (1,...,1,t,1,...,1) \in ({\bf C}^*)^n$ ($t$ is on
the $i$ -th place). Then we have $x= \lim_{t\rightarrow 0} \gamma(t)
\cdot e$.
Hence the visible cone at $x$ of the (torus - invariant) complex
 $Z$ equals
$\lim_{t\rightarrow 0} \gamma (t)\cdot K$. But the latter limit
 will be the union of the 2-plane $L$ (which is preserved under
all $\gamma(t)$ and some other part which will lie inside $H_i$.
It remains to show that the intersection $L\cap H_i$ is precisely
the line in $H_i$ whose direction is the value at $x$ of the
 direction field $\Sigma _i (Z)$. This checking is left to the reader.

\vskip 1cm

\centerline {\bf (4.3)
  Representation of the space $\overline{M_{0,n}} = G(2,n)//H$
as a blow-up.}

\vskip 1cm

In the previous sections we have constructed regular birational
morphisms $\sigma_i$ from $\overline{M_{0,n}}= G(2,n)//H$ to
 projective spaces. When such a morphism is found, it is always
 desirable to decompose it to simpler ones. Standard examples of
"simplest" regular birational morphisms are provided by blow- ups.
\vskip .3cm

\noindent (4.3.1)  Recall [24,25] that the {\it blow-up}
 (or monoidal transformation, or sigma- process) $Bl_YX$ is defined
for any smooth closed subvariety $Y$ (which is called the
{\it center}) in a smooth variety $X$ . This is a new smooth
 variety equpped with a canonical morphism $p$ to $X$.
The morphism $p$ is one-to one outside $Y$ and for any $y\in Y$
the preimage $p^{-1}(y)$ is canonically identified with the
 projectivization of the normal space $T_yX/T_yY$. If $Z\i X$
is another submanifold not contained in $Y$ then the
 {\it strict preimage } (or proper transform) of $Z$
is the closure $\tilde Z$ of $p^{-1}(Z-Y)$ in $Bl_YX$.
The subvariety $\tilde Y\i Bl_YX$ can be, in its turn blown up,
thus giving an iterated blow-up which is abusively denoted
 $Bl_ZBl_YX$. This blow-up does not, in general, coincide with
$Bl_YBl_Z(X)$ (it does, if $Y$ and $Z$ are disjoint).
Similar construction can be performed for several subvarieties
$Y_1,...,Y_r$.

\vskip .3cm

\noindent (4.3.2) Our aim in this section is to decompose the
 morphism $\sigma_i$ into a sequence of monoidal transformations
and therefore to give a "constructive" definition of
 $\overline{M_{0,n}}= G(2,n)//H$ as a iterated blow- up
 of a projective space. Proposition 4.10 suggests that in
 order to obtain $\overline {M_{0,n}}$ we should blow up
$n-1$ generic points $q_1,...,q_{n-1}$ in $P^{n-3}$, and all
 projective  subspaces spanned by them. However, an iterated
 blow- up depends on the ordering of the centers, so the question is delicate.
\vskip .3cm

\proclaim (4.3.3) Theorem. Choose $n-1$ generic points
$q_1,...,q_{n-1}$ in $P^{n-3}$.The variety
$\overline{M_{0,n}}$ can be obtained from $P^{n-3}$ by a
  series of blow ups of all the projective spaces spanned
by $q_i$. The order of these blow-ups can be taken as follows : \hfill\break
1) Points $q_1,...,q_{n-2}$ and all the projective subspaces
spanned by them in order of the increasing dimension;\hfill\break
2) The point $q_{n-1}$ , all the lines $<q_1, q_{n-1}>,...,
<q_{n-3}, q_{n-1}>$ and subspaces spanned by them in order of
 the increasing dimension;\hfill\break
3)The line $<q_{n-2},q_{n-1}>$, the planes $<q_i, q_{n-2}, q_{n-1}>,
i\neq n-3$ and all subspaces spanned by them in order of the
 increasing dimension. \hfill\break
etc. etc

\vskip .3cm

\noindent {\bf (4.3.4) Remark.} A representation of
 $\overline{M_{0,n}}$ as an iterated blow-up of the
Cartesian power $(P^1)^{n-3}$ was given by S.Keel.
Still another representation of $\overline{M_{0,n}}$
 as a blow-up of $P^{n-3}$, different from the one
given here, can be deduced from a more general construction
 of W.Fulton and R.MacPherson [15]. In Fulton-MacPherson
construction all the centers of blow-ups have codimension 2.

\vskip .3cm

\noindent (4.3.5) The rest of this section will be devoted
to the proof of Theorem 4.3.3.

In  Proposition 1.4.12 we have constructed some regular
birational morphisms $f_I$  of general Chow quotient of
Grassmannian $G(k,n)$ to the "secondary variety" of the
 product of two simplices $\Delta^{k-1}\times\Delta^{n-k-1}$.
 We want now to analyze these morphisms for the Chow quotient
 of $G(2,n)$ which is $\overline{M_{0,n}}$ in order to use them
as a halfway approximation to a required sequence of blow-ups.
\vskip .3cm

\noindent (4.3.6) Recall that for every two-element subset
 $I=\{i,j\}\i \{1,...,n\}$ the coordinate subspace
 ${\bf C}^I ={\bf C}e_i \oplus {\bf C}e_j \i {\bf C}^n$
is a fixed point for the torus action on $G(2,n)$ and hence
 our torus $H$ acts in the tangent space $T_I = T_{{\bf C}^I}G(2,n)$
and on its projectivization. To each $H$- orbit in $G(2,n)$ whose
 closure contains ${\bf C}^I$ is therefore associated an $H$
- orbit in the projective space $P(T_I)$. The map $f_I$ from
$G(2,n)//H = \overline {M_{0,n}}$ to $P(T_I)//H$ is induced
by this correspondence.
\vskip .3cm

\noindent (4.3.7) As we explained in  \S 0.2, the Chow quotient
 of a projective space by a torus $H$ is a toric variety
corresponding to the secondary polytope of the point
 configuration given by the characters of $H$ defining
 the action. In our case the space $T_I$ is identified
 with the space of 2 by $n-2$ matrices $||a_{ij}||, i\in I,
j\in \bar I$ and the action of a torus element $(t_1,...,t_n)$
 on such a matrix gives a new matrix $||t_i^{-1}t_ja_{ij}||,
 i\in I, j\in\bar I$. The $H$- characters are therefore
identified with vectors $e_j-e_i, i\in I, j\in\bar I$ of
${\bf Z}^n$. These vectors are vertices of the simplicial
 prism $\Delta^1\times \Delta^{n-3}$ which we shall also
denote $\Delta^I\times \Delta^{\bar I}$ to emphasize the dependence on $I$.

For the case of simplicial prisms $\Delta^1\times \Delta^k$
 triangulations have a complete and simple description.
\vskip .3cm

\noindent (4.3.8)
Note that the symmetric group $S_{k+1}$ acts on
$\Delta^1\times \Delta^k$ by permuting the vertices
 of the second factor and hence acts on the triangulations
 of $\Delta^1\times \Delta^k$.

Let us describe the {\it standard} triangulation of
$\Delta^1\times \Delta^k$ used in combinatorial topology
 [16]. It depends on the numberings of the vertices of factors.
 To fix these numberings  denote the vertices of our prism by
 pairs $(a,b)$ where $0\leq a\leq 1, 0\leq b\leq k$. The
triangulation $T_{st}$ consists of the simplices $\Delta_i,
0\leq i\leq k$, where $\Delta_i$ is the convex hull of
 $(0,j), j\leq k$ and $(1,j),j\geq k$. The characteristic
 function of this triangulation (i.e. the corresponding
vertex of the secondary polytope, see \S (0.2)) equals $\phi_{st}(i,j) = j+1$.

\vskip .3cm

\proclaim (4.3.9) Proposition.  There exist exactly $(k+1)!$
 regular triangulations of the prism $\Delta^1\times \Delta^k$.
 All they can be obtained from the standard one by action of $S_{k+1}$.

In fact, all the triangulation of the prism are regular,
 but we do not need this.

\noindent {\sl Proof:} Let $\Sigma$ be the secondary polytope
 of $\Delta^1\times \Delta^k$. Its vertices are functions
$\phi_T(i,j), i=0,1; j=0,...,k$ where $T$ runs over all the
 regular triangulations. Let us use the original interpretation
of the secondary polytope as the Newton polytope of the principal
 determinant [22,23]. In our situation this means the following.

 Consider a $2\times (k+1)$ -matrix
$$A = \pmatrix { a_{00} &a_{01}& ...&a_{0k} \cr
a_{10}&a_{11}&...&a_{1k}}$$
with indeterminate entries. Consider the polynomial $E(A) =
 (\prod_{p,j} a_{pj})\cdot \prod_{0\leq i < j\leq k}D_{ij}(A)$
where $D_{ij}(A) = a_{0i}a_{1j} - a_{0j}a_{1i}$ is the minor o
f $A$ on $i$ -th and $j$ - th column. Then, as shown in [23],
$\Sigma$ is the Newton polytope of $E$ i.e. the convex hull in
 $Mat(2\times (k+1), {\bf R})$ of integral points $\omega =
 ||\omega_{pj}||\in Mat(2\times (k+1), {\bf Z}_+)$ such that
the monomial $\prod_{p,j} a_{pj}^{\omega_{pj}}$ enters $E(A)$
 with non-zero coefficient.
On the other hand, $E(A)$ can be found explicitly by means of
the Vandermonde
determinant:
$$E(A) = (\prod_{p,j} a_{pj}) \cdot det \pmatrix {
a_{00}^k&a_{01}^k&...&a_{0k}^k\cr
a_{00}^{k-1}a_{10}&a_{01}^{k-1}a_{11}&...&a_{0k}^{k-1}a_{1k}\cr
a_{00}^{k-2}a_{10}^2&a_{01}^{k-2}a_{11}^2&...&a_{0k}^{k-2}a_{1k}^2\cr
\vdots&\vdots&\vdots&\vdots\cr
a_{10}^k&a_{11}^k&...&a_{1k}^k }$$
The exponent vector of any monomial of this polynomial is
 obtained from the vector $\phi_{st}$ described in (4.3.8),
 by a permutation of columns. Proposition is proven.

\vskip .3cm

\proclaim (4.3.10) Corollary. The secondary polytope of
$\Delta^1\times \Delta^k$ is linearly isomorphic to the
convex hull of the $S_{k+1}$ -orbit of the point $(1,2,...,k+1)
\in {\bf Z}^{k+1}$.

This polytope is known as the $k$- dimensional {\it permutohedron}
and denoted ${\cal P}_k$. It is a particular case of so-called
general hypersimplices associated to homogeneous spaces $G/P$ by
 I.M.Gelfand and V.V.Serganova [21].
\vskip .3cm

\noindent (4.3.11) The toric variety corresponding
to ${\cal P}_k$ will be called the $k$- dimensional
 {\it permutohedral space} and denoted $\Pi^k$. Here
 is one of the description of this space which generalizes
 to arbitrary $G/P$- hypersimplices.
\vskip .3cm

\proclaim (4.3.12) Proposition. {\rm [21]}$\quad$ The
permutohedral space $\Pi^k$ is isomorphic to the closure
 of a generic orbit of torus $({\bf C}^*)^{k+1}$ on the
space of complete flags of linear subspaces of ${\bf C}^{k+1}$.

We will be interested in a slightly different point of view on
 $\Pi^k$ realizing it as an explicit blow-up of a projective space $P^k$.
\vskip .3cm

\proclaim (4.3.13) Proposition. The permutohedral space $\Pi^k$
can be obtain from the projective space $P^k$ by the following
 sequence of blow-ups. First blow up
  $k+1$ generic points (the projectivizations of basis vectors)
then blow up the strict preimages of all coordinate lines joining
them, then the strict preimages of coordinate planes etc.

\noindent {\sl Proof:} Let $F_i$ be the space of $(1,2,...,i)$-
 flags in ${\bf C}^{k+1}$. Let $X=X_k$ be   the closure of a
 generic orbit of $({\bf C}^*)^k$ in $F_k$ and $X_i$- the projection
of $X$ to $F_i$. Then $X_1$ is the projective space $F_1= P^k$.
 It is straightforward to see that each projection $X_i\rightarrow
 X_{i-1}$ realizes $X_i$ as the blow-up of strict preimages of all
 $(i-1)$- dimensional projective subspace spanned by basis
 vectors of ${\bf C}^{k+1}$.
\vskip .3cm

\noindent {\bf (4.3.14) Remarks.} a)Note that the
orbit closure $\Pi^k = X \i F_k$ can be mapped as well to the
projective space of hyperplanes in ${\bf C}^{k+1}$.
 Considering the decomposition of this projection
through spaces of $(i,i+1,...,k)$- flags we find that
$X$ is represented as the blow- up of the dual projective
 space $P^{k\vee}$ similar to that of Proposition 4.3.13.
 The corresponding birational map from $P^k$ to  $P^{k\vee}$
is the standard Cremona inversion [24]. Thus the permutohedral
 space provides an explicit decomposition of the Cremona
inversion to sigma- processes and their inverses.

b) In the correspondence between convex polytopes and toric
varieties blowing up the closure of an orbit corresponds to
chiseling of the face corresponding to this orbit, see [49].
 Proposition 4.3.13 amounts to the following construction of
permutohedron from the simplex. First cut out all vertices,
then all edges etc.

\vskip .3cm
\noindent (4.3.15) Let us relate the regular birational
morphisms $\sigma_i :\overline{M_{0,n}}\rightarrow P^{n-3}_i$
 and $f_{ij}:\overline{M_{0,n}}\rightarrow \Pi^{n-3}_{ij}$.
\vskip .3cm

\proclaim (4.3.16) Proposition. There exist regular birational
morphisms $\tau_{ij}: \Pi^{n-3}_{ij}\rightarrow P^{n-3}_i$ such
 that the composite morphisms $\overline{M_{0,n}}\buildrel f_{ij}
\over \rightarrow \Pi^{n-3}_{ij}\buildrel \tau_{ij}\over\rightarrow
P^{n-3}_i$ coincide with $\sigma_i$.

\noindent {\sl Proof:} The choice of coordinates identifies the
tangent space to $G(2,n)$ at the fixed point $<e_i, e_j>$ with
the open Schubert cell in $G(2,n)$ consisting of all lines not
intersecting the span of points $e_m, m\neq i,j$. For a generalized
 Lie complex $Z$ the point $\sigma_i(Z)$ can be, by considerations
of \S4, be read off from the normal spaces to $Z$ at a generic point
of $P^{n-2}(e_i)$. Such a generic point can be contracted, by the
 action of the torus $H$, to the point $<e_i, e_j>$. Therefore,
our normal space in question can be recovered from the part of $Z$
 which can be contracted to this fixed point i.e. from $f_{ij}(Z)$.
\vskip .3cm

\proclaim (4.3.17) Proposition. The space $\overline {M_{0,n}}$
coincides with the closure of the open stratum $M_{0,n}$ in the
inverse limit of $\Pi^{n-3}_{ij}$ and $P^{n-3}_i$.

\noindent {\sl Proof:} First let us show that the natural map of
$\overline {M_{0,n}}$ into the said inverse limit is injective.
This means that if two generalized Lie complexes $Z, Z'$ induce
the same algebraic cycles in the projectivizations of all the
 tangent spaces $T_{<e_i,e_j>}G(2,n)$ then they coincide. This
is obvious since any $H$- orbit has in its closure some fixed
 point.

Hence we have a regular morphism (denote it $\psi$) of
$\overline {M_{0,n}}$ to the inverse limit in question
which is bijective on {\bf C} -points. To show that $\psi$
is in fact an isomorphism  of algebraic varieties, it suffices
 to show that the differential of $\psi$ does not annihilate non
-zero tangent vectors to
$\overline {M_{0,n}}$. This is done similarly to the proof of
Theorem 3.3.14.
\vskip .3cm

\proclaim (4.3.18) Proposition. The map
$$f_{ij}\times \pi_i :\overline {M_{0,n}}\rightarrow \Pi^{n-3}_{ij}
\times \overline{M_{0,n-1}}$$
is an embedding of algebraic varieties.

\noindent {\sl Proof:}  We shall check only the injectivity on {\bf C}
 -points leaving the injectivity on tangent vectors to the reader.
 Let ${\cal T}$ be a tree  bounding the endpoints $1,...,n$ and
 $M({\cal T})$ be the corresponding stratum of $\overline{M_{0,n}}$
 or, what is the same, the corresponding Chow stratum in $G(2,n)//H$.
 Let $Z$ be a generalized Lie complex  from $M({\cal T})$ and $C=Z_u$
 be the corresponding stable $n$-pointed curve of genus 0.

 The value of $f_{ij}(Z)$ depends only of the components of $Z$
 containing the fixed point $<e_i,e_j>$. The components correspond
 to vertices of ${\cal T}$ lying on the  path $[ij]$  -the shortest
edge  path joining $i$-th and $j$-th endpoints.
Denote these vertices, in natural order of movement from $i$
 to $j$, by   $v_1,...,v_r$. Let $s_\nu$ be the number of edges
meeting $v_\nu$. To the chain of vertices $v_\nu$ corresponds a
chain of irreducible components $C_1,...,C_r$ of $C$ and on each
$C_\nu$ we have $s_\nu$ marked points. The projective configurations
 of these groups of points are precisely what is taken into account
 by the map $f_{ij}$ on the curve from $M({\cal T})$. Now our
assertion means that the isomorphism class of the stable $n$-
pointed curve $C$ can be recovered from two groups of data:

 a) The isomorphism class of the stable  $(n-1)$-pointed curve $\pi_i(C)$ ;

b) The isomorphism class of the stable curve $C' = C_1\cup...
\cup C_r\i C$ pointed by $x_i, x_j$ and all the marked and double
 points of $C$ lying on $C'$.

This is obvious and Proposition 4.3.18 is proven.

\vskip .3cm

\noindent (4.3.19) Let us now connect the spaces
 $\overline {M_{0,n-1}}$ and $\Pi^{n-3}_{ij}$. We view the
latter space at the blow-up of $P^{n-3}_i$ at all vertices,
 edges etc of the coordinate simplex formed by points $q_m,
 m\neq j$. Projecting these points from $q_j$ gives a circuit
in the space $P^{n-4}_{ij}$ of lines in $P^{n-3}_i$ meeting $q_j$.
The space $P^{n-4}_{ij}$ is in the same relation to $\overline
 {M_{0,n-1}}= \pi_i(\overline{M_{0,n}})$ as $P^{n-3}_i$ was to
 $\overline{M_{0,n}}$. In particular, we have the regular
 birational morphism $\sigma_{j/i}$ from $\overline {M_{0,n-1}}$
to $P^{n-4}_{ij}$. On the other hand, consider the blow-up
 $Bl_{p_i}\Pi^{n-3}_{ij}$. It also posesses a projection to
 $P^{n-4}_{ij}$. Proposition 4.3.18 implies the following
 corollary:

\proclaim (4.3.20) Corollary. The space $\overline{M_{0,n}}$
coincides with the closure of $M_{0,n}$ in the fiber product
of $\overline{M_{0,n-1}}$ and $Bl_{p_j}\Pi^{n-3}_{ij}$ over $P^{n-4}_{ij}$.

This corollary can be reformulated as follows. Suppose we
knew the way of constructing $\overline {M_{0,n-1}}$ as an
 iterated blow-up of the projective space $P^{n-4}_{ij}$
whose centers are proper transforms of smooth subvarieties
 $Y_1,...,Y_r$. Then we have in $Bl_{p_j}P^{n-3}_j$ the
varieties $\tilde Y_\nu$ which are blow-ups of cones over
 $Y_\nu$ with apex $p_j$. The corollary means that if we
perform the sequence of blow-ups of $Bl_{p_j}\Pi^{n-3}_{ij}$
with centers in proper transforms of $\tilde Y_\nu$ then we
 obtain $\overline {M_{0,n}}$

In other words, the problem of recovering $\overline{M_{0,n}}$
 from the partial blow-up $Bl_{p_i}\Pi^{n-3}_{ij}$ is equivalen
t to the problem of recovering $\overline{M_{0,n-1}}$ from the
 projective space. This gives an inductive proof of Theorem 4.3.3.

\hfill\vfill\eject

\beginsection References

\vskip 1cm

\item{1.} M.F.Atyah, The geometry and physics of knots,
Cambridge Univ. Press, 1991.

\item{2.} H.Baker, Principles of Geometry, vol. 3-4,
Cambridge Univ. Press, 1925 (republished by F.Ungar Publ. New York, 1963).

\item{3.}D.Barlet, Espace analytique reduit des cycles
 analytiques compexes compacts, in: Lecture Notes in
Mathematics, {\bf 482}, p.1-158, Springer-Verlag, 1975.

\item{4.} D.Barlet, Note on the "Joint theorem", preprint 1991.

\item{5.} D.Bayer, I.Morrison, Gr\"obner bases and geometric
 invariant theory I, {\it J.Symbolic Computation}, {\bf 6} (1988), 209-217.

\item{6.} D.Bayer, M.Stillman, A theorem on refining division
 orders by the reverse lexicographic order, {\it Duke Math. J.},
{\bf 55} (1987), 321-328.

\item{7.} A.A.Beilinson, R.D.MacPherson, V.V.Schechtman, Notes
on motivic cohomology, \hfill\break
{\it Duke Math.J.}, {\bf 54} (1977), 679-710.

\item{8.} L.J.Billera, B.Sturmfels, Fiber polytopes, to appear.

\item{9.} L.M.Brown, The map of $S_m$ by means of its $n$ -
ic primals, {\it J. Lond. Math. Soc.}, {\bf 5} (1930), 168-176.

\item{10.} A.Byalynicki-Birula, A.J.Sommese. A conjecture
 about compact quotients by tori, {\it Adv. Studies in Pure Math.},
{\bf 8} (1986), 59-68.

\item{11.} A.B.Coble, Algebraic geometry and theta-functions,
(AMS Colloquium Publ., vol. 10) Providence RI 1928 (Third Ed., 1969).

\item{12.} P.Deligne, Resum\'e des premi\`ers expos\'es de A.
Grothendieck, in: SGA7, Exp. I, Lecture Notes in Math. {\bf 288},
p.1-24, Springer-Verlag, 1972.

\item{13.} P.Deligne, Theorie de Hodge II, {\it Publ. Math. IHES},
{\bf 40} (1971), 5-58.

\item{14.} I.Dolgachev and D.Ortland, Point sets in projective
spaces and theta functions,\hfill\break {\it Ast\'erisque},
 {\bf 165}, Soc. Math. France, 1988.

\item {15.} W.Fulton, R.D.MacPherson, The compactification
 of configuration spaces, Preprint, 1991.

\item{16.}P.Gabriel, M.Zisman. Calculus of fractions and homotopy theory,
Springer- Verlag, New York, 1967.

\item{17.} A.M.Gabrielov, I.M.Gelfand and M.V.Losik,
Combinatorial calculation of characteristic classes I,
 {\it Preprint 1975-12}, Institute of Applied Math.AN USSR,
 Moscow, 1975 (Reprinted in: Collected papers of I.M.Gelfand, vol.3,
 Springer, 1989).

\item{18.} I.M.Gelfand, General Theory of hypergeometric functions,
 {\it Sov. Math. Dokl.}, {\bf 33} (1986), 573-577 (Collected Papers,
 Vol.3, p.877-881).

\item{19.} I.M.Gelfand, R.M.Goresky, R.W.MacPherson, V.V.Serganova,
Combinatorial Geometries, Convex polyhedra and Schubert cells,
 {\it Adv. in Math}, {\bf 63} (1978), 301-316.

\item{20.}I.M.Gelfand and R.W. MacPherson, Geometry in
 Grassmannians and a generalization of the dilogarithm,
{\it Adv. in Math.}, {\bf 44} (1982), 279 -312.

\item{21.}I.M.Gelfand and V.V.Serganova, Combinatorial
geometries and torus strata on compact homogeneous spaces,
 {\it Uspekhi Mat. Nauk}, {\bf 42} (1987), 107-134.
(Reprinted in: Collected Papers of I.M.Gelfand, vol.3, p.926-958.)

\item{22.}
 I.M.Gelfand, A.V.Zelevinsky, M.M.Kapranov, Discriminants
 of polynomials in several variables and triangulations of
 Newton polytopes, {\it Leningrad Math. Journal}, {\bf 2} (1990), No.3, 1-62.

\item{23.} I.M.Gelfand, A.V.Zelevinsky, M.M.Kapranov, Newton
polytopes of principal A- determinants, {\it Sov. Math. Dokl.},
 {\bf 40} (1990), 278-281.

\item {24.} Ph. Griffiths, J.Harris, Principles of Algebraic
 Geometry, Wiley- Interscience, New York, 1978.

\item{25.} R.Hartshorne, Algebraic geometry (Graduate texts
in Math., {\bf 52}), Springer- Verlag, 1977.

\item{26.} Y.Hu, Geometry and topopogy of quotient varieties,
Thesis, MIT, 1991.

\item{27.} C.M.Jessop, A Treatise on the Line Complex, Cambridge
 Univ. Press, 1903 (republished by Chelsea Publ. Co.,New York, 1969).

\item{28.} F.Junker, Ueber symmetrische Funktionen von mehreren
 Reihen von Ver\"anderlichen, {\it Math. Ann.}, {\bf 43} (1893), 225- 270.

\item{29.} M.M.Kapranov, Veronese curves and Knudsen's moduli
 space $\overline{M_{0,n}}$, preprint, Northwestern University, 1991.

\item{30.} M.M.Kapranov, B.Sturmfels, A.V.Zelevinsky,
 Chow polytopes and general resultants, {\it Duke Math. J.}, to appear.

\item{31.} M.M.~Kapranov, B.~Sturmfels and A.V.~Zelevinsky,
Quotients of toric varieties,
{\it Math. Ann.} {\bf 290} (1991), 643-655.

\item{32.} S.Keel, Intersection theory of Moduli spaces
of stable n-pointed curves of genus zero, to appear in Trans. AMS.

\item{33.} A.N.Kirillov, N.Yu. Reshetikhin, Multiplicity
of weights in irreducible representations of a general
linear superalgebra, {\it Funct. Anal. Appl.} {\bf 22} (1988), 328-330.

\item{34.} A.N.Kirillov, N.Yu. Reshetikhin, The Bethe Ansatz
and the combinatorics of Young tableaux, {\it J. Soviet. Math.} {\bf 41}
(1988), 925-955.

\item{35.} F.C.Kirwan, Partial desingularizations of quotients
 of nonsingular varieties and their Betti numbers, {\it Ann. Math.},
 {\bf 122}, (1985), 41-85.

\item{36.} A.A.Klyachko, Orbits of the maximal torus on the flag
space, {\it Funct. Anal. Appl.}, {\bf 19}, No.1 (1985), 77-78.

\item{37.} F.F.Knudsen, The projectivity of moduli spaces of
 stable curves, II: the stacks $M_{g,n}$, {\it Math. Scand},
 {\bf 52} (1983), 161-199.

\item{38.} I.Macdonald, Symmetric functions and Hall polynomials,
 Clarendon Press, Oxford, 1979.

\item{39.} P.A.Macmahon,  Memoir on symmetric functions of
 the roots of systems of equations, {\it Phil. Trans.}, {\bf 181}
 (1890), 481-536 (Collected Papers, vol.2, p. 32-84, MIT Press, 1986).

\item{40.} R.D.MacPherson, The combinatorial formula of Gabrielov,
 Gelfand and Losik for the first Pontrjagin class,
 {\it Sem. Bourbaki}, No.497, 1976-77.

\item{41.} K.Matsumoto, T.Sasaki and M.Yoshida, The monodromy
 of the period map of a 4-parameter family of K3 surfaces and
 the hypergeometric function of type $G(3,6)$, - {\it Intern. J. Math.}, to
appear.

\item{42.} D.Mumford, J.Fogarty,  Geometric Invariant Theory,
Springer- Verlag, New York, 1982.

\item{43.} M.Nagata, On the normality of the Chow variety of
 positive 0-cycles in an algebraic variety, {\it Mem. College.
 Sci. Univ. Kyoto}, {\bf 29} (1955), 165-176.

\item{44.} A.Neeman, Zero-cycles in $P^n$, {\it Adv. in Math.},
 {\bf 89} (1991), 217-227.

\item{45.} K.Okonek, M.Schneider, H.Spindler, Vector bundles
on complex projective spaces (Progress in Math., No.3),
 Birkh\"auser-Verlag, Boston, 1980.

\item{46.} T.G.Room, The geometry of determinantal loci,
Cambridge Univ. Press, 1938.

\item{47.} E.Sernesi, Topics on families of projective schemes,
Queen's papers in Pure and Applied Math., No.73, 1986.

\item{48.} B.Sturmfels, Gr\"obner bases of toric varieties,
 {\it Tohoku Math. J.}, {\bf 43} (1991), 249-261.

\item {49.} T.Oda. Convex bodies and algebraic geometry,
Springer- Verlag, New York, 1988.

\item{50.} B.L van der Waerden, Einf\"uhrung in
die algebraische Geometrie, Springer-Verlag, 1955.

\vskip 2cm
{\sl Department of mathematics, Northwestern University,
 Evanston IL 60208}

e-mail: kapranov@math.nwu.edu

\bye